\documentclass[12pt,a4paper, oneside]{book}
\usepackage{etex}
\usepackage{setspace}
\usepackage[british]{babel}
\usepackage{microtype}
\usepackage{amssymb}
\usepackage{appendix}  
\usepackage{amsmath}
\usepackage{graphs} 
\usepackage{times}
\usepackage{anysize} 
\usepackage[Lenny]{fncychap}
\usepackage{lettrine}
\usepackage{fancyhdr}
\usepackage{url}
\usepackage{bigstrut,array,float, fancyvrb,marvosym,verbatim,combgames}
\marginsize{40mm}{20mm}{20mm}{20mm}

\usepackage{psgo,array,url}

\newenvironment{dedication}
{
   \cleardoublepage
   \thispagestyle{empty}
   \vspace*{\stretch{1}}
   \hfill\begin{minipage}[t]{\textwidth}
   \raggedright
}%
{
   \end{minipage}
   \vspace*{\stretch{3}}
   \clearpage
}

\makeatletter
{\catcode`\|=\z@ \catcode`\\=12 |gdef|bslash{\}}
\makeatother

\setgounit{.4cm}

\newtheorem{definition}{Definition}[section]
\newtheorem{lemma}{Lemma}[section]
\newtheorem{theorem}{Theorem}[section]
\newtheorem{conjecture}{Conjecture}[section]

\newtheorem{question}{Question}[section]

\newcommand{\qed}{\hfill \bf q.e.d.} %Right justified q.e.d
\newcommand{\lef}{\mathcal{L}}
\newcommand{\ri}{\mathcal{R}}
\newcommand{\n}{\mathcal{N}}
\newcommand{\pre}{\mathcal{P}}
\newcommand{\ti}{\mathcal{T}} 
\newcommand{\gfr}{G_F^{SR}}
\newcommand{\gfl}{G_F^{SL}}
\newcommand{\plus}{+_{\ell}}

\newcommand{\G}{\mathcal{G}}
\newenvironment{proof}{{\noindent\it Proof}: }{\qed}

\fancypagestyle{plain}{%
\fancyhf{} % clear all header and footer fields
\fancyhead[C]{\thepage} % except the center

}

\pagestyle{fancy}
\lhead{}
\chead{\thepage} % page number
\rhead{}
\lfoot{}
\cfoot{\itshape \nouppercase \leftmark} % chapter-title
\rfoot{}

\begin{document}
\doublespacing
%\singlespacing
\frontmatter
\bibliographystyle{plain}
\thispagestyle{empty}
\vspace*{6cm}
\begin{center}
\Huge{
SCORING PLAY COMBINATORIAL GAMES}
\end{center}
\vspace{1cm}
\begin{center}
\Large{
Fraser Ian Dowall Stewart}
\end{center}
\vfill
\begin{flushright}
\large{
Submitted for the Degree of Doctor of Philosophy\\
University of Dundee\\
October 2011}
\end{flushright}

\tableofcontents
\listoffigures
 \listoftables
\begin{dedication}
\huge{\textit{For everyone, everywhere, with autism, all of you make our lives more interesting}}
\end{dedication}
\chapter{Acknowledgements}

This thesis could not have been completed without the help and support of my parents.  For the last 10 years they have helped me achieve my dream of becoming a doctor and learning to speak Chinese fluently.  They helped me when I battled depression, and supported me financially after I was rejected from Dalhousie University to study combinatorial game theory.  Thank you Mum and Dad I couldn't have done it without you.

I would also like to thank my supervisor Dr. Keith Edwards.  He accepted me as his PhD student under the understanding that I would be researching harmonious chromatic colourings of graphs.  However when it wasn't working out he agreed to let me research combinatorial game theory instead.  He has been incredibly kind and patient with me, I feel very glad to have had him as my supervisor and I could have made things much easier for him than I did.  He also agreed to take me back as his student when I stopped my PhD for 2 years and effectively terminated my studies.  Without him this would not have been possible.  Thank you, you have my eternal gratitude.

\chapter{Declarations}

\vspace{3cm}

\section*{Candidate's Declaration}

I, Fraser Stewart, hereby declare that I am the author of this thesis; that I have consulted all references cited; that I have done all the work recorded by this thesis; and that it has not been previously accepted for a degree.

\vspace{3cm}

\section*{Supervisor's Declaration}

I, Keith Edwards, hereby declare that I am the supervisor of the candidate, and that the conditions of the relevant Ordinance and Regulations have been fulfilled.
\chapter{Abstract}

This thesis will be discussing scoring play combinatorial games and looking at the general structure of these games under different operators.  I will also be looking at the Sprague-Grundy values for scoring play impartial games, and demonstrating that there is an easily computable function that will solve a large range of octal games easily.  I will also be demonstrating that my theory can readily be applied to the scoring play game of Go and can lead to a much greater understanding of the game.
\mainmatter
\chapter{Introduction}

\section{Games and Outcome Classes}

Combinatorial game theory is the study of all two player, perfect information games, i.e. games where both players have complete information about the game.  For example poker is not a perfect information game as you do not know what cards are in your opponent's hand.  

Combinatorial game theory is the development of mathematical methods that can be used to find winning strategies in perfect information games. Intuitively a combinatorial game is one that we would like to have the following properties, this list is taken from the extras section of chapter 1 of Winning Ways \cite{WW}.

\begin{enumerate}
\item{The game has two players, often called Left and Right.}
\item{There are finitely many positions, and often a particular starting position.  The game also ends after finitely many moves.}
\item{There are clearly defined rules that specify the moves that either player can make from a given position to its options.}
\item{Left and Right move alternately, in the game as a whole.}
\item{Both player have complete information about the game.}
\item{Under normal play the last player to move wins.  Under mis\`ere play the last player to move loses.}
\item{There are no chance moves.}
\item{The rules of play are such that the game will always end because one of the player will be unable to move.  There are no games that are drawn by repetition of moves.}
\end{enumerate}

An example of a game that fully satisfies these axioms is the game Domineering, which has the following rules;

\begin{enumerate}
\item{The game is played on a checkerboard, with a pre-determined shape.}
\item{Each player takes it in turns to place tiles on the board that cover up exactly two squares.}
\item{Left places his tiles vertically, Right places his tiles horizontally.}
\item{The last player to move wins.}
\end{enumerate}

This game satisfies every axiom because; there are two players; at every point in the game each player has only finitely many possible moves (i.e. there are finitely many positions), and the game ends after finitely many moves; the rules are clearly defined; Left and Right move alternately; both players have complete information about the game; the last player to move wins; there is no repetition of moves; finally there are no chance moves. 

Mathematically a combinatorial game is defined as follows;

\begin{definition}\cite{LIP}  The options of a game $G$ are the sets of all games that a player can move to from the game $G$ and are denoted by:

\noindent $G^L=\{\hbox{All games that Left can move to from } G\}$,\\
\noindent $G^R=\{\hbox{All games that Right can move to from } G\}$

A game $G$ is written as $\{G^L|G^R\}$, where $G^L$ and $G^R$ are the options of Left and Right respectively.
\end{definition}

We abuse notation by letting $G^L$ and $G^R$ represent a set of options the the specific options themselves.  Note that when we write down a game we do not write out the braces for the sets $G^L$ and $G^R$, i.e. if $G^L=\{A,B,C,\dots\}$ and $G^R=\{D,E,F,\dots\}$, then $G=\{A,B,C,\dots|D,E,F,\dots\}$.  Also note that if $G^L$ or $G^R$ is empty then we write the game as either $\{.|G^R\}$ or $\{G^L|.\}$.  In other words we do not explicitly write down $\emptyset$, but use a ``dot'' to indicate that a player cannot move.

Examples of well known games that are studied by combinatorial game theorists are \emph{Chess} and \emph{Go}.  \emph{Chess} does not fully satisfy the definition because draws can be permitted, either by infinite play or a stalemate, which do not satisfy axioms 6 and 8.  \emph{Go} also does not fully satisfy this definition because under certain rules repetition of moves is permitted, which also does not satisfy axiom 8, and the winner is determined by the score not who moves last, for further discussion of Go see Section~\ref{Go}.  The theory demonstrated in this chapter will only cover games that fully satisfy the definition. 

\subsection{Outcome Classes}

Outcome classes indicate who wins each game, under optimal play, i.e. if a player has a winning move from a given position then we say that the player wins moving first.  This does not mean that they cannot lose playing first, if they make a bad move.  The standard symbols used to represent the four possible outcome classes are $\lef,\ri,\pre$ and $\n$, which stand for ``Left player win'', ``Right player win'', ``Previous (second) player win'' and ``Next (first) player win'' respectively.  They are defined as follows:

\begin{definition}\cite{LIP}  We define the following;
\begin{itemize}
\item{$\lef=\{G|\hbox{Left wins playing first or second in } G\}$.}
\item{$\ri=\{G|\hbox{Right wins playing first or second in } G\}$.}
\item{$\pre=\{G|\hbox{The second player to move wins in } G\}$.}
\item{$\n=\{G|\hbox{The first player to move wins in } G\}$.}
\end{itemize}
\end{definition}

We will sometimes refer to a position that is a Left win as an $\lef$ position, and similarly for $\ri, \pre$ and $\n$ positions.

\begin{definition}\cite{LIP}
The game tree of a game $G=\{G^L|G^R\}$ is a tree with a root node, and every node has children either on the Left or the Right, and are the Left and Right options of $G$.
\end{definition}

\begin{definition}\cite{LIP} The depth of a game tree, is the length of its longest descending path.
\end{definition}

The only game with depth 0 is $\{.|.\}$.  The games with depth 1 are $\{\{.|.\}|.\}$, $\{.|\{.|.\}\}$ and $\{\{.|.\}|\{.|.\}\}$.  So all games where a player has an option to a game of depth 1, will have depth 2 and so on.  The precise number of games of depth $n$ is currently unknown and is an open problem in Combinatorial Game Theory.

\begin{definition}\cite{LIP}
A short game, is a game that has a game tree of bounded depth.
\end{definition}

\begin{lemma}\label{short}\cite{LIP}
Let $G$ be a short game and let $H$ be an option of $G$, then $\hbox{depth}(H)<\hbox{depth}(G)$.
\end{lemma}

\begin{proof}
Since $G$ is a short game, the length of its longest path is finite, then this implies that the length of the longest path of $H$ much also be finite.  Since $H$ is an option of $G$, then this means that the longest path of $G$ must be at least one greater than the longest path of $H$.
\end{proof}

The reason that we need this result, is because normally when doing a proof by induction, the induction is based on the length of the game trees under consideration.

Normal play theory can be extended to include games that are not short; however for the purposes of this text we will not be considering these games.

\begin{theorem}\label{n}\cite{LIP}
\begin{enumerate}
\item{Under normal play rules a game is a Left win, if and only if Left has an option to either an $\lef$ or a $\pre$ position, and Right has no options at all or only has options to $\n$ or $\lef$ positions.}
\item{Under normal play rules a game is a Right win, if and only if Right has an option to either an $\ri$ or a $\pre$ position, and Left has no options at all or only has options to $\n$ or $\ri$ positions.}
\item{Under normal play rules a game is a $\pre$ position if and only if Left either has no options or only has options to $\n$ or $\ri$ positions and Right either has no options or only has options to $\n$ or $\lef$ positions.}
\item{Under normal play rules a game is an $\n$ position if and only if Left can move to a $\pre$ or an $\lef$ position and Right can move to a $\pre$ or an $\ri$ position.}
\end{enumerate}
\end{theorem}

\begin{proof}
The proof of this is by induction on the depth of game trees.  The case for the induction will be games of depth 0, there is only one, that is $\{.|.\}$, which, trivially, is a $\pre$ position.

Consider the game $G=\{G^L|G^R\}$.  By Lemma~\ref{short} the depth of the game $G$ is greater than that of $G^L$ and $G^R$.  So let $G$ be a game of depth $n+1$ and assume that the theorem holds for all games up to depth $n$.  The following statements are clear:

Left can win going first if and only if there is a $G^L\in \lef\cup\pre$.

Right can win going first if and only if there is a $G^R\in \ri\cup\pre$.

Left can win going second if and only if $G^R=\emptyset$ or for all $g^R\in G^R$, $g^R\in \lef\cup\n$.

Right can win going second if and only if $G^L=\emptyset$ or for all $g^L\in G^L$, $g^L\in \ri\cup\n$.

Therefore, $G\in \lef\Leftrightarrow\exists G^L\in\lef\cup\pre$ and $G^R=\emptyset$ or $G^R\in \lef\cup\n$.

$G\in\ri\Leftrightarrow\exists G^R\in \ri\cup\pre$ and $G^L=\emptyset$ or $G^L\in\ri\cup\n$.

$G\in\n\Leftrightarrow\exists G^L\in\lef\cup\pre$ and $G^R\in \ri\cup\pre$,

$G\in\pre\Leftrightarrow\forall G^L$, $G^L\in \ri\cup\n$ and $\forall G^R$, $G^R\in\lef\cup\n$.

This completes the proof.

\end{proof}

\newpage

\begin{theorem}\label{m}\cite{LIP}
\begin{enumerate}
\item{Under mis\`ere play rules a game is a Left win, if and only if Left has no options to either an $\lef$ or a $\pre$ position, and Right has at least one option and these are to $\n$ or $\lef$ positions.}
\item{Under mis\`ere play rules a game is a Right win, if and only if Right has no options or has an option to either an $\ri$ or a $\pre$ position, and Left has at least one option and these are to to $\n$ or $\ri$ positions.}
\item{Under mis\`ere play rules a game is a $\pre$ position if and only if Left only has options to $\n$ or $\ri$ positions and Right only has options to $\n$ or $\lef$ positions.}
\item{Under mis\`ere play rules a game is an $\n$ position if and only if Left has no options or has an option to a $\pre$ or an $\lef$ position and Right has no options or has an option to a $\pre$ or an $\ri$ position.}
\end{enumerate}
\end{theorem}

\begin{proof}
The proof of this is by induction on the depth of game trees.  The case for the induction will be games of depth 0, there is only one, that is $\{.|.\}$, which, trivially, is a $\n$ position.

Consider the game $G=\{G^L|G^R\}$.  By Lemma~\ref{short} the depth of the game $G$ is greater than that of $G^L$ and $G^R$.  So let $G$ be a game of depth $n+1$ and assume that the theorem holds for all games up to depth $n$.  The following statements are clear:

Left can win going first if and only if $G^L=\emptyset$ or there is a $G^L\in \lef\cup\pre$.

Right can win going first if and only if $G^R=\emptyset$ or there is a $G^R\in \ri\cup\pre$.

Left can win going second if and only if there is at least one Right option and all $G^R\in \lef\cup\n$.

Right can win going second if and only if there is at least one Left option and all $G^L\in \ri\cup\n$.

Therefore, $G\in \lef\Leftrightarrow\exists G^L\in\lef\cup\pre$ or $G^L=\emptyset$ and $G^R\neq\emptyset$ $G^R\in \lef\cup\n$.

$G\in\ri\Leftrightarrow\exists G^R\in \ri\cup\pre$ or $G^R=\emptyset$ and $G^L\neq\emptyset$ $G^L\in\ri\cup\n$.

$G\in \n\Leftrightarrow G^L=\emptyset$ or $G^R=\emptyset$ or $\exists G^L\in \lef\cup\pre$ or $G^R\in\ri\cup\pre$.

$G\in\pre\Leftrightarrow\forall G^L$, $G^L\in\ri\cup\n$ and $\forall G^R$, $G^R\in\lef\cup\n$.

This completes the proof.

\end{proof}

\begin{theorem}\cite{LIP}
All short Combinatorial Games under both normal play and mis\`ere play rules are in one of the four outcome classes, namely:
\begin{enumerate}
\item{Left win}.
\item{Right win}.
\item{Next player win}.
\item{Previous player win}.
\end{enumerate}
\end{theorem}

\begin{proof}
The proof of this is by induction on the depth of game trees.  The four games of depth 1 or less are $\{.|.\}$, $\{\{.|.\}|.\}$, $\{.|\{.|.\}\}$ and $\{\{.|.\}|\{.|.\}\}$, all of which are in one of the four outcome classes.  Assume that the theorem holds for games up to depth $n$.

Let the game $G=\{G^L|G^R\}$ have depth $n+1$.  By Lemma~\ref{short} we know that the games $G^L$ and $G^R$ are of depth $n$ or less, if they exist.  By the inductive hypothesis we know that all of these games in one of the four outcome classes and by Theorems~\ref{n} and \ref{m}  we know that the game $G$ must be in one of the four outcome classes.

This finishes the proof.
\end{proof}

\subsection{Hackenbush}\label{Hackenbush}

The standard example of a Combinatorial Game is a game called Red-Blue Hackenbush.  This game is used extensively throughout the book Winning Ways \cite{WW} to describe all parts of the theory.  The game is played on a graph, with edges coloured red and blue, that are connected to a ``ground'' that is defined arbitrarily before the game begins.  The rules are as follows:

\begin{enumerate}
\item{The players take it in turns to remove edges}.
\item{Left removes Blue edges, Right removes Red edges}.
\item{Any edges not connected to the ground are also removed}.
\item{Under normal play the last player to move wins, under mis\`ere play the last player to move loses}.
\end{enumerate}

\begin{figure}[htb]
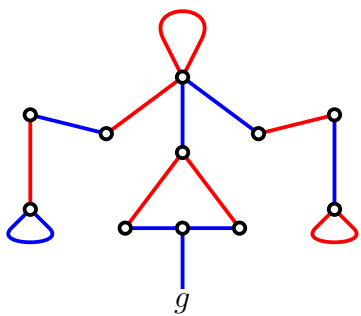

\begin{center}
\begin{graph}(4,4)

\graphnodecolour{1}\graphlinewidth{0.05}
\roundnode{a}(2,0)\roundnode{b}(2,1)\roundnode{c}(1.25,1)\roundnode{d}(2.75,1)
\roundnode{e}(2,2)\roundnode{f}(2,3)\roundnode{g}(1,2.25)\roundnode{h}(0,2.5)
\roundnode{i}(0,1.25)\roundnode{j}(3,2.25)\roundnode{k}(4,2.5)\roundnode{l}(4,1.25)

\edge{a}{b}[\graphlinecolour(0,0,1)]\edge{b}{c}[\graphlinecolour(0,0,1)]
 \edge{b}{d}[\graphlinecolour(0,0,1)] \edge{c}{e}[\graphlinecolour(1,0,0)]
\edge{d}{e}[\graphlinecolour(1,0,0)]
\edge{e}{f}[\graphlinecolour(0,0,1)]\edge{f}{g}[\graphlinecolour(1,0,0)]
\edge{g}{h}[\graphlinecolour(0,0,1)]\edge{h}{i}[\graphlinecolour(1,0,0)]
\edge{f}{j}[\graphlinecolour(0,0,1)]\edge{j}{k}[\graphlinecolour(1,0,0)] \edge{k}{l}[\graphlinecolour(0,0,1)]
\loopedge{f}(0.25,0.5)(-0.25,0.5)[\graphlinecolour(1,0,0)]
\loopedge{i}(0.25,-0.25)(-0.25,-0.25)[\graphlinecolour(0,0,1)]
\loopedge{l}(0.25,-0.25)(-0.25,-0.25)[\graphlinecolour(1,0,0)]

\nodetext{a}{$g$}

\end{graph}
\end{center}
\caption{An example of a Red-Blue Hackenbush position}\label{hex1}
\end{figure}

\begin{figure}[htb]
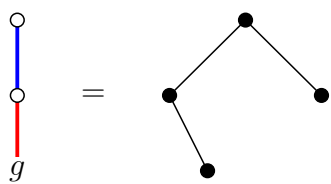

\begin{center}
\begin{graph}(4, 2)

\roundnode{1}(0,0)[\graphnodecolour{1}]\roundnode{2}(0,1)[\graphnodecolour{1}]\roundnode{3}(0,2)[\graphnodecolour{1}]

\edge{1}{2}[\graphlinewidth{0.05}\graphlinecolour(1,0,0)]\edge{2}{3}[\graphlinewidth{0.05}\graphlinecolour(0,0,1)]

\roundnode{4}(2,1)\roundnode{5}(2.5,0)\roundnode{6}(3,2)\roundnode{7}(4,1)

\edge{4}{5}\edge{4}{6}\edge{6}{7}

\freetext(1,1){$=$}\nodetext{1}{$g$}

\end{graph}
\end{center}
\caption{A Red-Blue Hackenbush position with its game tree.}\label{hex3}
\end{figure}

An example a of Red-Blue Hackenbush positions is given in Figures~\ref{hex1}.  An example of a Hackenbush position with its game tree is also given in Figure~\ref{hex3}.  The vertices that are labelled with a ``$g$'' are the vertices that are connected to the ground.

\begin{definition}
A grounded edge, is an edge that is incident with the ground.
\end{definition}

The general strategy of this game is to always remove your highest edge.  This is because it is undesirable to remove grounded edges since the player with the most grounded edges usually has the advantage.  We will be discussing this game in much more detail in chapter 5, however it was shown by Berlekamp in \cite{WW} that in general this game is NP-hard.  In short this means that it is difficult to determine precisely what the best move is or what the winning strategy is, for a further discussion of what NP-hard means please see chapter 5.

\section{Structure of Games in Normal Play}

\subsection{Disjunctive Sum of Games}

It was observed that many games naturally split up into smaller components.  This means you are really playing two or more smaller games at the same time.  This gave rise to the disjunctive sum, and it is defined as follows:

\begin{definition}\cite{LIP}
The disjunctive sum of two games $G$ and $H$ is,

$$G+H=\{G^L+H, G+H^L|G^R+H, G+H^R\}$$
\end{definition}

This is another abuse of notation, where the comma means set union and $G^L+H$ means $H$ added to all elements of $G^L$.  However since we always write down the options that a player has as a list, it makes more sense to write the disjunctive sum in this way.

A good example to illustrate this is from the game of domineering.  Consider the disjunctive sum of a game played on a $2\times 2$ board and a game played on a $1\times 2$ board.  The disjunctive sum of these games is shown in the diagram below.

\begin{center}
$\Domineering{oo\\oo\\}+\Domineering{o\\o\\}=\combgame{\{\Domineering{^o\\vo\\}+\Domineering{o\\o\\},\Domineering{o^\\ov\\}+\Domineering{o\\o\\}, 
\Domineering{oo\\oo\\}+\Domineering{^\\v\\} | \Domineering{<>\\oo\\}+\Domineering{o\\o\\}, \Domineering{oo\\<>\\}+\Domineering{o\\o\\}\}}$
\end{center}

In general this means a player can play in $G$ or $H$ on his turn, but not both.  Also note that if $G=\{.|G^R\}$, i.e. $G^L=\emptyset$ then $G^L+H=\emptyset$, i.e. $G+H=\{G+H^L|G^R+H,G+H^R\}$, which means that if a player can't move in one component, then he must move in one of the other components.  Under these definitions of outcome classes and addition, it is possible to construct the following addition table, which can be found in \cite{LIP}.

\begin{table}[htb]
\begin{center}
\begin{tabular}{c|cccc}
$G+H$&$G\in\pre$&$G\in\n$&$G\in\lef$&$G\in\ri$\\\hline
$H\in\pre$&$\pre$&$\n$&$\lef$&$\ri$\\
$H\in\n$&$\n$&$\lef,\ri,\n,\pre$&$\lef,\ri,\n,\pre$&$\lef,\ri,\n,\pre$\\
$H\in\lef$&$\lef$&$\lef,\ri,\n,\pre$&$\lef$&$\lef,\ri,\n,\pre$\\
$H\in\ri$&$\ri$&$\lef,\ri,\n,\pre$&$\lef,\ri,\n,\pre$&$\ri$\\
\end{tabular}
\end{center}
\caption{Outcome class of $G+H$ under normal play}
\end{table}
\subsection{Greater than and Equal to}

When playing games and analysing them, we would like to know when one option is better than another.  The previous definitions and theorems tell us Left can win going first if he has a good move, but they do not help us determine what that move is.

The following definitions can be found in \cite{LIP}.

\begin{definition}\cite{LIP} We define the following:
\begin{itemize}
\item{$G=H$ if $G+X$ has the same outcome as $H+X$ for all games $X$.}
\item{$G\geq H$ if Left wins $H+X$ implies that Left wins $G+X$ for all games $X$.}
\item{$G\leq H$ if Right wins $H+X$ implies that Right wins $G+X$ for all games $X$.}
\item{$-G=\{-G^R|-G^L\}$.}
\end{itemize}
\end{definition}

What $-G$ means is that we play the same position as $G$, only now Left is playing with Right's options from $G$ and Right is playing with Left's options from G.  A good example is to consider any Hackenbush position.  The ``negative'' of a Hackenbush position is simply to change blue edges for red edges and vice-versa.

The definitions found in On Numbers and Games \cite{ONAG} are stated differently, but the two definitions are equivalent.  While these definitions are the same for mis\`ere play, as we will show in Section~\ref{Mis}, the definitions of equality and greater than do not give a partial order on mis\`ere games, and that while the definition of $-G$ is the same for both mis\`ere play this \emph{does not} correspond to an inverse.

\begin{definition}\cite{LIP}
Two games $G$ and $H$ are equivalent if $G=H$.
\end{definition}

This definition allows to have well defined equivalence classes with a unique representative, depending on the operator being used.

\subsection{Domination and Reversibility}

If we have a game $G=\{A,B,C,\dots|D,E,F,\dots\}$, then we would like to know what the best options for both players are.  That is, no player will want to move to an inferior option, i.e. Left will not choose $A$ over $B$ if $B\geq A$ in the partial order.

\begin{definition}\cite{LIP}
For a game $G=\{A,B,C,\dots|D,E,F,\dots\}$, $A$ dominate B if $A\geq B$, and $D$ dominates $E$ if $D\leq E$.

$A$ is reversible if $A^R\leq G$ and $D$ is reversible if $D^L\geq G$.
\end{definition}

If $A$ is a Left option of $G$ and $A^R\leq G$, then if Left plays to $A$, Right will move to $A^R$ because it is at least as good, if not better than the original game $G$ was.  So Left may as well move to the Left option of $A^R$ directly.  This is what it really means for an option to be reversible.

\begin{definition}\cite{LIP}
A game $G$ is in canonical form if it has no dominated or reversible options.
\end{definition}

\begin{theorem}\cite{LIP}
The canonical form of a game is unique.
\end{theorem}

\begin{proof}  Let $G$ and $H$ be two games such that $G=H$ and neither $G$ nor $H$ has a dominated or reversible option.

So first let Right win $H+X$, since $G=H$, this implies that Right wins $G+X$.  However if Left moves to $G^L+X$ then Right cannot win $G^{LR}+X$, since if he did, then this would mean Right wins $H+X$ implies that Right wins $G^{LR}+X$, i.e. $G\leq H$ and $G$ would have a reversible option, so therefore Right wins $G^L+X^R$.  This implies that Left does not win $H^L+X^R$, since if he did then $H$ would have a dominated option.  Therefore Right wins $G^L+X^R$ if and only if Right wins $H^L+X^R$, i.e. for all $g^L\in G^L$ there is an $h^L\in H^L$ such that $g^L\leq h^L$, and for all $h^L\in H^L$ there is a $g^{L'}\in G^L$ such that $h^L\leq g^{L'}$.

So that means $g^L\leq h^L\leq g^{L'}$, however $g^L$ and $g^{L'}$ must be identical, otherwise $g^L$ is a dominated option.  So, every Left option of $G$ is equal to a Left option of $H$, i.e. $G^L\subseteq H^L$, and by a symmetrical argument $H^L\subseteq G^L$, i.e. $G^L=H^L$, and similarly $H^R=G^R$.  Therefore $G$ and $H$ are identical and the proof is finished. 
\end{proof}

This canonical form is the representative taken for the equivalence class in which it lies.

\subsection{Values}

In normal play we can assign values to certain games, and these values tell us how much of an advantage a player has.  For example is a game has value 1, then we would say that Left has a 1 move advantage.  Sometimes it is not possible to give games numerical values at all.  When we can and when we can not has been well defined, and it is as follows:

\begin{definition}\cite{ONAG}
A game $X=\{X^L|X^R\}$ is a number if and only if $\forall x^L\in X^L$ and $\forall x^R\in X^R$, $x^L<x^R$.
\end{definition}

So if we have a game that has this property, how do we know what its actual value is?  The answer is the simplicity rule and it works as follows:

\begin{definition}\cite{WW}
A dyadic rational is any rational number of the form $\frac{j}{2^i}$, where $j\in \mathbb{Z}$ and $i\in \mathbb{Z}^+$.
\end{definition}

The value of any game that is a number is the simplist dyadic rational that lies strictly between the options of $G$, i.e. the first dyadic rational with the smallest $i$.

For example the value of the games $\{0|1\}=\frac{1}{2}$, $\{\frac{1}{2}|1\}=\frac{3}{4}$, $\{0|2\}=\{0|3\}=1$.

There are also several identities that let us know what the value of a game is.  They are as follows:

\begin{eqnarray}
\{n|.\}&=&n+1, \hbox{if } n\geq0\\
\{n|n+1\}&=&=n+\frac{1}{2}, \hbox{if } n\geq 0\\
\{\frac{p-1}{2^q}|\frac{p+1}{2^q}\}&=&\frac{p}{2^q}, \hbox{if } p \hbox{ is odd}
\end{eqnarray}

We call the game $\{.|.\}=0$.  Examples of games that are not numbers are $\{0|0\}=*$, we call this game ''star``, $\{*|0\}=\downarrow$ and $\{0|*\}=\uparrow$, which are called ``down'' and ``up'' respectively.  We give the value 0 to all games that are $\pre$ positions, positive values to Left win games, and negative values to Right win games.  All other games are ''fuzzy`` or confused with 0; for example $*$ is neither positive, negative or zero.

\begin{theorem}\cite{LIP} Normal play games form a partially ordered abelian group under the disjunctive sum.
\end{theorem}

\begin{proof}  The proof of this will be split into two parts, the first part will demonstrate that normal play games are partially ordered under the disjunctive sum, the second part will show that they form an abelian group.

\noindent\textbf{Part 1: Partial Order}

To show that we have a partially ordered set we need 3 things.

\begin{enumerate}
\item{Transitivity: If $G\geq H$ and $H\geq J$ then $G\geq J$.}
\item{Reflexivity: For all games $G$, $G\geq G$.}
\item{Antisymmetry: If $G\geq H$ and $H\geq G$ then $G=H$.}
\end{enumerate}

1. Let $G\geq H$ and $H\geq J$. $G\geq H$ means that if Left wins $H+X$ this implies that Left wins $G+X$ for all games $X$.  $H\geq J$ means that if Left wins $J+X$ this implies that Left wins $H+X$ for all games $X$.  Since $G\geq H$, then this implies that Left wins $G+X$  for all games $X$, therefore if Left wins $J+X$ this implies that Left wins $G+X$ for all games $X$ and $G\geq J$.

2. Clearly $G\geq G$, since if Left wins $G+X$ then Left wins $G+X$ for all games $X$.

3. First let $G\geq H$ and $H\geq G$.  Since $G\geq H$ then if Left wins $H+X$ this implies that Left wins $G+X$ for all games $X$, however if Left wins $G+X$ then this implies that Left wins $H+X$ for all games $X$.  In other words Left wins $G+X$ if and only if Left wins $H+X$ for all games $X$, and by symmetry Right wins $G+X$ if and only if Right wins $H+X$ for all games $X$.  Therefore $G+X$ has the same outcome as $H+X$ for all games $X$ i.e. $G=H$.

\noindent\textbf{Part 2: Abelian Group}

To show that the set forms an abelian group under the disjunctive sum, we need to show that we have closure, associativity, commutativity, identity and inverses.  So we will deal with each of those properties separately.

\textit{Closure}:  $G+H=\{G^L+H, G+H^L|G^R+H, G+H^R\}$ is a pair of sets of games, by induction, and so $G+H$ is a game.

\textit{Associativity}: Since the Right options follow by symmetry we will focus on the Left options only:

\begin{eqnarray*}
\{(G+H)+J\}^L\\
&=&\{(G+H)^L+J, (G+H)+J^L\}\\
&=&\{(G^L+H)+J, (G+H^L)+J, (G+H)+J^L\}\hbox{ (by induction)}\\
&=&\{G+(H+J)\}^L
\end{eqnarray*}

\textit{Commutativity}: By definition, $G+H=\{G^L+H, G+H^L|G^R+H, G+H^R\}$ and by induction, all the simpler games are commutative, for example, $g^L+H=H+g^L$ for all $g^L\in G^L$.  So,

\begin{eqnarray*}
G+H&=&\{G^L+H, G+H^L|G^R+H, G+H^R\}\\
&=&\{H+G^L, H^L+G|H+G^R, H^R+G\} \hbox{ (by induction)}\\
&=&H+G
\end{eqnarray*}

\textit{Identity}: We take the identity to be any game $I$ such that $G+I$ has the same outcome as $G$ for all games $G$.  In particular for all $I\in \pre$, $G+I$ has the same outcome as $G$ for all games $I$.

There are two cases to consider since the final two cases follow by symmetry, either Left wins moving first on $G$ or Right wins moving second on $G$.

Case 1: Left wins moving first on $G$.  So consider $G+I$ where $I\in \pre$.  Left moving first on $G+I$ will play his winning move on $G$ and move to $G^L+I$, since $I\in \pre$ then if Right has no move on $I$ he must also play on $G$, and since Left wins moving first on $G$, this implies that Left will move last on $G$ meaning that Right must play first on $I$, and will lose (since he has no move on $I$).  If Right does have a move on $I$, then Left will respond in $I$, since $I\in \pre$ this implies that Left moves last on $I$, and therefore Right must eventually move first on $G^L$ and will lose.  Therefore if Left wins moving first on $G$, this implies that Left wins moving first on $G+I$, for all $I\in \pre$.

Case 2: Right wins moving second on $G$.  Consider $G+I$, since $I\in\pre$ this implies that Right moves most last playing second on both $I$ and $G$, therefore which ever game Left chooses to move on, Right will always be able to respond on the same game, which guarantees that he will move last on $G+I$, and therefore win, moving second.

Therefore $G+I$ has the same outcome as $G$ for all games $G$ and for all $I\in\pre$.

\textit{Inverses}:  For all games $G$ we take the game $-G$ to be its inverse, i.e. $G+(-G)\in\pre$ for all games $G$.

Since $-G=\{-G^R|-G^L\}$, this means that Right's options on $G$ and identical to Left's options on $-G$ and vice-versa.  Therefore whichever component the first player chooses on $G+(-G)$ the second player can play the exact same move on the opposite component, e.g. if Left moves to $G^L+(-G)$, Right can move to $G^L+(-G^L)$ and so on.  So by playing this strategy the second player can always guarantee that they will move last, and therefore win $G+(-G)$, therefore $G+(-G)\in \pre$.

This completes the proof.
\end{proof}

\section{Impartial Games}\label{impartial}

\begin{definition}\cite{WW}
An impartial game is one where both players have the same options, and the options are also impartial.
\end{definition}

The theory of impartial games has been well established in normal play, and has had similar developments made in mis\`ere play due to the work of Plambeck and Siegel.  However I will not be discussing the theory for mis\`ere play here, instead I will only be describing the theory in normal play.

\subsection{Nim}\label{nim}

Nim is the best example of an impartial game, it has been studied in many different variants by combinatorial game theorists \cite{PG, PS, RS}.  The reason that Nim is an important game to study is because, as I will show in Theorem~\ref{impt}, under normal play, for every impartial game $G$, there is a Nim heap of size $n$ such that $G=n$.  The rules of Nim are as follows:

\begin{enumerate}
\item{The players play with heaps of beans.}
\item{On a player's turn he may remove as many beans as he wants from any one heap.}
\item{The game ends when there are no beans left.}
\item{Under normal play rules the last player to move wins, under mis\`ere play rules the last player to move loses.}
\end{enumerate}

The best way to describe the solution to this game is to look at an example.  Suppose we are playing normal play Nim with heaps of size 6, 7 and 3.  We first write them out base 2, and we get $6=(1,1,0)$, $7=(1,1,1)$ and $3=(0,1,1)$.  Then we write this out in tabular form to get:

\begin{table}[htb]
\begin{center}
\begin{tabular}{c|ccc}
3&0&1&1\\
7&1&1&1\\
6&1&1&0\\\hline
&0&1&0\\
\end{tabular}
\end{center}
\caption{Example of a Nim position}
\end{table}

The winning move is to make the number of 1's in each column even, or so that the sum of the sizes of the remaining heaps add up to $0$ base 2.  The reason for why this is the winning move will be described in the following sections.

In normal play the very first impartial game is, $\{.|.\}=0$.  The next one is $\{0|0\}=*$ and represents a Nim heap of size 1.  After that we have $\{0,*|0,*\}=*2$, which represents a Nim heap of size 2, and in general we have the following formula:

$$\{0,*,\dots,*(n-1)|0,*,\dots,*(n-1)\}=*n$$

\begin{definition}\cite{WW}
$\hbox{mex}\{a,b,c,\dots\}=m$, where $m$ is the least non-negative integer that does not appear in the set $\{a,b,c,\dots\}$.
\end{definition}

For example, $\hbox{mex}\{2,3,10,132260\}=0$, $\hbox{mex}\{0,1,2,9,\dots\}=3$.

We then get an evaluation scheme for a general position, in normal play Nim.  It is expressed in the following theorem:

\begin{theorem}\cite{LIP}\label{mex}
If $G=\{*a,*b,*c,\dots|*a,*b,*c,\dots\}$, then the value of $G$ is $*m$, where $m=\hbox{mex}\{a,b,c,\dots\}$.
\end{theorem}

\begin{proof} To prove this we will show that $G-*m\in \pre$.

If either player moves either component to $*k$ for $k<n$, there is a matching move in the other component.  In particular, since the mex of $\{a,b,c,\dots\}$ is $m$, this means that $*k$ is an option of $G$ and $*m$.  Hence, the second player can respond to $*k-*k\in \pre$.

The only other moves are from $G-*m$ to $*k-*m$ for $k>m$.  In this case $*m$ is an option from $*k$, so the second player responds locally to $*n-*n\in\pre$.

\end{proof}

Under this evaluation scheme for normal play games, other games may have the same value as that of a Nim heap of size $m$, but we will always refer to it as being $*m$, i.e. $*m$ is the value of the game.  The values $*m$ are called ``Nimbers'', after the game Nim where they came from, or Sprague-Grundy values.

\begin{theorem}\cite{LIP} An impartial game $G=\{*a,*b,*c,\dots|*a,*b,*c,\dots\}$ is a $\pre$ position if and only if $\hbox{mex}\{a,b,c,\dots\}=0$; it is an $\n$ position otherwise.
\end{theorem}

\begin{proof}  If $G=\{*a,*b,*c,\dots|*a,*b,*c,\dots\}$ is an impartial game, then $G=*m$, where $m=\hbox{mex}\{a,b,c,\dots\}$, since $*0=0=\{.|.\}\in\pre$, if $m=0$ then this implies that $G\in\pre$.  If $m\neq 0$ then this implies that $0\in\{a,b,c,\dots\}$, i.e. both players have an option to move to the game $0$, which implies that $G\in\n$.  Therefore the theorem is proven.
\end{proof}

Using this theory it is possible to determine when the disjunctive sum of two impartial games is a $\pre$ position or an $\n$ position.  This uses something called ``Nim addition''.  First you write the Sprague-Grundy values of the two games base 2.  Then you add them together using the ``Exclusive Or'' function.

The general strategy for impartial games is encapsulated in the following theorem:

\begin{theorem}\cite{LIP}\label{impt}
(a) If $(a,b,c,\dots)$ is a position in Nim, then this is a $\pre$ position if and only if $a\oplus b\oplus c\dots=0$.

(b) Any impartial game $G$ is equivalent to a Nim-heap, by $G=*a$ where,

$$a=\hbox{mex}\{G'|G'\hbox{ are the options of }G\}$$

(c) For non-negative integers $k$ and $j$, $*k+*j=*(k\oplus j)$.
\end{theorem}

\begin{proof}  We will prove each part of this separately.

(a) Suppose that $a\oplus b\oplus\dots = 0$.  Without loss of generality, suppose that the first player removes $r$ counters from heap $a$.  Since the binary expansion of $a-r$ is not the same as $a$ then $(a-r)\oplus b\oplus\dots\neq 0$.

Now suppose that $q=a\oplus b\oplus\dots\neq 0$.  Let $q_jq_{j-1}\dots q_0$ be the binary expansion of $q$; that is, each bit $q_i$ is either 1 or 0 and $q_j=1$.  Then one of the heaps, again without loss of generality, say $a$, must have a 1 in position $j$ in its binary expansion.

We will show that there is a move from this heap of size $a$ which produces a position whose is 0 base 2; in particular, reducing $a$ to $x=q\oplus a$.

\noindent\textit{Step 1.} The move is legal: 

In changing $a$ to $q\oplus a$ the left most bit in $a$ that is changed is a 1 (to a 0) and therefore, $q\oplus a<a$, and so the move is legal. (The move \emph{reduces} the size of the heap.)

\noindent{Step 2.} Sum is 0 base 2: 

The resulting position is

\begin{eqnarray*}
(q\oplus a)\oplus b\dots &=& ((a\oplus b\dots)\oplus a)\oplus b\oplus\dots\\
&=&(a\oplus a)\oplus(b\oplus b)\oplus\dots\\
&=&0\\
\end{eqnarray*}

(b) Let $G$ be impartial. By induction the options of $G$ are equivalent to Nim-heaps and by Theorem~\ref{mex} we know how to find the equivalent Nim-heap for $G$.

(c) From part (a), we have that Nim with heaps of size $k$, $j$ and $k\oplus j$ is a $\pre$ position.  The values of the individual Nim-heaps are, respectively $*k$, $*j$ and $*(k\oplus j)$.  The fact that they form a $\pre$ position means that $*k+*j=*(k\oplus j)$.
\end{proof}

What these theorems tell us is how to win any impartial game $G$ simply by knowing its options.  That is we can easily evaluate the value of an impartial game $G$, and if it is $0$ then there is no winning move from that position for the first player, and if it is not $0$, then there will always be an option whose value is $0$.

\section{Structure of Games in Mis\`ere Play}\label{Mis}

From the definition of a combinatorial game we know that mis\`ere play means that the last player to move loses.  This small difference changes things dramatically.  When look at game outcome classes for mis\`ere, we do not just get the reverse of normal play, i.e. interchange $\lef$'s for $\ri$'s and $\pre$'s for $\n$'s.

For example, in normal play a Nim heap of size 1 is an $\n$ position and a $\pre$ position in mis\`ere play.  However a Nim heap of size 2 or greater is an $\n$ position in both mis\`ere and normal play.

In the first section I explained that for normal play games there is a nice addition table, and that the games formed a partially ordered abelian group under the disjunctive sum, and that a Left win played under the disjunctive sum with another Left win was always a Left win regardless of what the games were.  For mis\`ere play, Ottaway \cite{O} proved the following theorems;

\begin{theorem}\cite{O} Undere mis\`ere rules, for any three outcome classes $\mathcal{X}$, $\mathcal{Y}$ and $\mathcal{Z}$, there exists a $G\in \mathcal{X}$, $H\in\mathcal{Y}$ such that $G+H\in \mathcal {Z}$.
\end{theorem}

In particular that it is possible to find two Left win games under mis\`ere rules so that when they are played under the disjunctive sum, the resulting game is a Right win.  We also showed that under normal play, the entire set of $\pre$ position games, act as an identity set under the disjunctive sum.  The following theorem shows that for mis\`ere play games, the identity set is trivial and contains only one element, namely $\{.|.\}$.

\begin{theorem}\cite{O} Under mis\`ere rules, there is no, non-trivial, game $X$ such that $G+X$ has the same outcome as $G$, for all games $G$.
\end{theorem}

\section{Normal and Mis\`ere Play Under Different Operators}

In \cite{ONAG}, Conway gives definitions for other operators that can be used for combinatorial games, such as the conjunctive sum, where players move on all components every turn, and the selective sum, where players choose which components to play on on their turn.  They are defined as follows.

\begin{definition}  The conjunctive sum and the selective sum are defined respectively;

\begin{itemize}
\item{$G\bigtriangleup H=\{G^L\bigtriangleup H^L|G^R\bigtriangleup H^R\}$}
\item{$G\triangledown H=\{G^L\triangledown H, G\triangledown H^L, G^L\triangledown H^L|G^R\triangledown H, G\triangledown H^R, G^R\triangledown H^R\}$}
\end{itemize}
\end{definition}

It is important though to look for other natural ways of playing combinatorial games, since under normal play an entire number system was discovered \cite{ONAG} by examining the structure of normal play games under the disjunctive sum.  So while mis\`ere and scoring play games may not yield any apparent structure under the disjunctive sum, we may yet be surprised to learn what the structure is if we were to use a different operator.

\begin{definition}\cite{FS2}  The sequential join of two games $G$ and $H$ is
defined as follows:
$$G\rhd H=\begin{cases}

&\{G^L\rhd H| G^R\rhd H\} \text{, if $G\neq \{|\}$}\\

&H \text{, Otherwise}

\end{cases}$$

\end{definition}

In \cite{FS2}, Stewart proved the following two theorems about the sequential join.

\begin{theorem}\label{SJN} The addition table for the sequential join under normal play rules is given in Table~\ref{tn}.

\begin{table}[htb]
\begin{center}
\begin{tabular}{c|cccc}
$G\rhd H$&$G\in \pre$&$G\in \n$&$G\in \lef$&$G\in \ri$\\\hline
$H\in \pre$&$\pre$&$\n$&$\lef$&$\ri$\\
$H\in \n$&$\lef,\ri,\n,\pre$&$\lef,\ri,\n,\pre$&$\lef,\ri,\n,\pre$&$\lef,\ri,\n,\pre$\\
$H\in \lef$&$\lef,\pre$&$\lef,\n$&$\lef$&$\lef,\ri,\n,\pre$\\
$H\in \ri$&$\ri,\pre$&$\ri,\n$&$\lef,\ri,\n,\pre$&$\ri$\\
\end{tabular}
\end{center}
\caption{Outcome class table of $G\rhd H$ under normal play}\label{tn}
\end{table}

\end{theorem}

\begin{theorem}\label{SJM} The addition table for the sequential join under mis\`ere rules is given in Table~\ref{tm}.

\begin{table}[htb]
\begin{center}
\begin{tabular}{c|cccc}
$G\rhd H$&$G\in \n$&$G\in \pre$&$G\in \lef$&$G\in \ri$\\\hline
$H\in \n$&$\n$&$\pre$&$\lef$&$\ri$\\
$H\in \pre$&$\lef,\ri,\n,\pre$&$\lef,\ri,\n,\pre$&$\lef,\ri,\n,\pre$&$\lef,\ri,\n,\pre$\\
$H\in \lef$&$\lef,\n$&$\lef,\pre$&$\lef$&$\lef,\ri,\n,\pre$\\
$H\in \ri$&$\ri,\n$&$\ri,\pre$&$\lef,\ri,\n,\pre$&$\ri$\\
\end{tabular}
\end{center}
\caption{Outcome class of $G\rhd H$, under mis\`ere play}\label{tm}
\end{table}

\end{theorem}

Both proofs can be found in either \cite{FS} or \cite{FS2}.  As shown both mis\`ere and normal play games form a non-trivial monoid under the sequential join, with a well defined identity set, and interestingly the structure of both sets of games seems to be remarkably similar.  We will look at this in relation to scoring play games further in chapter 4.

While there is no relationship between the different tables themselves, what these tables do give us is an interesting comparison in the relationship between the four outcome classes under different operators.  

\chapter{Scoring Play Games}

After normal and mis\`ere play games, scoring play games are a third and totally overlooked way of playing combinatorial games.  With scoring play we are no longer interested in who moves last, but who has accumulated the most points during play. 

An example of a scoring play game is the ancient Chinese game of Go.  This game is played on a grid, and players take it in turns to place stones on the intersections.  The winner is the player who captures the largest area of the board, and the most of their opponent's stones, combined.

Another example is an ancient African game called Mancala.  In this game both players have six cups and one kala, a special cup assigned to each player.  There are beans placed in each of the player's cup, and the players take it in turns to sow these beans into the other cups, and their kala, the winner is the player who manages to gather most beans into his or her kala.

This chapter will lay the foundations for the remainder of the thesis.  When we discuss a new set of objects, such as scoring play games, we must make sure the definitions are the correct ones.  This chapter will be looking at the general structure of scoring play games under the disjunctive sum, since it is by far the most common and widely used operator in combinatorial game theory.

Intuitively we want all scoring play games to have the following four properties;

\begin{enumerate}
\item{The rules of the game clearly define what points are and how players either gain or lose them.}
\item{When the game ends the player with the most points wins.}
\item{For any two games $G$ and $H$, $a$ points in $G$ are equal to $a$ points in $H$, where $a\in \mathbb{R}$. For example in the game Go you get one point for each of your opponent's stones that you capture, and for each piece of area you successfully take.  In Mancala you get one point for each bean you place in your Kala, so when comparing these games we would like one point in Mancala to be worth one point in Go.}
\item{At any stage in a game $G$ if Left has $L$ points and Right has $R$ points, then the score of $G$ is $L-R$, where $L,R\in \mathbb{R}$.}
\end{enumerate}

Mathematically we define scoring play games in the following way, using a variant of the Conway \cite{ONAG} definition for a combinatorial game;

\begin{definition} A scoring play game $G=\{G^L|G^S|G^R\}$, where $G^L$ and $G^R$ are sets of games and $G^S\in\mathbb{R}$, the base case for the recursion is any game $G$ where $G^L=G^R=\emptyset$.

\noindent $G^L=\{\hbox{All games that Left can move to from } G\}$\\
\noindent $G^R=\{\hbox{All games that Right can move to from } G\}$,

and for all $G$ there is an $S=(P,Q)$ where $P$ and $Q$ are the number of points that Left and Right have on $G$ respectively.  Then $G^S=P-Q$, and for all $g^L\in G^L$, $g^R\in G^R$, there is a $p^L,p^R\in\mathbb{R}$ such that $g^{LS}=G^S+p^L$ and $g^{RS}=G^S+p^R$.

$\gfl$ and $\gfr$ are called the final scores of $G$ and are the largest scores that Left and Right can achieve when $G$ ends, moving first respectively, if both players play their optimal strategy on $G$.

\end{definition}

\begin{definition}
The game tree of a scoring play game $G=\{G^L|G^S|G^R\}$ is a tree with a root node, and every node has children either on the Left or the Right that are the Left and Right options of $G$.  All nodes are numbered, and are the scores of the game $G$ and all of its options.
\end{definition}

It is important to note that we will only be considering finite games, i.e. for any game $G$, the game tree of $G$ has finite depth and finite width.  This means that $\gfl$ and $\gfr$ are always computable, and cannot be infinite, or unbounded.  

There is also the case where a game may have a form of aggregate scoring.  For example players may play two games in sequence, and the winner would be the player who gets the most points over both games.  This gives scoring play games an additional dynamic, where in the event of a tie after two games, the winner may be determined by the player who managed to accumulate more points in one of the games.

However as far as this thesis is concerned, I will not be considering games of this type.  We will only look at games where the winner is determined after one game ends.  Games with aggregate scoring would be an interesting area to look at for further research.

There are two conventions that I will be using throughout this thesis.  The first is that in all examples given we will take the initial score of the game to be $0$, unless stated otherwise.  The second is that if for a game $G$, $G^L=G^R=\emptyset$, I will simply write $G$ as $G^S$, rather than $\{.|G^S|.\}$.   For example the game $G=\{\{.|0|.\}|1|\{.|2|.\}\}$, will be written as $\{0|1|2\}$.  The game $\{.|n|.\}$, will be written as $n$, and so on.  This is simply for convenience and ease of reading.

\subsection{Games Examples}\label{ch2ex}

Before I continue I will give an example of a scoring play game to demonstrate how to use the above notation.  So consider the game Toad and Frogs from Winning Ways \cite{WW}, played under scoring play rules.  The rules are as follows;

\begin{enumerate}
\item{The game is played on a horizontal one dimensional grid.}
\item{Left moves Toads and Right moves Frogs.}
\item{Toads move from left to right and Frogs move from right to left, and on a player's turn they may only move their piece one space.}
\item{If a Toad is next to a Frog and the following space is empty, then the Toad may "jump" the Frog into the next space, and vice versa.  Toads may only jump Frogs and vice versa.}
\item{The player who makes the most jumps wins.}
\end{enumerate}

In Winning Ways, the rules are that the last player to move wins, by changing the winning condition it changes it to a scoring play game.  Some game positions with their corresponding game values are as follows.  In the Figures~\ref{tf1} and \ref{tf2} $B$ represents a blank square, $T$ represents Toads and $F$ represents Frogs.  The numbers in brackets are the current score.  The size of the grid is equal to the number of letters in the diagram, e.g. $TBF$ is three squares long, $TTBF$ is four squares and so on.

The game in Figure~\ref{tf1} has value $\{\{.|0|\{-1|-1|.\}\}|0|\{\{.|1|1\}|0|.\}\}$, and the game in Figure~\ref{tf2} has value $\{0|0|\{\{.|0|0\}|-1|.\}|-1|.\}|0|\{\{0|0|.\}|1|.\}|1|\{\{.|2|2\}|1|.\}|0|.\}$.  Both of these games are in ``canonical form'', that is neither has a dominated or reversible option.  For more details see Section~\ref{canfor}.

\begin{figure}[htb]
\begin{center}
\begin{pspicture}(7,6)
\put(0,3){$TBF=$}\put(4,6){\cgtree{
  {TBF(0)}(+{BTF(0)}(|{FTB(-1)}({FBT(-1)}|)) |+{TFB(0)}({BFT(1)}(|{FBT(1)})|))
}}
\end{pspicture}
\end{center}
\caption{$TBF=\{\{.|0|\{-1|-1|.\}\}|0|\{\{.|1|1\}|0|.\}\}$}\label{tf1}
\end{figure}

\begin{figure}[htb]
\begin{center}
\begin{pspicture}(8,10)
\put(-1,4){$TTBF=$}\put(6,10){\cgtree{
{TTBF(0)}(++{TBTF(0)}({BTTF(0)}|{TFTB(-1)}({TFBT(-1)}({BFTT(0)}(|{FBTT(0)})|)|))|++{TTFB(0)}({TBFT(1)}(+{BTFT(1)}(|{FTBT(0)}({FBTT(0)}|))|+{TFBT(1)}({BFTT(2)}(|{FBTT(2)})|))|))}}
\end{pspicture}
\end{center}
\caption{$TTBF=\{0|0|\{\{.|0|0\}|-1|.\}|-1|.\}|0|\{\{0|0|.\}|1|.\}|1|\{\{.|2|2\}|1|.\}|0|.\}$}\label{tf2}
\end{figure}

\section{Outcome Classes}

Under scoring play the outcome classes are a little different to the outcome classes in normal and mis\`ere play.  In combinatorial game theory we would like to know who wins under optimal play, e.g. if $G\in\lef$, then that means Left has a winning strategy moving first or second, if he plays his optimal strategy.

In scoring play this is just as important, since there may be instances when computing the final score might be too difficult, but determining the outcome class may be relatively easy.  It is also important, since if we have a game $K$ that is played as the disjunctive sum of two smaller components $G$ and $H$, and we know the outcome of both $G$ and $H$, then can we also say what the outcome of $K$ is.  In normal play this is certainly the case, but as has been shown \cite{O}, it is not the case for mis\`ere play.

\newpage

Before we define what the outcome classes are, we first define the following;

\begin{definition}
\item{$L_>=\{G|\gfl>0\}$ , $L_<=\{G|\gfl<0\}$ , $L_= =\{G|\gfl=0\}$.}
\item{$R_>=\{G|\gfr>0\}$ , $R_<=\{G|\gfr<0\}$ , $R_= =\{G|\gfr=0\}$.}
\item{$L_\geq=L_>\cup L_=$ , $L_\leq = L_<\cup L_=$.}
\item{$R_\geq=R_>\cup R_=$ , $L_\leq = R_<\cup R_=$.}
\end{definition}

So what this means is that $L_>$ is the set of all games that Left wins moving first, $L_<$ is the set of all games that Left loses playing first, $L_=$ is the set of all games that are a tie when Left moves first, $L_\geq$ is the set of all games that Left does not lose, and may win moving first, and so on.

Since we would like to classify every game by an outcome class it is also important that every game belongs to exactly one outcome class.  Again this is true for normal play and mis\`ere play, however we cannot use the same definitions as before, since under scoring play, a final score of $0$ means a tied game, i.e. nobody wins.  Therefore this will give us five outcome classes, and they are defined as follows.

\begin{definition}
The outcome classes of scoring play games are defined as follows:

\begin{itemize}
\item{$\lef=(L_>\cap R_>)\cup(L_>\cap R_=)\cup(L_=\cap R_>)$}
\item{$\ri=(L_<\cap R_<)\cup(L_<\cap R_=)\cup(L_=\cap R_<)$}
\item{$\n=L_>\cap R_<$}
\item{$\pre=L_<\cap R_>$}
\item{$\ti=L_=\cap R_=$}
\end{itemize}
\end{definition}

What these mean is that $\lef$ is the set of games that Left wins moving first and wins or ties moving second, and vice versa, and similarly for $\ri$.  $\n$ is the set of games that the first player to move wins, $\pre$ is the set of games that the second player to move wins, and $\ti$ is the set of games that ends in a tie regardless of who moves first.

The reason that I chose the outcome classes in this way, is because if you have a game $G=\{1|0|0\}$, then it is more natural to say that belongs to the outcome $\lef$, since Right cannot win, but Left can if he moves first.  In this way we also keep the usual convention of calling a game $G\in \n$ a ``next player win'' and  a game $H\in \pre$ a ``previous player win''.

This also creates an interesting distinction in that while $\lef$ means the set of games where Left can win moving first or second in both normal and mis\`ere play, in scoring play, it means that if Left wins moving first he doesn't lose, and may win, moving second, and vice-versa.  Another distinction is the addition of the outcome class $\ti$, which of course does not exist in either normal or mis\`ere play.

To demonstrate the different outcome classes I will give an example, so consider the game scoring play Hackenbush.  The rules of scoring play Hackenbush are as follows;

\begin{enumerate}
\item{The game is played on a graph with coloured edges that is connected to a ground defined arbitrarily before the game begins.}
\item{Players take it in turns to remove edges, any edges that are disconnected from the ground are also removed.}
\item{Left removes blue edges and Right removes red edges.}
\item{The player who disconnects the most edges from the ground wins.}
\end{enumerate}

\begin{figure}[htb]
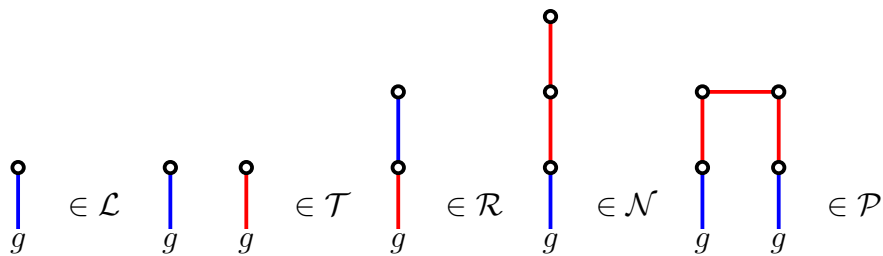

\begin{center}
\begin{graph}(11,3)

\graphnodecolour{1}\graphlinewidth{0.05}
\roundnode{1}(0,0)\roundnode{2}(0,1)\roundnode{3}(2,0)\roundnode{4}(2,1)\roundnode{5}(3,0)\roundnode{6}(3,1)
\roundnode{7}(5,0)\roundnode{8}(5,1)\roundnode{9}(5,2)\roundnode{10}(7,0)\roundnode{12}(7,1)\roundnode{13}(7,2)\roundnode{14}(7,3)\roundnode{18}(9,0)\roundnode{19}(9,1)\roundnode{20}(9,2)\roundnode{21}(10,0)\roundnode{22}(10,1)
\roundnode{23}(10,2)

\edge{1}{2}[\graphlinecolour(0,0,1)]\edge{3}{4}[\graphlinecolour(0,0,1)]\edge{5}{6}[\graphlinecolour(1,0,0)]\edge{7}{8}[\graphlinecolour(1,0,0)]\edge{8}{9}[\graphlinecolour(0,0,1)]\edge{10}{12}[\graphlinecolour(0,0,1)]\edge{12}{13}[\graphlinecolour(1,0,0)]\edge{13}{14}[\graphlinecolour(1,0,0)]\edge{18}{19}[\graphlinecolour(0,0,1)]\edge{19}{20}[\graphlinecolour(1,0,0)]\edge{21}{22}[\graphlinecolour(0,0,1)]\edge{22}{23}[\graphlinecolour(1,0,0)]\edge{20}{23}[\graphlinecolour(1,0,0)]

\freetext(1,0.5){$\in\lef$}\freetext(4,0.5){$\in\ti$}\freetext(6,0.5){$\in\ri$}\freetext(8,0.5){$\in\n$}\freetext(11,0.5){$\in\pre$}

\nodetext{1}{$g$}\nodetext{3}{$g$}\nodetext{5}{$g$}
\nodetext{7}{$g$}\nodetext{10}{$g$}\nodetext{18}{$g$}
\nodetext{21}{$g$}

\end{graph}
\end{center}
\caption{Examples of games in different outcome classes.}\label{exout}
\end{figure}

As with all Hackenbush diagrams, nodes labelled with a ``$g$'' are grounded nodes.  Examples of games in the five different outcome classes are given in Figure~\ref{exout}.

\begin{theorem}
Every game $G$ belongs to exactly one outcome class.
\end{theorem}

\begin{proof}  This is clear since every game belongs to exactly one of $L_>$, $L_<$, $L_=$ and exactly one of $R_>$, $R_<$, $R_=$.  Therefore every game belongs to exactly one of the nine possible intersections of $L_>$, $L_<$, $L_=$ and $R_>$, $R_<$, $R_=$.  Since each outcome class is simply the union of one or more of these then each game can only be in exactly one outcome class.
\end{proof}

\section{The Disjunctive Sum}\label{disjunctive}

As I mentioned earlier, the disjunctive sum is by far the most commonly used operator in combinatorial game theory.  This is because many well known games such as Go naturally break up into the disjunctive sum of two or more components.  For scoring play the disjunctive sum needs to be defined a little differently, this is because in scoring games when we combine them together we have to sum the games and the scores separately.  

For this reason I will be using two symbols $\plus$ and $+$.  The $\ell$ in the subscript stands for ``long rule'', this comes from \cite{ONAG}, and means that the game ends when a player cannot move on all components on his turn.  The ``short rule'' means that the game ends when a player cannot move on at least one component on his turn.  In this thesis every operator that we consider will be played under the ``long rule'' since it is the most common rule to play a combinatorial game with. 

Mathematically what this means is that if $G=G_1*G_2*\dots*G_n$, then $G^L=\emptyset$ if and only if $G^L_i=\emptyset$ for all $i=1,\dots,n$ under the long rule, and $G^L=\emptyset$ if there it at least one $j$ such that $G^L_j=\emptyset$, where $j=1,\dots,n$ under the short rule, for any operator $*$.

I will also be describing what it means for each operator individually, but this difference is crucial when studying combinatorial games since the ending condition is of course very important.

So the definition is as follows;

\begin{definition}  The disjunctive sum is defined as follows:
$$G\plus H=\{G^L\plus H,G\plus H^L|G^S+H^S|G^R\plus H,G\plus H^R\},$$
\noindent where $G^S+H^S$ is the normal addition of two real numbers.
\end{definition}

As with the disjunctive sum of normal and mis\`ere play games we abuse notation by making the comma mean set union, and $G^L\plus H$ means take the disjunctive sum of all $g^L\in G^L$ with $H$.  Note also that if $H^L=\emptyset$ and $G^L\neq \emptyset$ then $G^L\plus H,G\plus H^L=G^L\plus H$ under the long rule and $\emptyset$ under the short rule.

\subsection{Disjunctive Sum Examples}

Now I will give some examples of games played under the disjunctive sum to demonstrate how to compute the disjunctive sum of two games.  The first two examples are from scoring play Hackenbush.

\begin{center}
\begin{graph}(3,1)

\graphnodecolour{1}\graphlinewidth{0.05}

\roundnode{1}(0,0)\roundnode{2}(0,1)\roundnode{3}(3,0)\roundnode{4}(3,1)

\edge{1}{2}[\graphlinecolour(0,0,1)]\edge{3}{4}[\graphlinecolour(0,0,1)]

\freetext(1.5,0.5){$\plus$}

\nodetext{1}{$g$}\nodetext{3}{$g$}

\end{graph}
\end{center}

\noindent\begin{eqnarray*}
&=&\{1|0|.\}\plus\{1|0|.\}\\
&=&\{\{2|1|.\}|0|.\}
\end{eqnarray*}

\begin{center}
\begin{graph}(3,3)

\graphnodecolour{1}\graphlinewidth{0.05}

\roundnode{1}(0,0)\roundnode{2}(0,1)\roundnode{3}(0,2)
\roundnode{4}(3,0)\roundnode{5}(3,1)

\edge{1}{2}[\graphlinecolour(0,0,1)]\edge{2}{3}[\graphlinecolour(1,0,0)]\edge{4}{5}[\graphlinecolour(1,0,0)]

\freetext(1.5,1){$\plus$}

\nodetext{1}{$g$}\nodetext{4}{$g$}

\end{graph}
\end{center}

\noindent\begin{eqnarray*}
&=&\{2|0|\{0|-1|.\}\}\plus\{.|0|-1\}\\
&=&\{2\plus\{.|0|-1\}|0|\{0|-1|.\}\plus\{.|0|-1\},-1\plus\{2|0|\{0|-1|.\}\}\}\\
&=&\{\{.|2|1\}|0|\{\{.|0|-1\}|-1|\{-1|-2|.\}\},\{1|-1|\{-1|-2|.\}\}\}
\end{eqnarray*}

The final example is the disjunctive sum of two games where the initial scores are not zero.  So consider $G=\{3|2|-4\}$ and $H=\{1.5|0.5|1\}$, then $G\plus H=\{3\plus \{1.5|0.5|1\}, 1.5\plus \{3|2|-4\}|2.5|-4\plus \{1.5|0.5|1\}, 1\plus \{3|2|-4\}\}$, which is\\ $\{\{4.5|3.5|4\}, \{4.5|3.5|-2.5\}|2.5|\{-2.5|-3.5|-3\},\{4|3|-3\}\}$.

It is also important to note that if $n=\{.|n|.\}$, then $n\plus G=\{G^L\plus n|G^S+n|G^R\plus n\}$, this is easy to see since neither player can move on $n$, so you simply add the score $n$ to all scores on $G$.

\section{Greater Than and Equal To}

We would also like to know when one game is ``better'', than another one.  That is, given several options to play, which one is my best move.  In normal play and mis\`ere play the definitions of ``$\geq$'' and ``$\leq$'', were relatively easy to define, since players either win or lose, however, for scoring play we have to take into account tied scores.  So for this reason I will re-define ``$\geq$'' and ``$\leq$''.

\begin{definition} We define the following:
\begin{itemize}
\item{$-G=\{-G^R|-G^S|-G^L\}$.}
\item{For any two games $G$ and $H$, $G=H$ if $G\plus X$ has the same outcome as $H\plus X$ for all games $X$.}
\item{For any two games $G$ and $H$, $G\geq H$ if $H\plus X\in O$ implies $G\plus X\in O$, where $O=L_\geq$, $R_\geq$, $L_>$ or $R_>$, for all games $X$.}
\item{For any two games $G$ and $H$, $G\leq H$ if $H\plus X\in O$ implies $G\plus X\in O$, where $O=L_\leq$, $R_\leq$, $L_<$ or $R_<$, for all games $X$.}
\item{$G\cong H$ means $G$ and $H$ have identical game trees.}
\item{$G\approx H$ means $G$ and $H$ have the same outcome.}
\end{itemize}
\end{definition}

\subsubsection{An Example}

\begin{figure}[htb]
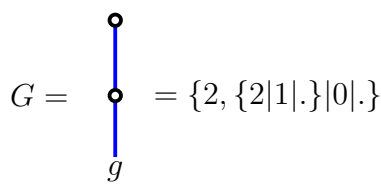

\begin{center}
\begin{graph}(3,2)

\graphnodecolour{1}\graphlinewidth{0.05}

\roundnode{1}(1,0)\roundnode{2}(1,1)\roundnode{3}(1,2)

\edge{1}{2}[\graphlinecolour(0,0,1)]\edge{2}{3}[\graphlinecolour(0,0,1)]

\freetext(0,1){$G=$}\freetext(3,1){$=\{2,\{2|1|.\}|0|.\}$}

\nodetext{1}{$g$}

\end{graph}
\end{center}
\caption{The game $G=\{2,\{2|1|.\}|0|.\}$}\label{g1}
\end{figure}

\begin{figure}[htb]
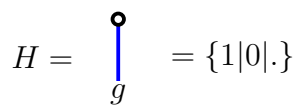

\begin{center}
\begin{graph}(2.5,1)

\graphnodecolour{1}\graphlinewidth{0.05}

\roundnode{1}(1,0)\roundnode{2}(1,1)

\edge{1}{2}[\graphlinecolour(0,0,1)]

\freetext(0,0.5){$H=$}\freetext(2.5,0.5){$=\{1|0|.\}$}

\nodetext{1}{$g$}

\end{graph}
\end{center}
\caption{The game $H=\{1|0|.\}$}\label{h1}
\end{figure}

I will now give examples to demonstrate the definition of ``$\geq$'', since the other definitions are exactly the same as the ones given for normal and mis\`ere games, and so can be easily understood.  So consider the games $G$ and $H$ given in Figures~\ref{g1} and \ref{h1}.

To show that $G\geq H$ we must prove that if $H\plus X\in O$ then $G\plus X\in O$ for all games $X$, where $O=L_>,R_>,L_\geq,R_\geq$.  If $H\plus X\in O$, then Left can do at least as well on $G\plus X$ by simply copying his strategy from $H\plus X$, since if Left removes the single edge at any point during $H\plus X$, then Left can remove the bottom edge on $G\plus X$ at the same point and Left will get 2 points instead of 1, thus earning at least 1 more point than he did on $H\plus X$.  So therefore $G\geq H$.

\begin{theorem}
$G\geq H$ if and only if $H\leq G$
\end{theorem}

\begin{proof}  First let $G\geq H$, and let $G\plus X\in O$ for some game $X$, where $O$ is one of $L_\leq$, $R_\leq$, $L_<$ or $R_<$.  This means that $H\plus X\not\in O'$, where $O'$ is one of $L_\geq$, $R_\geq$, $L_>$ or $R_>$, since if it was this would mean that $G\plus X\in O'$, since $G\geq H$, therefore $H\plus X\in O$, and hence $H\leq G$.

A completely identical argument can be used for $H\leq G$, and hence $G\geq H$ if and only if $H\leq G$ and the theorem is proven.
\end{proof}

The reason why we need this theorem is because we want to make sure that the definition of ``$\geq$'' and ``$\leq$'' are compatible.  This is a natural and desirable property for these two definitions to have.  So for example the game $\{1|1|1\}\geq \{0|0|0\}$ and $\{0|0|0\}\leq \{1|1|1\}$.

\begin{theorem} Scoring play games are partially ordered under the
disjunctive sum.
\end{theorem}

\begin{proof}  To show that we have a partially ordered set we need
3 things.

\begin{enumerate}
\item{Transitivity: If $G\geq H$ and $H\geq J$ then $G\geq J$.}
\item{Reflexivity: For all games $G$, $G\geq G$.}
\item{Antisymmetry: If $G\geq H$ and $H\geq G$ then $G=H$.}
\end{enumerate}

1.  Let $G\geq H$ and $H\geq J$.  $G\geq H$ means that if $H\plus X\in O$ this implies $G\plus X\in O$, where $O=L_\geq$, $R_\geq$, $L_>$ or $R_>$, for all games $X$.  $H\geq J$, means that if $J\plus X\in O$ this implies that $H\plus X\in O$.  Since $G\geq H$, then this implies that $G\plus X\in O$, therefore $J\plus X\in O$ implies that $G\plus X\in O$ for all games $X$, and $G\geq J$.

2. Clearly $G\geq G$, since if $G\plus X\in O$ then $G\plus X\in O$, where  $O=L_\geq$, $R_\geq$, $L_>$ or $R_>$, for all games $X$.

3.  First let $G\geq H$ and $H\geq G$.  $G=H$ means that $G\plus X\approx H\plus X$ for all $X$.  So first let $G\plus X\in L_=$, then this implies that $H\plus X\in L_\geq$, since $H\geq G$.  However $H\plus X\in L_=$, since if $H\plus X\in L_>$, then this implies that $G\plus X\in L_>$, since $G\geq H$, therefore $G\plus X\in L_=$ if and only if $H\plus X\in L_=$.

An identical argument can be used for all remaining cases, therefore $G\plus X\approx H\plus X$ for all games $X$, i.e. $G=H$.
\end{proof}

Under mis\`ere play combinatorial games do not form a partial order, in particular it is possible to find three games $G$, $H$ and $J$ such that $G\geq H$, $H\geq J$ and $J\geq G$ \cite{O}.  So this theorem tells us that these games are a ``nicer'' set under the disjunctive sum than the set of mis\`ere play games, at least in that respect.

Unfortunately however these games are still not as ``nice'' as the set or normal play games under the disjunctive sum.  The following two theorems show that there is very little comparability within the set of scoring play games, and that the identity set is trivial, so these games are not a group under the disjunctive sum.

\begin{theorem}\label{outcome}
For any outcome classes $\mathcal{X}$, $\mathcal{Y}$ and $\mathcal{Z}$, there is a game $G\in\mathcal{X}$ and $H\in\mathcal{Y}$ such that $G\plus H\in\mathcal{Z}$.
\end{theorem}

\noindent\begin{proof}  Consider the games $G=\{\{\{d|c|e\}|b|.\}|a|.\}$ and $H=\{.|f|\{.|g|h\}\}$.  The final scores of $G$ are $\gfl=a$ and $\gfr=b$, and the final scores of $H$ are $H^{SL}_F=f$ and $H^{SR}_F=g$.  Now consider the game $G\plus H$ as shown in the figure.

\begin{figure}[htb]
\begin{center}
\begin{graph}(4.5,3.5)

\roundnode{1}(0,0)\roundnode{2}(1,0)\roundnode{3}(0.5,1)\roundnode{4}(1,2)
\roundnode{5}(1.5,3)\roundnode{6}(3.5,2.5)\roundnode{7}(4,1.5)\roundnode{8}(4.5,0.5)

\edge{1}{3}\edge{2}{3}\edge{3}{4}\edge{4}{5}\edge{6}{7}\edge{7}{8}

\freetext(2.5,1){$\plus$}

\nodetext{1}{$d$}\nodetext{2}{$e$}\nodetext{3}{$c$}\nodetext{4}{$b$}
\nodetext{5}{$a$}\nodetext{6}{$f$}\nodetext{7}{$g$}\nodetext{8}{$h$}

\end{graph}
\end{center}
\caption{The game $G\plus H$, $G=\{\{\{d|c|e\}|b|.\}|a|.\}$ and $H=\{.|f|\{.|g|h\}\}$.}
\end{figure}

The final scores of $G\plus H$ are $(G\plus H)^{SL}_F=e+g$ or $d+h$ and $(G\plus H)^{SR}_F=e+h$.  Since $e$, $d$ and $h$ can take any value we can select them so that: $e+g$, $d+h$ and $e+h>0$ and $G\plus H\in \lef$; $e+g$, $d+h$ and $e+h<0$ and $G\plus H\in \ri$; $e+g$, $d+h>0$ and $e+h<0$ and $G\plus H\in \n$; $e+g$, $d+h<0$ and $e+h>0$ and $G\plus H\in \pre$ or finally $e+g=d+h=e+h=0$ and $G\plus H\in \ti$.

Since the outcomes of $G$ and $H$ depend on the values of $a$, $b$, $f$ and $g$, we can select them so that $G$ and $H$ can be in any outcome class, and thus the theorem is proven.

\end{proof}

This theorem suggests that finding comparable games will be very difficult.  In particular what this theorem says is that there are two games $G, H\in\lef$ such that $G\plus H\in \ri$.  In other words, we have a games $G$ and $H$ where Left has a very big advantage, but if we play them under the disjunctive sum Right now has a big advantage.

This is possibly the worst case of all of them, in comparison, under normal play, for all $G,H\in \le$, $G+H\in\le$.  This is a lot better in terms of comparability between games, since in many cases it allows to easily determine the winner of a game by knowing the outcome classes of the individual components.

Under normal play combinatorial games form an abelian group under the disjunctive sum.  The identity that is used is the set $\pre$, that is if $I\in \pre$ then $G + I\approx G$ for all games $G$, where ``$+$'' here means the disjunctive sum.  In this case the entire set $\pre$ has a single unique representative, the game $\{.|.\}$.  A direct consequence of this is that $G=H$ if and only if $G + (-H)\in \pre$.

Under mis\`ere play, the identity set contains only one element, which is the same game $\{.|.\}$.  That is, if $G\not\cong \{.|.\}$, then $G\neq \{.|.\}$.  This was proven by Paul Ottaway \cite{O}.  Which means that there cannot be an equivalent method for determining if two games are equivalent under mis\`ere play.

For scoring play games, we have a very similar theorem.  That is our identity set contains only one element, namely the game $\{.|0|.\}$, which we call $0$.  It should be clear that $0\plus G\approx G$ for all games $G$, and so $0$ is the identity.  The theorem and proof are given below.

\begin{theorem}\label{identity}
For any game $G$, if $G\not\cong 0$ then $G\neq 0$.
\end{theorem}

\begin{proof}  The proof of this is very simple, first let $G^L\neq\emptyset$, since the case $G^R\neq\emptyset$ will follow by symmetry.  Next let $P=\{.|a|b\}$, and note that $P^{SL}_F=a$, since Left has no move on $P$.  So let $a>0$, if $G=0$, then this means that $(G\plus P)^{SL}_F\approx P$, however since $G$ is a combinatorial game, we know from the definition that $G$ has both finite depth, and finite width.  So we can choose $b<0$ such that $|b|$ is greater than any score on the game tree of $G$.

Therefore when Left moves first on $G\plus P$ he must move to the game $G^L\plus P$, and Right will respond by moving to $G^L\plus b$, and therefore $(G\plus P)^{SL}_F<0$, by choice of $b$, therefore $G\plus P\not\approx P$, and $G\neq 0$.  Hence the theorem is proven.
\end{proof}

What is interesting is that unlike mis\`ere games, some scoring games do have an inverse, namely the set of games $\{.|n|.\}$, where $n$ is a real number.  It should be clear that these are the only games which are invertible under scoring play, and any other non-trivial game cannot be inverted.

\section{Canonical Forms}\label{canfor}

Canonical forms are important, because if we can show that these games can be split up into equivalence classes with a unique representative for each class, then it makes these games much easier to analyze and compare.  We don't have to consider each game individually, but only the equivalence class to which it belongs.

\begin{theorem}\label{eq}
There exist two games $G$ and $H$ such that $G\not\cong H$ and
$G=H$.
\end{theorem}

\begin{proof} Consider the following games $G$ and $H$

\vspace{0.1in}

\begin{figure}[htb]
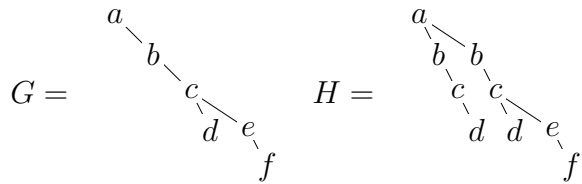

\begin{center}
\begin{graph}(7,2)

\graphnodesize{0.1} \graphlinewidth{0.01}

\roundnode{1}(1,2)\roundnode{2}(1.5,1.5)\roundnode{3}(2,1)\roundnode{4}(2.25,0.5)
\roundnode{5}(2.75,0.5)\roundnode{6}(3,0)\roundnode{7}(5,2)\roundnode{8}(5.25,1.5)
\roundnode{9}(5.75,1.5)\roundnode{10}(5.5,1.0)\roundnode{11}(6,1)\roundnode{12}(5.75,0.5)
\roundnode{13}(6.25,0.5)\roundnode{14}(6.75,0.5)\roundnode{15}(7,0)

\edge{1}{2}\edge{2}{3}\edge{3}{4}\edge{3}{5}
\edge{5}{6}\edge{7}{8}\edge{7}{9}\edge{8}{10}
\edge{10}{12}\edge{9}{11}\edge{11}{13}\edge{11}{14}\edge{14}{15}

\nodetext{1}{$a$}\nodetext{2}{$b$}\nodetext{3}{$c$}\nodetext{4}{$d$}
\nodetext{5}{$e$}\nodetext{6}{$f$}\nodetext{7}{$a$}\nodetext{8}{$b$}
\nodetext{9}{$b$}\nodetext{10}{$c$}\nodetext{11}{$c$}\nodetext{12}{$d$}
\nodetext{13}{$d$}\nodetext{14}{$e$}\nodetext{15}{$f$}

\freetext(0,1){$G=$}\freetext(4,1){$H=$}

\end{graph}
\end{center}
\caption{Two games $G$ and $H$, where $G\not\cong H$, but $G=H$.}
\end{figure}

\noindent where $a,b,c,d,e,f\in \mathbb{R}$.

This example is a variant of a similar example used to prove the
same theorem for mis\`ere games in \cite{GAMO}.

For any two games $G$ and $H$, $G=H$ if $G\plus X\approx H\plus X$ for all games $X$.  The easiest way to prove this is to show that $G\geq H$ and $H\geq G$.  Right can do at least as well playing $H\plus X$ as he can playing $G\plus X$, by simply copying his strategy from $G\plus X$ and not playing the left hand string on $H$.  Right cannot do better on $H\plus X$ than he can on $G\plus X$, since the string on the left hand side of $H$ can be copied on $G\plus X$ by simply not moving to $e$. So therefore if $H\plus X\in O$ then this implies that $G\plus X\in O$ where $O=L_\geq$, $R_\geq$, $L_>$ or $R_>$, i.e. $G\geq H$.

Left can also do at least as well playing $H\plus X$ as he can playing $G\plus X$, since if Right can achieve a lower final score playing the left hand string on $H\plus X$, then he can also do so by choosing not to move to $e$ on $G\plus X$.  Similarly if Right copies his strategy from $G\plus X$ onto $H\plus X$ then their final scores will be the same.  So if $G\plus X\in O$ then this implies that $H\plus X\in O$ where $O=L_\geq$, $R_\geq$, $L_>$ or $R_>$, i.e. $H\geq G$.  So therefore $G=H$ and the proof is finished.
\end{proof}

For both normal and mis\'ere play games, the following methods are used to reduce a game to its canonical form.  We will now demonstrate that they also can be applied to scoring games, and that if we reduce a game using these methods, the games form equivalence classes with unique representatives.

\begin{definition}
Let $G=\{A,B,C,\dots|G^S|D,E,F,\dots\}$, if $A\geq B$ or $D\leq E$ we
say that $A$ dominates $B$ and $D$ dominates $E$.
\end{definition}

\begin{definition}
Let $G=\{A,B,C,\dots|G^S|D,E,F,\dots\}$, an option $A$ is reversible if
$A^R\leq G$.  An option $D$ is also reversible if $D^L\geq G$.
\end{definition}

\begin{theorem}\label{dom}
Let $G=\{A,B,C,\dots|G^S|D,E,F,\dots\}$, and let $A\geq B$, then\\
$G'=\{A,C,\dots|G^S|D,E,F,\dots\}=G$.  By symmetry if $D\leq E$ and\\
$G''=\{A,B,C,\dots|G^S|D,F,\dots\}$ then $G''=G$.
\end{theorem}

\begin{proof}  Let $G=\{A,B,C,\dots|G^S|D,E,F,\dots\}$ such that $A\geq B$, further let\\ $G'=\{A,C,\dots|G^S|D,E,F,\dots\}$.  First suppose that $G\plus X\in O$, where $O=L_\geq$, $R_\geq$, $L_>$ or $R_>$ if Left moves to $B\plus X$.  This implies that $G'\plus X\in O$, since $A\geq B$.  Hence if $G\plus X\in O$  this implies that $G'\plus X\in O$, and since the Right options of $G$ and $G'$, this implies that $G'\geq G$.

Next suppose that $G'\plus X\in O'$ where $O'=L_\leq$, $R_\leq$, $L_<$ or $R_<$.  This implies that $G\plus X\in O'$, since the only option in $G^L$ that is not in $G'^L$ is $B$ and $B\leq A$, therefore $G'\leq G$, and $G=G'$.  So this means that the option $B$ may be disregarded and the proof is finished.

\end{proof}

\subsubsection{Examples of Games with Dominated Options}

Games with dominated options are easy to find, for example consider the following game $G=\{1,2|0|.\}$, then clearly there is no reason why Left would choose to move to $1$ over $2$, since he is guaranteed more points moving to $2\plus X$ than $1\plus X$ for all games $X$.  In other words $1$ is a \emph{dominated} option and $G=\{2|0|.\}$.

Another thing to note is that just because $G\geq H$ for some $G$ and $H$, $H$ may not be a dominated option.  To demonstrate this, consider the game $G'$ in Figure~\ref{domex};

\begin{figure}[htb]
\begin{center}
\begin{graph}(1,3)

\graphnodecolour{1}\graphlinewidth{0.05}

\roundnode{1}(1,0)\roundnode{2}(1,1)\roundnode{3}(1,2)\roundnode{4}(1,3)

\edge{1}{2}[\graphlinecolour(0,0,1)]\edge{2}{3}[\graphlinecolour(0,0,1)]
\edge{3}{4}[\graphlinecolour(0,0,1)]

\freetext(0,1){$G'=$}

\nodetext{1}{$g$}

\end{graph}
\end{center}
\caption{The game $G'=\{3,\{3,\{3|2|.\}|1|.\},\{3|2|.\}|0|.\}$}\label{domex}
\end{figure}

This game has the games $G$ and $H$ from Figure~\ref{g1} and \ref{h1} as two of its options, and as I already established $G\geq H$.  But $H$ is \emph{not} a dominated option.  The reason for that is because now $G^S\neq H^S\neq 0$, so Left cannot simply copy his strategy from $H\plus X$ on $G\plus X$, since there may be an instance where Right ends $H\plus X$ before Left can play on $H$, and if Left copied this on $G\plus X$ he would end up with 1 point less, and therefore may not win.

This is why we have to be very careful when examining scoring play games, and not simply look at a games options, but also the \emph{scores} before we decide if an option is dominated or not.

\begin{theorem}\label{rev1}
Let $G=\{A,B,C,\dots|G^S|D,E,F,\dots\}$, and let $A$ be reversible with
Left options of $A^R=\{W,X,Y,\dots\}$.  If
$G'=\{W,X,Y,\dots,B,C,\dots|G^S|D,E,F,\dots\}$, then $G=G'$.  By symmetry  if $D$ is reversible with Left options $\{T, S, R,\dots\}$ and\\ $G''=\{A,B,C,\dots|G^S|T,S,R,\dots,D,E,F,\dots\}$, then $G=G''$.  
\end{theorem}

\begin{proof} Let $G=\{A,B,C,\dots|G^S|D,E,F,\dots\}$, where the Left options of $A^R=\{W,X,Y,\dots\}$ and let $G'=\{W,X,Y,\dots,B,C,\dots|G^S|D,E,F,\dots\}$, further let $A^R\leq G$.  If $G\plus X\in O$, where $O=L_\geq$, $R_\geq$, $L_>$ or $R_>$, when Left does not move to $A$ on $G$, then clearly $G'\plus X$ is also in $O$, since all other options for Left on $G$ are available for Left on $G'$.  

So consider the case where $G\plus X\in O$ if Left moves to $A\plus X$, then this implies that $A^R\plus X$ must also be in $O$.  This means that $G'\plus X\in O$ because $A^{RL}\subset G'^L$, and since all other options on $G'$ are the same as $G$, then $A^R\plus X\in O$ implies that $G'\plus X\in O$.  Hence if $G\plus X\in O$  then this implies that $G'\plus X\in O$, for all games $X$, i.e. $G'\geq G$.

Next assume that $G\plus X\in O'$, where $O'=L_\leq$, $R_\leq$, $L_<$ or $R_<$, for all games $X$.  However $A^R\leq G$, i.e. $G\plus X\in O'$ implies that $A^R\plus X\in O'$, and since $A^{RL}\subset G'^L$, and all other options on $G'$ are identical to options on $G$, this means that $G\plus X\in O'$, implies that $G'\plus X\in O'$, for all games $X$, i.e. $G'\leq G$.  Therefore $G=G'$ and the theorem is proven.
\end{proof}

\subsubsection{Reversible Options}

Unfortunately games with reversible options are not so easy to find, and I have yet to find a game that fits the definition for a reversible option.  So I leave the following as an open problem.

\begin{question}
Are there any scoring play games with a reversible option?
\end{question}

Although I know this is an open problem, Theorems~\ref{rev1} and \ref{cong} tell us that if such games do exist, then the reduction will still work for scoring play games, i.e. it preserves equality and is a reduction to a games canonical form.  However for scoring play games considering reversibility may in fact be unnecessary.

\begin{theorem}\label{cong}
For any two games $G$ and $H$ if $G=H$ and $G$ and $H$ are in
canonical form then $G\cong H$.
\end{theorem}

\begin{proof} Let $G$ and $H$ be two games such that $G=H$ and neither $G$ nor $H$ has a dominated or reversible option.  

So first let $H\plus X\in O$, where $O=L_<, R_<, L_\leq$ or $R_\leq$, since $G=H$, this implies that $G\plus X\in O$.  However if Left moves to $G^L\plus X$ then $G^{LR}\plus X$ cannot be in $O$, since if it was, this would mean $H\plus X\in O$, implies $G^{LR}\plus X\in O$, i.e. $G^{LR}\leq H$, and $G$ would have a reversible option, which means that $G^L\plus X^R\in O$.  This implies that $H^L\plus X^R\not\in O'$, where $O'=L_>, R_>, L_\geq$ or $R_\geq$, since if it were then $H$ would have a dominated option.  Therefore $G^L\plus X^R\in O$ if and only if $H^L\plus X^R\in O$, i.e. for all $g^L\in G^L$ there is an $h^L\in H^L$ such that $g^L\leq h^L$, and for all $h^L\in H^L$ there is a $g^{L'}\in G^L$ such that $h^L\leq g^{L'}$.

So that means $g^L\leq h^L\leq g'^L$, however $g^L$ and $g^{L'}$ must be identical, otherwise $g^L$ is a dominated option.  So, every Left option of $G$ is equal to a Left option of $H$, i.e. $G^L\subseteq H^L$, and by a symmetrical argument $H^L\subseteq G^L$, i.e. $H^L=G^L$, and similarly $H^R=G^R$.  Therefore $H\cong G$ and the proof is finished. 
\end{proof}

\begin{theorem}
Let $G\to G_1\to G_2\to\dots \to G_n$ represent a series of
reductions on a game $G$ to a game $G_n$, which is in canonical
form.  Further let $G\to G_1'\to G_2'\to\dots\to G_m'$ represent a
different series of reductions on $G$ to a game $G_m'$ which is also
in canonical form, then $G_n\cong G_m'$
\end{theorem}

\begin{proof}  Since each reduction preserves equality, then $G_n=G_m'$ and they are both in canonical form.  By Theorem~\ref{cong} $G_n\cong G_m'$, and so the theorem is proven.
\end{proof}

\section{Summary}

To summarise, what I have done in this chapter is defined a new class of combinatorial games to include scores and also extended the standard definitions such as the disjunctive sum, outcome classes, ``$\ge$'' and ``$\le$'' to this new class. I then showed that under these definitions these games form a partially ordered set and do form equivalence classes with a unique representative, or canonical form.

This is a significant improvement on the previous methods of dealing with games with scores.  The old
combinatorial game theory method, \cite{rjn}, was to simply play the game as usual, then at the end of the game itself, the players take it in turns to make extra ``move'' which consist only of losing one point per move; thus the player
with the most points is the last to move and so wins the game. This method has some merit, but it has a few serious flaws.  The first is that it does not include tied scores.  By converting scoring play games to normal play we still only get four outcome classes, namely $\lef, \ri, \pre$ and $\n$.  Scoring play theory also includes a fifth outcome class, namely $\ti$.

The second serious flaw is that under normal play if $G,H\in\lef$, then $G\plus H\in\lef$.  However as I have shown this is \emph{not} true for scoring play games.  It is much better to simply consider scoring play games as their own class of games, rather than trying to use normal play theory to analyse them.  This method can lead to a much better understanding of scoring play games in general, and tell us the true structure and behaviour of these games.

\chapter{Impartial Games}

In the previous chapter I discussed the structure of scoring play games under the disjunctive sum.  I showed that the games form a partially ordered set under this operator, and that they form equivalence classes with a unique representative using the usual rules of domination and reversibility. 

In this chapter I will be looking at \emph{impartial games}.  In Sections \ref{nim} and \ref{impartial} I explained some of the history of impartial games for normal and mis\`ere play games.  Impartial games and the game of Nim in particular has been so extensively studied for both normal and mis\`ere play games that when we define a new set of games it is a natural subset of games to look at and study.

The definition of an impartial scoring play game is less intuitive than for normal and mis\`ere play games.  The reason for this is because we have to take into account the score, for example, consider the game $G=\{4|3|2\}$.  On the surface the game does not appear to fall into the category of an impartial game, since Left wins moving first or second, however this game is impartial since both players move and gain a single point, i.e. they both have the same options.

So for this thesis I will be using the following definition for an impartial game.

\begin{definition}
A scoring play game $G$ is impartial if it satisfies the following;

\begin{enumerate}
\item{$G^L=\emptyset$ if and only if $G^R=\emptyset$.}
\item{If $G^L\neq \emptyset$ then for all $g^L\in G^L$ there is a $g^R\in G^R$ such that $g^L\plus -G^S=-(g^R\plus -G^S)$.}
\end{enumerate}
\end{definition} 

\begin{figure}[htb]
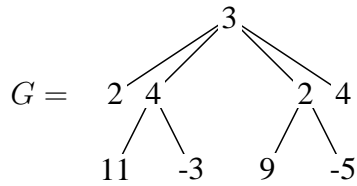

\begin{center}
\begin{graph}(4,2)

\roundnode{1}(1,0)\roundnode{2}(2,0)\roundnode{3}(3,0)\roundnode{4}(4,0)
\roundnode{5}(1,1)\roundnode{6}(1.5,1)\roundnode{7}(3.5,1)\roundnode{8}(4,1)
\roundnode{9}(2.5,2)

\edge{6}{1}\edge{6}{2}\edge{7}{3}\edge{7}{4}
\edge{9}{5}\edge{9}{6}\edge{9}{7}\edge{9}{8}

\nodetext{1}{11}\nodetext{2}{-3}\nodetext{3}{9}\nodetext{4}{-5}
\nodetext{5}{2}\nodetext{6}{4}\nodetext{7}{2}\nodetext{8}{4}
\nodetext{9}{3}

\freetext(0,1){$G=$}

\end{graph}
\end{center}
\caption{The impartial game $G=\{2,\{11|4|-3\}|3|,4,\{9|2|-5\}\}$}\label{ex1}
\end{figure}

\begin{figure}[htb]
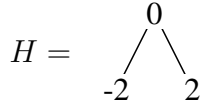

\begin{center}
\begin{graph}(2,1)

\roundnode{1}(1,0)\roundnode{2}(1.5,1)\roundnode{3}(2,0)

\edge{1}{2}\edge{2}{3}

\nodetext{1}{-2}\nodetext{2}{0}\nodetext{3}{2}

\freetext(0,0.5){$H=$}

\end{graph}
\end{center}
\caption{The impartial game $H=\{-2|0|2\}$}\label{ex2}
\end{figure}

Some examples of impartial games are shown in Figures~\ref{ex1} and \ref{ex2}. It should be clear that the game in Figure~\ref{ex2} satisfies the definition of an impartial game since $G^S=0$ and $2=-(-2)$.  The game in Figure~\ref{ex1} also satisfies the definition since $2\plus -3=-(4\plus -3)$, $\{11|4|-3\}\plus -3=-(\{9|2|-5\}\plus -3)=\{8|1|-6\}$ and $11\plus -4=-(-3\plus -4)=9\plus -2=-(-5\plus -2)=7$.  

It is also important to note that an impartial game need not necessarily be symmetrical, in terms of the structure of the game tree.  However if a game is in canonical form then it will be symmetrical.

\begin{definition}
Two impartial games $G$ and $H$ are impartially equivalent if $G\plus X\approx H\plus X$ for all impartial games $X$.
\end{definition}

%Add in Examples

Note that two games $G$ and $H$ that satisfy the definition of impartial equivalence may not actually be equivalent, i.e it may not be true that $G\plus X\approx H\plus X$ for all games $X$.  However, in this chapter I will only be considering impartial games.

\begin{theorem}
Impartial scoring play games form a non-trivial commutative monoid under the disjunctive sum.
\end{theorem}

\begin{proof}  A monoid is a semi-group that has an identity.  Scoring games in general have an identity, but it is a set with only one element, namely $\{.|0|.\}$, see Theorem~\ref{identity}.  I will show that if we restrict scoring play games to impartial games, then there is an identity set that contains more than one element.

First I will define a subset of the impartial games as follows;

$$I=\{i|G\plus i\approx G, \hbox{ for all impartial games }G\}$$  

To show that we have a non-trivial monoid we have to show that $I$ contains more than one element.  So consider the following impartial game, $i=\{\{0|0|0\}|0|\{0|0|0\}\}$.

\begin{figure}[htb]
\begin{center}
\begin{graph}(3,2)

\roundnode{1}(0,0)\roundnode{2}(0.5,1)\roundnode{3}(1,0)\roundnode{4}(1.5,2)
\roundnode{5}(2,0)\roundnode{6}(2.5,1)\roundnode{7}(3,0)

\edge{1}{2}\edge{2}{3}\edge{4}{2}\edge{4}{6}\edge{5}{6}\edge{6}{7}

\nodetext{1}{0}\nodetext{2}{0}\nodetext{3}{0}\nodetext{4}{0}
\nodetext{5}{0}\nodetext{6}{0}\nodetext{7}{0}

\end{graph}
\end{center}
\caption{The game $\{\{0|0|0\}|0|\{0|0|0\}\}$}
\end{figure}

To show that $i\plus G\approx G$ for all impartial games $G$, there are 3 cases to consider $\gfl>0$, $\gfl<0$ and $\gfl=0$, since the cases for Right follow by symmetry.  First let $\gfl>0$, if Left has no move on $G$, then neither does Right, since $G$ is impartial, i.e. $G=G^S$, so they will play $i$ and the final score will still be $G^S$.

So let Left have a move on $G$, Left will choose his best move on $G$.  Right has two choices, either continue to play $G$, or play $i$, and attempt to change the parity of $G$, i.e. force Left to make two consecutive moves on $G$.  However Left will simply respond by also playing $i$, and then it will be Right's turn to move on $G$ again.  Thus $(G\plus i)^{SL}_F>0$.

Next let $\gfl<0$, this means that no matter what Left does, he will lose playing $G$ on $G\plus i$, since Right will simply respond in $G$, until $G$ is finished, then they will play $i$, which does not change the final score of $G$.  Again if Left tries to change the parity of $G$, by playing $i$, Right will also play $i$, and it will be Left's turn to move on $G$ again.  Therefore $(G\plus i)^{SL}_F<0$.

Finally let $\gfl=0$.  This means that Left's best move will be a move that eventually ties $G$.  The only reason Left or Right would choose to move on $i$ is again to change the parity of $G$ and potentially win, i.e. forcing your opponent to move twice on $G$.  However if Left, say, moves on $i$ then Right will simply respond by also playing $i$ and it will be Left's turn to move on $G$ again, and similarly if Right moves on $i$, meaning that $(G\plus i)^{SL}_F=0$.

So therefore the set of impartial games is a non-trivial monoid and the theorem is proven.

\end{proof}

\begin{conjecture}
Not all impartial scoring play games have an inverse.
\end{conjecture}

To prove this one needs to show that given an impartial game $G$ for all impartial games $Y$ there is an impartial game $P$ such that $G\plus Y\plus P\not\approx P$.  This is very difficult to show, however it is extremely likely that this conjecture is true because for normal play games the inverse of any game $G$ is $-G$, and as I will now show there are impartial games $H$ where $-H$ is not the inverse.

So consider the game $G=\{2,\{1|2|3\}|0|-2,\{-3|-2|-1\}\}$, in this case $-G=G$.  If $G$ is the inverse of itself then $G\plus G\plus 0\approx 0$, in other words, $G\plus G\in \ti$.  However $G\plus G\in \pre$, this is easy to see since if Left moves first and moves to $2\plus G$, then Right can respond by moving to $2\plus \{-3|-2|-1\}$ and Left must move to $2\plus -3$ and loses.  If Left moves to $\{1|2|3\}\plus G$, then Right will move to $\{1|2|3\}\plus -2$ and Left must move to $1\plus -2$ and again loses.  Obviously the opposite will be true if Right moves first on $G\plus G$.  So $G\plus G\plus 0\not\approx 0$ and $G\plus G\not\in I$.

So, because $-G$ is not the inverse of $G$ in this case then it is very unlikely that any other impartial game could be $G$'s inverse, and while I can not prove it, I think this simple example shows that the conjecture is probably true.

It is also worth noting that impartial scoring games can belong to any of the five outcomes for scoring games, i.e. $\lef, \ri, \pre, \n$ and $\ti$.  This is in stark contrast to both normal play and mis\'ere play games, where impartial games can only belong to either $\pre$ or $\n$.  

It is easy to see that this is true by considering the impartial game of the form $\{a|G^S|b\}$.  Clearly when $G^S=0$ then $b=-a$ and the outcome can only be $\n, \pre$ or $\ti$.  However we can set $G^S\neq 0$ and either large enough that both $a$ and $b$ are greater than zero, or less than zero, depending on if we make $G^S$ a very large negative or positive number.  In these cases the outcome will either be $\lef$ or $\ri$.

\section{Nim}

Nim is one of the most well studied games in combinatorial game theory, and it is the standard impartial game under normal play.  So when studying impartial games, one of the first and most fundamental games to look is Nim.  

In Section~\ref{nim}, I gave the basic rules for the game of Nim played under normal play rules.  For scoring play Nim we will take the standard rules to be as follows;

\begin{enumerate}
\item{The initial score is $0$.}
\item{The game is played on heaps of beans, and on a player's turn he may remove as many beans as he wishes from any one heap.}
\item{A player gets $1$ point for each bean he removes.}
\item{The player with the most points wins.}
\end{enumerate}

It should be clear that the best strategy for this game is simply to remove all the beans from the largest possible heap, and keep doing so until the game ends.  

Another thing to note is that, under normal play, for every single impartial game $G$ there is a Nim heap of size $n$ such that $G=n$. This not the case with scoring play games, but as I will show in the next section, these games are still relatively easy to solve, regardless of the rules and of the scoring method.

\subsection{Scoring Sprague-Grundy Theory}\label{SGT}

Sprague-Grundy theory is a method that is used to solve any variation of a game of Nim.  The function for normal play $\G(n)$ is defined in a such a way that if for a given heap $n$, played under some rules, if $\G(n)=m$ then this means that the original heap $n$ is equivalent to a Nim heap of size $m$. 

For scoring play games this function is going to be defined slightly differently.  Rather than telling us equivalence classes of different games, it will tell us the final scores of games.  While this may not be as powerful as normal play Sprague-Grundy theory, it is still a very useful function and can be used to solve many different variations of scoring play Nim.

One of the standard variations that have been used widely in books such as Winning Ways \cite{WW}, are a group of games called octal games.  These games cover a very large portion of Nim variations, including all subtraction games.  For scoring games we will use the following definition;

\begin{definition}
A scoring play octal game $O=(n_1n_2\dots n_k, p_1p_2\dots p_k)$, is a set of rules for playing Nim where if a player removes $i$ beans from a heap of size $n$ he gets $p_i$ points, $p_i\in \mathbb{R}$, and he must leave $a,b,c\dots$ or $j$ heaps, where $n_i=2^a+2^b+2^c+\dots+2^j$.
\end{definition}

By convention I will say that a Nim heap $n\in O$ means that $n$ is played under the rule set $O$.  I will now define the function that will be the basis of my theory;

\begin{definition}
Let $n\in O=(t_1t_2\dots t_f, p_1p_2\dots p_t)$ and $m\in P=(s_1s_2\dots s_e, q_1q_2\dots q_t)$;

\begin{itemize}
\item{$\G_s(0)=0$.}
\item{$\G_s(n)=\max_{k,i}\{p_k-\G_s(n_1\plus n_2\plus \dots\plus n_{i})\}$, where $n_1+n_2+\dots +n_{i}=n-k$, $t_k=\Sigma_{i\in S_k}2^i$.}
\item{$\G_s(n\plus m)=\max_{k,i,l,j}\{p_k-\G_s(n_1\plus n_2\plus \dots\plus n_i\plus m),q_l-\G_s(n\plus m_1\plus m_2\plus\dots\plus m_j)\}$, where $n_1+n_2+\dots +n_i=n-k$, $t_k=\Sigma_{i\in S_k}2^i$, $m_1+m_2+\dots m_j=m-l$ and $s_l=\Sigma_{j\in R_l}2^j$.}
\end{itemize}
\end{definition}

The first thing to prove is that this function gives us the information we want, namely the final score of a game.  So we have the following theorem;

\begin{theorem}\label{final}
$\G_s(n)=n^{SL}_F=-n^{SR}_F$ and $\G_s(n\plus m)=(n\plus m)^{SL}_F=-(n\plus m)^{SR}_F$.
\end{theorem}

\begin{proof}  The proof of this will be by induction on all heaps $n_1,n_2,\dots, n_i,m_1\dots,m_j$, such that $n_1+n_2\dots + n_i, m_1+\dots +m_j\leq K$ for some integer $K$, the base case is trivial since $\G_s(0\plus 0\plus 0\dots \plus 0)=0$ regardless of how many $0$'s there are.

So assume that the theorem holds for all $n_1,n_2,\dots, n_i,m_1\dots,m_j$, such that $n_1+n_2\dots + n_i, m_1+\dots +m_j\leq K$ for some integer $K$, and consider $\G_s(n\plus m)$, where $n+m=K+1$.

$\G_s(n\plus m)=\max_{k,i,l,j}\{p_k-\G_s(n_1\plus n_2\plus \dots\plus n_i\plus m),q_l-\G_s(n\plus m_1\plus m_2\plus\dots\plus m_j)\}$, but since $n_1+n_2\dots+n_i+m$ and $n+m_1+m_2\dots+m_j\leq K$, then by induction $\max_{k,i,l,j}\{p_k-\G_s(n_1\plus n_2\plus \dots\plus n_i\plus m),q_l-\G_s(n\plus m_1\plus m_2\plus\dots\plus m_j)\}= \max_{k,i,l,j}\{p_k-(n_1\plus n_2\dots\plus n_i\plus m)^{SL}_F, q_l-(n\plus m_1\plus\dots m_j)^{SL}_F\}=(n\plus m)^{SL}_F$, and the theorem is proven.
\end{proof}

\subsubsection{Subtraction Games}

Subtraction games are a very widely studied subset of octal games.  A subtraction game is a game of Nim where there is a pre-defined set of integers and a player may only remove those numbers of beans from a heap.  This set is called a subtraction set.  From our definition of an octal game this means that each $n_i$ is either $0$ or $3$.  In this section I will also take each $p_i=i$ if $n_i=3$ and 0 otherwise, i.e. if a player removes $i$ beans from a heap he gets $i$ points.

\begin{lemma}
Let $S$ be a finite subtraction set, then for all $s\in S$, $\G_s(s+2ik)=s$ and $\G_s(s+(2i+1)k)=k-s$ for all $i\in\mathbb{Z}^+$, where $k=\max\{S\}$.
\end{lemma}

\begin{proof}  We will split the proof of this into three parts;

\noindent\textbf{Part 1}: For all $i\in\mathbb{Z}^+$, $\G_s(r+2ik)\leq r$

The first thing to show is that for each $0\leq r\leq k$, $\G_s(r)\leq r$ and $\G_s(r+2ik)\leq r$ for all $i\in\mathbb{Z}^+$.  First let $r\leq k$, $\G_s(r)=\max_j\{j-\G_s(r-j)\}$ and since each $j$ in the set is less than or equal to $r$, and each $\G_s(r-j)\geq 0$, this implies that $\G_s(r)\leq r$.

Next let $\G_s(r+2ik)\leq r$ for smaller $i$, and consider $\G_s(r+2ik)=\max_j\{j-\G_s(r+2ik-j)\}$.  If $j\leq r$, then since $\G_s(r+2ik-j)\geq 0$, we have $j-\G_s(r+2ik-j)\leq j\leq r$.  If $j>r$, then $\G_s(r+2ik-j)=\G_s(r+k-j+(2i-1)k)\geq k-(r+k-j)=j-r$, by induction, therefore $j-\G_s(r+2ik-j)\leq j-(j-r)=r$.  So therefore $\G_s(r+2ik)\leq r$ for all $i$.

\noindent\textbf{Part 2}: For all $i\in\mathbb{Z}^+$, $\G_s(r+(2i+1)k)\geq k-r$

We also need to show that for each $0\leq r\leq k$, $\G_s(r+(2i+1)k)\geq k-r$ for all $i\in \mathbb{N}$.  Clearly $\G_s(r+k)\geq k-\G_s(r)\geq k-r$.  Again let $\G_s(r+(2i+1)k)\geq k-r$ for smaller $i$, then $\G_s(r+(2i+1)k)\geq k-\G_s(r+2ik)$ and from above we know that $\G_s(r+2ik)\leq r$ and hence $\G_s(r+(2i+1)k)\geq k-\G_s(r+2ik)\geq k-r$ for all $i$.

\noindent\textbf{Part 3}: For all $s\in S$ and $i\in\mathbb{Z}^+$, $\G_s(s+2ik)\geq s$ and $\G_s(s+(2i+1)k)\leq k-s$.

Let $s\in S$, then $\G_s(s)\geq s-\G_s(0)=s$, since we know from part 1 that $\G_s(s)\leq s$, this means that $\G_s(s)=s$.  So consider $\G_s(s+k)=\max_j\{j-\G_s(s+k-j)\}$, if $j\leq s$ then $j-\G_s(s+k-j)\leq j-k+\G(s-j)\leq j-k+s-j\leq s-k\leq k-s$.  If $j>s$ then $j-\G_s(s+k-j)\leq j-s+\G_s(k-j)\leq j-s+k-j=k-s$.  From part 2 we know that $\G_s(s+k)\geq k-\G_s(s)=k-s$, so $\G_s(s+k)=k-s$.

So assume that the theorem holds up to $i\geq 1$, and consider $\G_s(s+(2i+1)k)=\max_j\{j-\G_s(s+(2i+1)k-j)\}$.  If $j\leq s$ then  $j-\G_s(s+(2i+1)k-j)\leq j-k+\G_s(s+2ik-j)$, and from part 2 we know that $\G_s(s+2ik-j)\leq s-j$ therefore $j-k+\G_s(s+2ik-j)\leq j-k+s-j\leq s-k\leq k-s$.  

If $j>s$ then $j-\G_s(s+(2i+1)k-j)=j-\G_s(s+k+2ik-j)\leq j-s+\G_s(k-j+2ik)\leq j-s+k-j$, by induction, which is equal to $k-s$.

Finally consider $\G_s(s+(2i+2)k)\geq k-\G_s(s+(2i+1)k)$, and from before we know that $\G_s(s+(2i+1)k)\leq k-s$, therefore $k-\G_s(s+(2i+1)k)\geq k-(k-s)=s$.  So therefore $\G_s(s+(2i+2)k)=s$ and the lemma is proven.
\end{proof}

The obvious question to ask is does the lemma hold for all $n$?  The answer is no.  While it is clear that the function is eventually periodic, for subtraction games at least, there are many examples where simply taking the largest number of beans, as in the lemma, is not always the best move.  For example consider a game with subtraction set $\{4,5\}$.  The table of this games $G_s(n)$ values up to $n=15$ are given in Table~\ref{o1}.

\begin{table}[htb]
\begin{center}
\begin{tabular}{|c|c|c|c|c|c|c|c|c|c|c|c|c|c|c|c|c|}
\hline
$n$&0&1&2&3&4&5&6&7&8&9&10&11&12&13&14&15\\\hline
$\G_s(n)$&0&0&0&0&4&5&5&5&5&1&0&0&0&3&4&5\\\hline
\end{tabular}
\end{center}
\caption{A game with subtraction set $\{4,5\}$.}\label{o1}
\end{table}

In particular consider the value of $\G_s(13)=\max\{4-\G_s(9), 5-\G_s(8)\}=4-\G_s(9)=3$.  For this game if you have a heap size of 13 taking 4 beans and gaining 4 points is preferable to taking 5 beans and gaining 5 points.  This is a very simple example to illustrate the point that we cannot say playing greedily would always work.  In other words we need to show that if $n$ is large enough then taking the largest number of beans available \emph{is} the best strategy.  So I make the following conjecture;

\begin{conjecture}
Let $S$ be a finite subtraction set, then there exists an $N$ such that $\G_s(n+k)=k-\G_s(n)$ for all $n\geq N$, where $k=\max\{S\}$.
\end{conjecture}

It seems plausible that this conjecture is true, given the lemma, however it is also possible that there is an $n$ such that $\G_s(n+2ik)=J$ and $\G_s(n+(2i+1)k)=k-j$, where $J>j$.  What we have seen from the data is that often if $n\not\in S$ the values of $\G_s(n+2ik)$ and $\G_s(n+(2i+1)k)$ will alternate as in the lemma, but then you will reach an $i$ where the values change, and this switch might happen several times before it settles down.

I have checked this lemma for all possible subtraction games up to $k=20$ and found no counterexample, so it is far more likely that the conjecture is true than it is false.

A proof of the conjecture or a counterexample would be a very big step forward in understanding how the function operates. 

\subsubsection{Taking-no-Breaking Games}

Taking-no-breaking games are a more general version of subtraction games, and cover a fairly wide range of octal games.  The rules of these games are fairly basic, when a player removes a certain number of beans from a heap, he will have one of three options.

\begin{enumerate}
\item{Leave a heap of size zero, i.e. remove the entire heap.}
\item{Leave a heap of size strictly greater than zero.}
\item{Leave a heap of size greater than or equal to zero.}
\end{enumerate}

From the definition of an octal game this means that each $n_i$ is either $0,1,2$ or $3$, also an octal game $O=(n_1n_2\dots n_k, p_1p_2\dots p_k)$ is finite if $k$ is finite.

It should be clear that for a fixed $m\in P$ and finite $O$, where $P$ and $O$ are two taking no breaking games, then the function $\G_s(n\plus m)$ must always be eventually periodic.  The reason is that we always compute each value from a finite number of previous values, and since $O$ is finite this implies that $\G_s(n\plus m)$ is bounded, and both of these facts together mean that the function will be eventually periodic.

The real question that one needs to answer however is not ``is it periodic?'', but ``what is the period?''.  I believe we can answer that question for a particular class of taking-no-breaking games, that is the class of games where if you remove $i$ beans you get $i$ points.  I make the following conjecture;

\begin{conjecture}\label{period}
Let $O=(n_1n_2\dots n_t, p_1p_2\dots p_t)$ and $P=(m_1m_2\dots m_l, q_1q_2\dots q_l)$ be two finite taking-no-breaking octal games such that, there is at least one $n_s\neq 0$ or $1$, and if $n_i$ and $m_j=1,2$ or $3$ then $p_i=i$ and $q_j=j$, and $p_i=q_j=0$, otherwise, then for all $m$ there exists an $N$ such that;

$$\G_s(n+2k\plus m)=\G_s(n\plus m)$$

\noindent for all $n\geq N$and $k$ is the largest entry in $O$ such that $n_k\neq 0,1$.
\end{conjecture}

There is very strong evidence that this conjecture will hold.  Since $m$ is a constant it changes the value of $\G_s(n\plus m)$, but not the period.  We have checked the conjecture for many examples and not yet found a counterexample, which suggests that it is probably true.

Unfortunately proving it is surprisingly difficult.  The conjecture says that if $n$ is large enough, then your best move is to simply remove the maximum available beans from the heap $n$, so a proof would need to show that for any given $m$, there are only finitely many places where moving on $m$ or removing fewer than $k$ beans from $n$ is a better move.

There are several problems with this, the first is that the function $\G_s(n\plus m)$ only tells us the maximum possible value from the set of possible values.  This makes it very difficult to do a proof that first shows $\G_s(n+2k\plus m)\geq \G_s(n\plus m)$ and vice-versa.  The second is understand \emph{why} removing a lower number of beans would be better than playing greedily in some instances.

The last problem is induction is hard because what may hold for lower values may not hold at higher values, making a proof by induction difficult.  However since the function is recursively defined an inductive proof seems to be more natural than a deductive proof.

I believe that a proof of this theorem would also help in finding the period, and proving it for the more general case, where $i$ beans are work $k$ points, $k\in \mathbb{R}$.

In Tables~\ref{t1} and \ref{t2}, I give two examples which support Conjecture~\ref{period}.  If you examine both tables you will notice that it is eventually periodic in the direction of $n$ and $m$, with periods $k$ and $k'$ respectively, where $k$ is the largest entry such that $n_k\neq 0,1$ and $k'$ the largest entry such that $m_{k'}\neq 0,1$.

\begin{table}[htb]
\begin{center}
\begin{tabular}{c|ccccccccccccc}
$n\plus m$&0&1&2&3&4&5&6&7&8&9&10&11&12\\\hline
0&0&1&0&3&2&3&0&1&2&3&2&1&0\\
1&1&0&1&2&3&2&1&0&1&2&3&2&1\\
2&0&1&0&3&2&3&0&1&2&3&2&1&0\\
3&3&2&3&0&1&2&3&2&1&0&1&2&3\\
4&2&3&2&1&0&1&2&3&2&1&0&1&2\\
5&3&2&3&2&1&0&1&2&3&2&1&2&1\\
6&0&1&0&3&2&1&0&1&2&3&2&1&0\\
7&1&0&1&2&3&2&1&0&1&2&3&2&1\\
8&2&1&2&1&2&3&2&1&0&1&2&3&2\\
9&3&2&3&0&1&2&3&2&1&0&1&2&3\\
10&2&3&2&1&0&1&2&3&2&1&0&1&2\\
11&1&2&1&2&1&2&1&2&3&2&1&0&1\\
12&0&1&0&3&2&1&0&1&2&3&2&1&0\\
\end{tabular}
\end{center}
\caption{$n\in (123, 123)\plus m\in (123,123)$}\label{t1}
\end{table}

\begin{table}[htb]
\begin{center}
\begin{tabular}{c|rrrrrrrrrrrrr}
$n\plus m$&0&1&2&3&4&5&6&7&8&9&10&11&12\\\hline
0&0&1&0&3&2&3&0&1&2&3&2&1&0\\
1&1&0&1&2&3&2&1&0&1&2&3&2&1\\
2&2&1&2&1&2&1&2&1&2&1&2&1&2\\
3&3&2&3&0&1&2&3&2&1&0&1&2&3\\
4&4&3&4&1&2&1&4&3&2&1&2&3&4\\
5&$-3$&4&$-3$&6&5&6&$-3$&$-2$&5&6&5&$-2$&$-3$\\
6&4&$-3$&4&5&6&$-1$&4&3&4&$-1$&0&3&4\\
7&$-3$&4&$-3$&6&5&6&$-3$&$-2$&5&6&5&$-2$&$-3$\\
8&4&$-3$&4&5&6&$-1$&4&3&4&$-1$&0&3&4\\
9&$-3$&4&$-3$&6&5&6&$-3$&$-2$&5&6&5&$-2$&$-3$\\
10&4&$-3$&4&5&6&$-1$&4&3&4&$-1$&0&3&4\\
11&$-3$&4&$-3$&6&5&6&$-3$&$-2$&5&6&5&$-2$&$-3$\\
12&4&$-3$&4&5&6&$-1$&4&3&4&$-1$&0&3&4\\
\end{tabular}
\end{center}
\caption{$n\in (123,123)\plus m\in (3111, 1234)$}\label{t2}
\end{table}

Of course it is natural to ask what happens in the general case, unfortunately in the general case the conjecture doesn't hold.  To see why consider the game $O=(3333, 2222)$.  The values of $G_s(n)$ are given in the following table;

\begin{table}[htb]
\begin{center}
\begin{tabular}{|c|c|c|c|c|c|c|c|c|c|c|c|}
\hline
$n$&0&1&2&3&4&5&6&7&8&9&10\\\hline
$\G_s(n)$&0&2&2&2&2&0&2&2&2&2&0\\\hline
\end{tabular}
\end{center}
\end{table}

This game has period 5, which does not correspond to a possible value of $k$, i.e. 1,2,3 or 4.    While all taking-no-breaking games are periodic as we can see from the example, it is not clear what the period is, since we can take our $p_i$'s to be any real number.  So I make the following conjecture;

\begin{conjecture}Let $O=(n_1n_2\dots n_t, p_1p_2\dots p_k)$ and $P=(m_1m_2\dots m_l, q_1q_2\dots q_l)$ be two finite taking-no-breaking octal games, then there exists a $t$ such that;

$$\G_s(n+t\plus m)=\G_s(n\plus m)$$

\end{conjecture}

\subsubsection{Taking-and-Breaking}

Another type of Nim games we can examine are taking-and-breaking games.  That is games where after the player removes some beans from a heap he must break the remainder into two or more heaps. This is more general than taking-no-breaking games, since taking-no-breaking games are a subset of taking-and-breaking games.

There are several problems with examining taking and breaking scoring games.  The first is that we cannot even say that the function $\G_s(n\plus m)$ is bounded.  The reason is that with each iteration you are increasing the number of heaps, which may increase the value of the function as $n$ increases.  So we cannot put a bound on the function as we could with subtraction games and taking-no-breaking games.

Another problem is that if we were to say examine the game $O=(26, p_1p_2)$, which means take one bean bean and leave one non-empty heap, or take two beans and leave either two non-empty heaps, or one non-empty heap, the number of computations required to find $\G_s(n)$ increases exponentially with $n$.  Since a heap of size $n-2$ may be broken into two smaller heaps $n_1$ and $n_2$, we must therefore also compute the value of $\G_s(n_1\plus n_2)$.  

However if $n_1-2=$ or $n_2-2$ may also be broken into two smaller heaps, say $n_1'$, $n_1''$, $n_2'$ and $n_2''$ then we must compute the value of $\G_s(n_1'\plus n_1''\plus n_2)$ and $\G_s(n_1\plus n_2'\plus n_2'')$.  This process will continue until we have heaps that are too small to be broken up.  So this means that computing $\G_s(n)$ for a taking and breaking game is a lot harder, than for a taking-no-breaking game, simply due to the number of computations involved.

So I leave the following open problem.

\begin{question}
Let $O=(n_1n_2\dots n_t, p_1p_2\dots p_k)$ and $P=(m_1m_2\dots m_l, q_1q_2\dots q_l)$ be two finite octal games.  Does there exists a $t$ such that;

$$\G_s(n+t\plus m)=\G_s(n\plus m)?$$

\noindent If so what is the value of $t$?
\end{question}

While I feel that the answer to this is possibly ``yes'', there is very little evidence to support it which is why I leave it as an open problem and do not state it as a conjecture.  Studying these games would certainly be interesting and anything anyone could find out about them would be useful.

\section{Summary}

To summarise, in this chapter I looked at impartial scoring play games, and demonstrated that under scoring play rules an impartial game can belong to any outcome class.  I also looked at the game of Nim, and gave an equivalent Sprague-Grundy theory for solving many different variants of the game of Nim.

This is quite a remarkable result.  As I demonstrated in chapter 2, scoring play games do not have a particularly ``nice'' structure.  That is they do not form a group.  So one would expect that when studying impartial games it would be quite ``wild'', as it is for mis\`ere play games, and as a result rather difficult to analyse.

But as I've shown in this chapter, they are not so ``wild'', and there is a simple algorithm that can be used to solve many different variations of the game of Nim.  This is something that is \emph{not} true for mis\`ere play.  I am very confident that the conjectures I gave in this chapter are all true, and a proof would no doubt help us to understand the function $\G_s(n)$ a lot better.

\chapter{Different Operators}

In chapters 2 and 3 I looked at scoring play games under the disjunctive sum, however there are many ways to play combinatorial games and combine them together.  In this chapter we will be looking at the structure of scoring play games when you use different methods to combine them.

There are three main operators that I will be examining in this chapter.  The conjunctive sum, where players play all components on each turn, the selective sum, where players choose which components to play, and the sequential join, where we arrange the components in order and players must play them in that order.  Note that all operators will be played under the long rule, which was described in Section~\ref{disjunctive}.

There are three main questions that I will be attempting to answer about each operator.  Do scoring play games form a group or non-trivial monoid under the operator?  What are the sums of games from different outcome classes?  Is there a Sprague-Grundy theory under the operator?  These three questions are the most important questions, since they form the basis for any future analysis.

\section{The Conjunctive Sum}

The first operator that I will be looking at is the conjunctive sum.  Under this operator, players must move on all components, where they have a move, on their turn.  Mathematically it is defined as follows;

\begin{definition} The conjunctive sum is:

$$G\bigtriangleup H=\{G^L\bigtriangleup H^L|G^S+H^S|G^R\bigtriangleup H^R\}$$

\noindent where $G^S+H^S$ is the normal addition of two real numbers.
\end{definition}

Note that if $H^L=\emptyset$ and $G^L\neq\emptyset$ then $G^L\bigtriangleup H^L=G^L\bigtriangleup H$ under the long rule, and $G^L\bigtriangleup H^L=\emptyset$ under the short rule, and similarly for Right.

\subsection{Conjunctive Sum Examples}\label{conjex}

To highlight the differences between each of the different operators I will be using two examples.  So consider the following games;

The first example is from Toads and Frogs described in chapter 2 Section~\ref{ch2ex}.

\begin{figure}[htb]
\begin{center}
\begin{pspicture}(9,5)
\put(0,2){$BTF\bigtriangleup TFB=$}\put(4,4){\cgtree{
  {BTF(0)}(|{FTB(-1)}({FBT(-1)}|))
  }}
   \put(9,4){\cgtree{
  {TFB(0)}({BTF(1)}(|{FBT(1)})|)
  }}
\put(6.5,2){$\bigtriangleup$}
\end{pspicture}
\end{center}
\caption{$BTF\bigtriangleup TFB=\{.|0|\{-1|-1|.\}\}\bigtriangleup \{\{.|1|1\}|0|.\}$}
\end{figure}

This game is $BTF\bigtriangleup TFB=\{.|0|\{-1|-1|.\}\}\bigtriangleup \{\{.|1|1\}|0|.\}$\\ $=\{\{.|1|\{0|0|.\}\}|0|\{\{.|0|0\}|-1|.\}\}$.  

The second example I'll be looking is from Hackenbush.  So consider the game in Figure~\ref{h.ex}.  This game has the value $\{\{.|1|-1\}|0|\{\{.|0|-1\}|-1|-3\}\}$ under the conjunctive sum.

\begin{figure}[htb]
\begin{center}
\begin{graph}(1,1)

\graphnodecolour{1}\graphlinewidth{0.05}

\roundnode{4}(0,0)
\roundnode{5}(0,1)\roundnode{6}(1,0)\roundnode{7}(1,1)

\edge{4}{5}[\graphlinecolour(0,0,1)]
\edge{5}{7}[\graphlinecolour(1,0,0)]\edge{6}{7}[\graphlinecolour(0,0,1)]

\nodetext{4}{$g$}\nodetext{6}{$g$}

\end{graph}
\end{center}
\caption{The game $\{\{.|1|-1\}|0|\{\{.|0|-1\}|-1|-3\}\}$}\label{h.ex}
\end{figure}

\begin{theorem}
If $G\not\cong 0$ then $G\neq 0$.
\end{theorem}

\begin{proof}  The proof of this is very similar to the proof for Theorem~\ref{identity}.  First consider the game $G^L=G^R=\emptyset$, then clearly if $G^S\neq 0$ then $G\neq 0$.

Next consider the case where $G^L\neq \emptyset$, since the case $G^R\neq \emptyset$ follows by symmetry.  Let $P=\{.|a|b\}$, where $a=P^{SL}_F>0$.  Since $G$ is a combinatorial game, this means that the game tree has finite depth and finite width, we can choose $b$ to be more negative than any number on $G$.  On Left's first turn he must move to $G^L\bigtriangleup P$, regardless of whether Right can play $G$ or not, he will have to move on $P$ on his next turn.

Thus $(G\bigtriangleup P)^{SL}_F<0$, and therefore $G\bigtriangleup P\not\approx P$, and the theorem is proven.  
\end{proof}

As with the disjunctive sum, this theorem is saying that under the conjunctive sum the identity set is trivial, i.e. it contains only one element.  In other words the set of scoring play games do not form a group under the conjunctive sum.

\begin{theorem}
For any three outcome classes $\mathcal{X}$, $\mathcal{Y}$ and $\mathcal{Z}$, there is a game $G\in \mathcal{X}$ and $H\in \mathcal{Y}$ such that $G\bigtriangleup H\in \mathcal{Z}$.
\end{theorem}

\begin{proof}  To prove this consider the following game $G\bigtriangleup H$, where\\ $G=\{\{.|b|\{.|c|\{e|d|f\}\}\}|a|.\}$ and $H=\{.|g|\{\{\{k|j|.\}|i|.\}|h|.\}\}$, as shown in Figure~\ref{a}.

\begin{figure}[htb]
\begin{center}
\begin{graph}(4.5,4)
\roundnode{1}(0.5,4)\roundnode{2}(0,3)\roundnode{3}(0.5,2)\roundnode{4}(1,1)\roundnode{5}(1.5,0)
\roundnode{6}(0.5,0)\roundnode{7}(3,0)\roundnode{8}(3.5,1)\roundnode{9}(4,2)\roundnode{10}(4.5,3)
\roundnode{11}(4,4)

\edge{1}{2}\edge{2}{3}\edge{3}{4}\edge{4}{5}\edge{4}{6}
\edge{7}{8}\edge{8}{9}\edge{9}{10}\edge{10}{11}

\freetext(2.25,2){$\bigtriangleup$}

\nodetext{1}{$a$}\nodetext{2}{$b$}\nodetext{3}{$c$}\nodetext{4}{$d$}\nodetext{5}{$f$}
\nodetext{6}{$e$}\nodetext{7}{$k$}\nodetext{8}{$j$}\nodetext{9}{$i$}\nodetext{10}{$h$}
\nodetext{11}{$g$}
\end{graph}
\end{center}
\caption{$\{\{.|b|\{.|c|\{e|d|f\}\}\}|a|.\}\bigtriangleup \{.|g|\{\{\{k|j|.\}|i|.\}|h|.\}\}$}\label{a}
\end{figure}

In these games $G^{SL}_F=c$ and $G^{SR}_F=a$, $H^{SL}_F=g$ and $H^{SR}_F=i$, however $(G\bigtriangleup H)^{SL}_F=e+j$ and $(G\bigtriangleup H)^{SR}_F=e+k$.  Since the outcome classes of $G$ and $H$ depend on $a$, $c$, $g$ and $i$, and the outcome class of $G\bigtriangleup H$ depends on $e+j$ and $e+k$, then clearly we can choose $a$, $c$, $g$, $i$, $e$, $j$ and $k$, so that $G$ and $H$ can be in any outcome class and $G\bigtriangleup H$ can be in any outcome class and the theorem is proven.
\end{proof}

\subsection{Impartial Games}

\begin{theorem}
Impartial games form an abelian group under the conjunctive sum.
\end{theorem}

\begin{proof}  To prove this we only need to show that there is an identity set $I$ that contains more than one element, and that for any impartial game $G$, there is a $G^{-1}$ such that $G\bigtriangleup G^{-1}\in I$.  We will split the proof into 3 main parts, Identity, Inverses and remaining properties.

\noindent\textbf{Part 1:} Identity

For the first part of the proof we need to show that there is a non-trivial identity set, so let $I=\{G|G \hbox{ is impartial and }G\in \ti\}$. We then need to show that for all $G\in I$, $G\bigtriangleup P\approx P$ for all impartial games $P$.  There are three cases to consider, since the remaining follow by symmetry, $P\in L_>$, $P\in L_<$ or $P\in L_=$.

\noindent\textbf{Part 1(a)}: $P\in L_>$

So first let $P\in L_>$, and consider the game $G\bigtriangleup P$.  Since Left can achieve a score of $0$ on $G$, then all Left has to do is play his winning strategy on $P$, and $G\bigtriangleup P\in L_>$.

\noindent\textbf{Part 1(b)}: $P\in L_<$

Next let $P\in L_<$, and consider the game $G\bigtriangleup P$.  $G\in L_=$, and since both $G$ and $P$ are impartial, neither player can change the parity of either game, since they must both play both games on every turn.     So all Right has to do is play his winning strategy on $P$ and $G\bigtriangleup P\in L_<$.

\noindent\textbf{Part 1(c)}: $P\in L_=$

Finally let $P\in L_=$, and consider the game $G\bigtriangleup P$.  If both players always make their best moves on $G$ and $P$ then the final score of $G\bigtriangleup P$ will be 0, since $G\in \ti$ and $P\in L_=$, this also implies that if Left chooses a different move other than his best move either $G$ or $P$, then the final score will be $\leq 0$, and similarly for Right.  

This means that as long as Right keeps playing his best strategy, if Left chooses anything else Right can potentially win and similarly for Left.  In other words the best thing for both players to do is to play their best strategy on both $G$ and $P$ and the final score will be a tie, i.e.  $G\bigtriangleup P\in L_=$.

The cases for $R_>$, $R_<$ and $R_=$ follow by symmetry.

\noindent\textbf{Part 2}: Inverses

For the inverse of a game $G$, where $G^{SL}_F=n$ and $G^{SR}_F=p$, we let $H$ be a game where $H^{SL}_F=-n$ and $H^{SR}_F=-p$.  Note that $G\bigtriangleup H\in I$ if and only if $G\bigtriangleup H\in \ti$.

So consider the game $G\bigtriangleup H$ with Left moving first, since the case where Right moves first follows by symmetry.  If $G^L=\emptyset$, then this implies that $G^R=\emptyset$ since $G$ is impartial, which implies that $G^{SL}_F=G^{SR}_F=n$, so for the inverse let $H$ be a game such that $H^{SL}_F=H^{SR}_F=-n$.  However since $H$ is impartial the only game that satisfies that condition is the game $H=\{.|-n|.\}$, which is clearly the inverse of $G$.

If $G^L\neq \emptyset$, and Left and Right make their best move at every stage on both $G$ and $H$, then the final score of $G\bigtriangleup H$ will be $G^{SL}_F+H^{SL}_F=n-n=0$.  Using the same argument as the identity proof if Left or Right try a different strategy then the final score will be either $\leq 0$ or $\geq 0$ respectively, therefore $(G\bigtriangleup H)^{SL}_F=(G\bigtriangleup H)^{SR}_F=0$ and $H$ is the inverse of $G$.

\noindent\textbf{Part 3}: Remaining Properties

It is clear that the set is closed, since if $G$ and $H$ are impartial then $G\bigtriangleup H$ must also be impartial.  It is also clear that we have commutativity and associativity, since we must play on every component on every turn, then the order of the components is irrelevant.
\end{proof}

\subsection{Sprague-Grundy Theory}

Sprague-Grundy theory is much easier under the conjunctive sum than the disjunctive sum, the fact that impartial games are a group mean that we can easily solve any octal game simply by knowing each heap's $\G_s(n)$ value.  So first we define the following;

\begin{definition}
Let $n\in O=(t_1t_2\dots t_f, p_1,\dots p_f)$ and $m\in P=(s_1s_2\dots s_e, q_1\dots q_e)$;

\begin{itemize}
\item{$\G_s(0)=0$.}
\item{$\G_s(n)=\max_{k,i}\{p_k-\G_s(n_1\bigtriangleup n_2\bigtriangleup \dots\bigtriangleup n_{i})\}$, where $n_1+n_2+\dots +n_{i}=n-k$ and $t_k=\Sigma_{i\in S_k}2^i$.}
\item{$\G_s(n\plus m)=\max_{k,i,l,j}\{p_k+q_l-\G_s(n_1\bigtriangleup n_2\bigtriangleup \dots\bigtriangleup n_i\bigtriangleup m_1\bigtriangleup m_2\bigtriangleup\dots\bigtriangleup m_j)\}$, where $n_1+n_2+\dots +n_i=n-k$,$t_k=\Sigma_{i\in S_k}2^i$, $m_1+m_2+\dots m_j=m-l$ and $s_l=\Sigma_{j\in R_l}2^j$.}
\end{itemize}
\end{definition}

\begin{theorem}
$$\G_s(n\bigtriangleup m)=\G_s(n)+\G_s(m)$$
\end{theorem}

\begin{proof}  We will prove this by induction.  The base case is trivial, since $\G_s(0\bigtriangleup 0)=\G_s(0)+\G_s(0)=0$.

So assume that the theorem holds for all values up to $\G_s(n\bigtriangleup m)$, and consider $\G_s(n+1\bigtriangleup m)$, since the case $\G_s(n\bigtriangleup m+1)$ follows by symmetry.  $\G_s(n+1\bigtriangleup m)=\max_{k,i,l,j}\{p_k+q_l-\G_s(n_1\bigtriangleup n_2\bigtriangleup \dots\bigtriangleup n_i\bigtriangleup m_1\bigtriangleup m_2\bigtriangleup\dots\bigtriangleup m_j)\}$, where $n_1+n_2+\dots +n_i=n+1-k$, $t_k=\Sigma_{i\in S_k}$, $m_1+m_2+\dots m_j=m-l$ and $s_l=\Sigma_{j\in R_l}$.  However each $n_i'<n+1$ and $m_j'<m$, therefore $\max_{k,i,l,j}\{p_k+q_l-\G_s(n_1\bigtriangleup n_2\bigtriangleup \dots\bigtriangleup n_i\bigtriangleup m_1\bigtriangleup m_2\bigtriangleup\dots\bigtriangleup m_j)\}=\max_{k,i,l,j}\{k+l-\G_s(n_1)+\G_s(n_2)+\dots+\G_s(n_i)+\G_s(m_1)+\G_s(m_2)+\dots+\G_s(m_j)\}$, by induction.  

Therefore $\max_{k,i,l,j}\{p_k+q_ll-\G_s(n_1)+\G_s(n_2)+\dots+\G_s(n_i)+\G_s(m_1)+\G_s(m_2)+\dots+\G_s(m_j)\}=\max_{k,i}\{p_k-G_s(n_1)+\G_s(n_2)+\dots+\G_s(n_i)\}+\max_{m,j}\{q_l-\G_s(m_1)+\G_s(m_2)+\dots+\G_s(m_j)\}=\max_{k,i}\{p_k-\G_s(n_1\bigtriangleup n_2\bigtriangleup\dots\bigtriangleup n_i)\}+\max_{m,j}\{q_l-\G_s(m_1\bigtriangleup m_2\bigtriangleup\dots\bigtriangleup m_j)\}=\G_s(n+1)+\G_s(m)$, and the proof is finished.
\end{proof}

So unlike the disjunctive sum, all Nim games are very easy to analyse, and form a nicer structure than they do under the disjunctive sum. 

The reason for this is because it is \emph{not} true that $\G_s(n\plus m)=\G_s(n)+\G_s(m)$.  So unlike the disjunctive sum we can determine the value of $\G_s(n\bigtriangleup m)$ simply by knowing what $\G_s(n)$ and $\G_s(m)$ are.

\section{The Selective Sum}

The selective sum is a more general version of the disjunctive sum.  Rather than choosing a single component on each turn and playing that one only, the player can select any components he wishes to play and play those components on his turn.  It is defined as follows;

\begin{definition}The selective sum is:

$$G\triangledown H=\{G^L\triangledown H, G\triangledown H^L, G^L\triangledown H^L|G^S+H^S|G^R\triangledown H, G\triangledown H^R, G^R\triangledown H^R\}$$

\noindent where $G^S+H^S$ is the normal addition of two real numbers.
\end{definition}

As with the disjunctive sum the ``,'' means set union, and $G^L\triangledown H$, $G^L\triangledown H^L$, mean take the selective sum of all $g^L\in G^L$ with $H$ and all $h^L\in H^L$.  Again note that if $H^L=\emptyset$ and $G^L\neq\emptyset$ then $G^L\triangledown H, G\triangledown H^L, G^L\triangledown H^L=G^L\triangledown H$ under the long rule and $\emptyset$ under the short rule.

\begin{theorem}
If $G\not\cong 0$ then $G\neq 0$.
\end{theorem}

\begin{proof}  The proof of this is very similar to the proof for an equivalent theorem given in chapter 2.  First consider the game $G^L=G^R=\emptyset$, then clearly if $G^S\neq 0$ then $G\neq 0$.

Next consider the case where $G^L\neq \emptyset$, since the case $G^R\neq \emptyset$ follows by symmetry.  Let $P=\{.|a|b\}$, where $a=P^{SL}_F>0$.  Since $G$ is a combinatorial game, this means that the game tree has finite depth and finite width, we can choose $b$ to be more negative than any number on $G$.  On Left's first turn he must move to $G^L\triangledown P$, Right can then win by simply moving to $G^L\triangledown b$ on his turn, since the final score will be less than $0$, regardless of what Left does.

Thus $(G\triangledown P)^{SL}_F<0$, and therefore $G\triangledown P\not\approx P$, and the theorem is proven.     
\end{proof}

\begin{theorem}
For any three outcome classes $\mathcal{X}$, $\mathcal{Y}$ and $\mathcal{Z}$, there is a game $G\in \mathcal{X}$ and $H\in \mathcal{Y}$ such that $G\bigtriangleup H\in \mathcal{Z}$.
\end{theorem}

\begin{proof}  To prove this consider the following game $G\triangledown H$, where $G= \{\{c|b|.\}|a|.\}$ and $H=\{.|d|\{.|e|\{.|f|g\}\}\}$, as shown in the following diagram.  

\begin{figure}[htb]
\begin{center}
\begin{graph}(3.5,3)

\roundnode{1}(0,1)\roundnode{2}(0.5,2)\roundnode{3}(1,3)\roundnode{4}(2,3)
\roundnode{5}(2.5,2)\roundnode{6}(3,1)\roundnode{7}(3.5,0)

\edge{1}{2}\edge{2}{3}\edge{4}{5}\edge{5}{6}\edge{6}{7}

\freetext(1.5,1.5){$\triangledown$}

\nodetext{1}{$a$}\nodetext{2}{$b$}\nodetext{3}{$c$}\nodetext{4}{$d$}
\nodetext{5}{$e$}\nodetext{6}{$f$}\nodetext{7}{$g$}

\end{graph}
\end{center}
\caption{$\{\{c|b|.\}|a|.\}\triangledown \{.|d|\{.|e|\{.|f|g\}\}\}$}
\end{figure}

In these games $G^{SL}_F=b$ and $G^{SR}_F=a$, $H^{SL}_F=d$ and $H^{SR}_F=e$, however $(G\triangledown H)^{SL}_F=c+f$ and $(G\bigtriangleup H)^{SR}_F=c+g$.  Since the outcome classes of $G$ and $H$ depend on $a$, $b$, $d$ and $e$, and the outcome class of $G\triangledown H$ depends on $c+f$ and $c+g$, then clearly we can choose $a$, $b$, $c$, $d$, $e$, $f$ and $g$, so that $G$ and $H$ can be in any outcome class and $G\triangledown H$ can be in any outcome class and the theorem is proven.

\end{proof}

\subsection{Selective Sum Examples}\label{sse}

Following on from the examples from Section~\ref{conjex}, again consider the same games, this time played under the selective sum.

\begin{figure}[htb]
\begin{center}
\begin{pspicture}(9,5)
\put(0,2){$BTF\triangledown TFB=$}\put(4,4){\cgtree{
  {BTF(0)}(|{FTB(-1)}({FBT(-1)}|))
  }}
   \put(9,4){\cgtree{
  {TFB(0)}({BTF(1)}(|{FBT(1)})|)}}
\put(6.5,2){$\triangledown$}
\end{pspicture}
\end{center}
\caption{$BTF\triangledown TFB=\{.|0|\{-1|-1|.\}\}\triangledown\{\{.|1|1\}|0|.\}$}
\end{figure}

The game $BTF\triangledown TFB=\{.|0|\{-1|-1|.\}\}\triangledown\{\{.|1|1\}|0|.\}=$\\ $\{\{.|1|\{.|1|\{0|0.\}\},\{0|0|.\},\{\{.|0|0\}|0|\{0|0|.\}\}\}|0|\{\{\{.|0|0\}|0|\{0|0|.\}\},\{.|0|0\},\{\{.|0|0\}|\\-1|.\}|-1|.\}\}$.  By comparison under the conjunctive sum this game had value\\ $\{\{.|1|\{0|0|.\}\}|0|\{\{.|0|0\}|-1|.\}\}$, which is totally different.

The game in Figure~\ref{h.ex2} has value $\{\{.|1|-1,\{.|0|-1\}|0|\{\{.|0|1\}|-1|-3\}\}$ under the selective sum.  The reason for the difference in value from Figure~\ref{h.ex} is that when Left removes the blue edge, it leaves two red edges which are played under the selective sum not the conjunctive sum.  

\begin{figure}[htb]
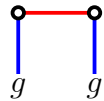

\begin{center}
\begin{graph}(1,1)

\graphnodecolour{1}\graphlinewidth{0.05}

\roundnode{4}(0,0)
\roundnode{5}(0,1)\roundnode{6}(1,0)\roundnode{7}(1,1)

\edge{4}{5}[\graphlinecolour(0,0,1)]
\edge{5}{7}[\graphlinecolour(1,0,0)]\edge{6}{7}[\graphlinecolour(0,0,1)]

\nodetext{4}{$g$}\nodetext{6}{$g$}

\end{graph}
\end{center}
\caption{The game $\{\{.|1|-1,\{.|0|-1\}|0|\{\{.|0|1\}|-1|-3\}\}$}\label{h.ex2}
\end{figure}

\subsection{Impartial Games}

\begin{theorem}
Impartial games form a non-trivial monoid under the selective sum.
\end{theorem}

\begin{proof}  As with the disjunctive sum to prove that we have a non-trivial monoid we simply need to define an identity set that contains more than the game $\{.|0|.\}$.

First I will define a subset of the impartial games as follows;

$$I=\{i|G\plus i\approx G, \hbox{ for all impartial games }G\}$$  

Again, in order to show that we have a \emph{non-trivial} monoid we have to show that $I$ contains more than one element.  So consider the following impartial game,\\ $i=\{\{0|0|0\}|0|\{0|0|0\}\}$.

\begin{figure}[htb]
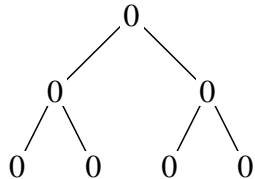

\begin{center}
\begin{graph}(3,2)

\roundnode{1}(0,0)\roundnode{2}(0.5,1)\roundnode{3}(1,0)\roundnode{4}(1.5,2)
\roundnode{5}(2,0)\roundnode{6}(2.5,1)\roundnode{7}(3,0)

\edge{1}{2}\edge{2}{3}\edge{4}{2}\edge{4}{6}\edge{5}{6}\edge{6}{7}

\nodetext{1}{0}\nodetext{2}{0}\nodetext{3}{0}\nodetext{4}{0}
\nodetext{5}{0}\nodetext{6}{0}\nodetext{7}{0}

\end{graph}
\end{center}
\caption{The game $\{\{0|0|0\}|0|\{0|0|0\}\}$}
\end{figure}

To show that $i\triangledown G\approx G$ for all impartial games $G$, there are 3 cases to consider $\gfl>0$, $\gfl<0$ and $\gfl=0$, since the cases for Right follow by symmetry.  First let $\gfl>0$, if Left has no move on $G$, then neither does Right, since $G$ is impartial, i.e. $G=G^S$, so they will play $i$ and the final score will still be $G^S$.

So let Left have a move on $G$, if Left chooses his best move on $G$,  then if Right plays $i$, then Left will respond in $i$ and Right must play $G$, which Left wins.  If Right tries to play on both $G$ and $i$, then either Right moves to a game where $G^L=G^R=\emptyset$, in which case Left moves on $i$ only and wins, or $G^L\neq \emptyset$, and Left also plays both $G$ and $i$ in order to maintain parity on $G$ and still wins.  Clearly if Right chooses to play $G$, then he will still lose, since Left also plays $G$ until it is finished and neither player can gain points on $i$.

Next let $\gfl<0$, this means that no matter what Left does, he will lose playing only $G$ on $G\triangledown i$, since Right will simply respond in $G$, until $G$ is finished, then they will play $i$, which does not change the final score of $G$.  Again if Left tries to change the parity of $G$, by playing $i$, Right will also play $i$, and it will be Left's turn to move on $G$ again.  If Left chooses to move on both $G$ and $i$, then as before Right will also move on $G$ and $i$ if $G^R\neq\emptyset$, and $i$ if $G^R=\emptyset$, but will win either way.

Finally let $\gfl=0$.  This means that Left's best move will be a move that eventually ties $G$.  So consider the game $G\bigtriangleup i$, Left's best move will be to move either on $G$ or $G$ and $i$, if Left moves on $i$ then this will give Right an opportunity to move first on $G$ and potentially win.  If Left moves on $G$ then Right can either play $G$, $i$ or $G$ and $i$.  If Right chooses to play $G$ then Left will simply respond in $G$ to force a tie, if Right plays $i$ then Left can either respond in $i$ and still tie, or play $G$ and potentially win.  If Right plays both $G$ and $i$, again Left can respond in both and tie, or play $i$ only and potentially win.  So therefore $(G\bigtriangleup i)^{SL}_F=0$.

Therefore the set of impartial games is a non-trivial monoid under the selective sum and the theorem is proven.

\end{proof}

\begin{conjecture}\label{cs1}
Not every impartial game is invertible under the selective sum.
\end{conjecture}

As with the disjunctive sum, it is quite likely that this conjecture is true, but proving it is very difficult.  Like the disjunctive sum we need to show that there exists an impartial game $G$ such that for all impartial games $Y$ there exists an impartial game $P$ such that $G\triangledown Y\triangledown P\not\approx P$.  In other words there are no impartial games $Y$ that would serve as an inverse for $G$.

However since the selective sum is a more general version of the disjunctive sum, if there are no inverses for the disjunctive sum, then it is even less likely that there would be inverses under the selective sum.

\subsection{Sprague-Grundy Theory}

As with the other operators I will define my function in the most general possible sense.

\begin{definition}
Let $n\in O=(t_1t_2\dots t_f, p_1,\dots p_f)$ and $m\in P=(s_1s_2\dots s_e, q_1\dots q_e)$;

\begin{itemize}
\item{$\G_s(0)=0$.}
\item{$\G_s(n)=\max_{k,i}\{p_k-\G_s(n_1\plus n_2\triangledown \dots\triangledown n_{i})\}$, where $n_1+n_2+\dots +n_{i}=n-k$ and $t_k=\Sigma_{i\in S_k}2^i$.}
\item{$\G_s(n\plus m)=\max_{k,i,l,j}\{p_k - \G_s(n_1\triangledown n_2\triangledown\dots\triangledown n_i\triangledown m),\\ q_l-\G_s(n\triangledown m_1\triangledown m_2\triangledown \dots\triangledown m_j), p_k+q_l-\G_s(n_1\triangledown n_2\triangledown \dots\triangledown n_i\triangledown m_1\triangledown m_2\triangledown\dots\triangledown m_j)\}$, where $n_1+n_2+\dots +n_i=n-k$, $t_k=\Sigma_{i\in S_k}2^i$, $m_1+m_2+\dots m_j=m-l$ and $s_l=\Sigma_{j\in R_l}2^j$.}
\end{itemize}
\end{definition}

\begin{theorem}
Suppose $O_1,\dots,O_v$ are octal games, and there are natural numbers $N_1,\dots,N_v$ such that for each $i=1,\dots,v$, $G_s(n)\geq 0$ for all $n\in O_i$ and $n\leq N_i$. Then if $n_i\in O_i$ and $n_i\leq N_i$ for each $i=1,\dots,v$, $\G_s(n_1\triangledown\dots\triangledown n_v)=\Sigma_{i=1}^v\G_s(n_i)$.
\end{theorem}

\begin{proof}  I will prove this by induction on $n_1+\dots+n_j$ for some $j$.  The base case is clearly trivial since $\G_s(0\triangledown\dots\triangledown 0)=0$ regardless of how many $0$'s there are.

So for the inductive step assume that the result holds for all $n_1+\dots+n_j\leq K$ and I will choose and $n$ and $m$ such that $n+m=K+1$, and $G_s(n)$ and $G_s(m)\geq 0$.  The reason I only choose two games $n$ and $m$ is because it makes the proof easier and it will also be clear that the same argument can be extended to any number of games.

$\G_s(n\triangledown m)=\max_{k,i,l,j}\{p_k - \G_s(n_1\triangledown\dots\triangledown n_i\triangledown m), q_l-\G_s(n\triangledown m_1\triangledown \dots\triangledown m_j), p_k+q_l-\G_s(n_1\triangledown \dots\triangledown n_i\triangledown m_1\triangledown\dots\triangledown m_j)\}$, and since $n_1+\dots n_i+m$, $m_1+\dots+m_j+n$ and $n_1\dots n_i+m_1+\dots + m_j\leq k$, then by induction, $\max_{k,i,l,j}\{p_k - \G_s(n_1\triangledown\dots\triangledown n_i\triangledown m), q_l-\G_s(n\triangledown m_1\triangledown \dots\triangledown m_j), p_k+q_l-\G_s(n_1\triangledown \dots\triangledown n_i\triangledown m_1\triangledown\dots\triangledown m_j)\}=\\ \max\{p_k - \G_s(n_1\triangledown\dots\triangledown n_i)-\G_s(m), q_l-\G_s(n)-\G_s(m_1\triangledown \dots\triangledown m_j),\\ p_k+q_l-\G_s(n_1\triangledown \dots\triangledown n_i) -\G_s(m_1\triangledown\dots\triangledown m_j)\}=\\ \max\{\G_s(n)-\G_s(m), \G_s(m)-\G_s(n), \G_s(n)+\G_s(m)\}$.  

However since we know that both $\G_s(n)$ and $\G_s(m)\geq 0$, then\\ $\max\{\G_s(n)-\G_s(m), \G_s(m)-\G_s(n), \G_s(n)+\G_s(m)\}=\G_s(n)+\G_s(m)$, as previously stated it is clear that exactly the same argument can be used for any number games and so the theorem is proven.
\end{proof}

Note that this theorem will not hold if either $\G_s(n)$ or $\G_s(m)<0$, since in that case it might be better to move on $n$ or $m$ but not both $n$ and $m$, but this is still quite a strong result and tells us quite a lot about Nim variants played under the selective sum.  In the general case I make the following conjecture.

\begin{conjecture}
Let $O=(n_1n_2\dots n_t, p_1p_2\dots p_t)$ and $P=(m_1m_2\dots m_l, q_1q_2\dots q_l)$ be two finite taking-no-breaking octal games such that, there is at least one $n_s\neq 0$ or $1$, and if $n_i$ and $m_j=1,2$ or $3$ then $p_i=i$ and $q_j=j$, and $p_i=q_j=0$, otherwise, then for all $m$ there exists an $N$ such that;

$$\G_s(n+2k\triangledown m)=\G_s(n\triangledown m)$$

\noindent for all $n\geq N$ and $k$ is the largest entry in $O$ such that $n_k\neq 0,1$.
\end{conjecture}

As with Conjecture \ref{cs1} a proof of Conjecture~\ref{period}, i.e. the same conjecture under the disjunctive sum, would almost certainly yield a proof under the selective sum since the two operators are so similar in nature.  As with the disjunctive sum I have tried this for over 20 different examples and have yet to find a counter-example to this conjecture, so it is likely to be true.  But for exactly the same reasons as the disjunctive sum, this theorem is very difficult to prove.

An example of the conjecture is given in Table~\ref{ts1};

\begin{table}[htb]
\begin{center}
\begin{tabular}{c|ccccccccccccc}
$n\triangledown m$&0&1&2&3&4&5&6&7&8&9&10&11&12\\\hline
0&0&1&2&3&4&$-1$&2&3&0&1&2&1&0\\
1&1&2&3&4&5&2&3&4&1&2&3&2&1\\
2&2&3&4&5&6&3&4&5&2&3&4&3&2\\
3&3&4&5&6&7&4&5&6&3&4&5&4&3\\
4&2&3&4&5&6&1&0&1&2&3&2&1&2\\
5&1&2&3&4&5&0&1&2&3&4&3&2&3\\
6&0&1&2&3&4&1&2&3&4&3&2&3&4\\
7&1&2&3&4&5&2&3&4&5&4&3&4&5\\
8&2&3&4&5&6&3&2&3&4&5&4&3&4\\
9&3&4&5&6&7&2&1&2&3&4&3&2&3\\
10&2&3&4&5&6&1&0&1&2&3&2&3&1\\
11&1&2&3&4&5&0&1&2&3&4&3&2&3\\
12&0&1&2&3&4&1&2&3&4&5&4&3&2\\
\end{tabular}
\end{center}
\caption{$n\in (3311, 1234)\triangledown m\in (333, 123)$}\label{ts1}
\end{table}

What is interesting about this is that changing the operator does not appear to change the period, and in fact I make an even stronger conjecture;

\begin{conjecture}
Let $O=(n_1n_2\dots n_t, p_1p_2\dots p_t)$ and $P=(m_1m_2\dots m_l, q_1q_2\dots q_l)$ be two finite octal games, then if $\G_s(n\plus m)$ eventually has period $p$, $\G_s(n\triangledown m)$ also eventually has period $p$.
\end{conjecture}

So in other words what this conjecture says is that if these values are eventually periodic under the disjunctive sum, then not only are they eventually periodic under the selective sum, but they have the same period.

\section{The Sequential Join}

The final operator that I will look at is the sequential join.  This operator says that we arrange all of the components in order, and then play them in that order.

\begin{definition}  The sequential join of two games $G$ and $H$ is:
$$G\rhd H=\begin{cases}

&\{G^L\rhd H|G^S+H^S|G^R\rhd H\} \text{ if $G\neq \{.|G^S|.\}$}\\

&\{G^S\rhd H|G^S+H^S|G^S\rhd H\} \text{ otherwise}

\end{cases}$$

\end{definition}

In this case the long and short rules don't apply since the ending condition is given by the definition.  That is if $G^L=\emptyset$ and $H^L\neq\emptyset$ then $G^L\rhd H=\emptyset$. 

\begin{theorem}\label{seqid}
Scoring play games form a non-trivial monoid under the sequential join.
\end{theorem}

\begin{proof}  To prove this we first define a set  $I=\{i|i\rhd G\approx G\rhd i\approx G\hbox{ for all games }G\}$, and show that $I$ contains more than one element namely $\{.|0|.\}$.  So consider the game  $i=\{\{0|0|0\}|0|\{0|0|0\}\}$, as shown in the figure.
\begin{figure}[htb]
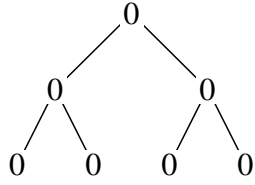

\begin{center}
\begin{graph}(3,2)

\roundnode{1}(0,0)\roundnode{2}(0.5,1)\roundnode{3}(1,0)\roundnode{4}(1.5,2)
\roundnode{5}(2,0)\roundnode{6}(2.5,1)\roundnode{7}(3,0)

\edge{1}{2}\edge{2}{3}\edge{4}{2}\edge{4}{6}\edge{5}{6}\edge{6}{7}

\nodetext{1}{0}\nodetext{2}{0}\nodetext{3}{0}\nodetext{4}{0}
\nodetext{5}{0}\nodetext{6}{0}\nodetext{7}{0}

\end{graph}
\end{center}
\caption{The game $i=\{\{0|0|0\}|0|\{0|0|0\}\}$}
\end{figure}

So first consider the game $i\rhd G$, if Left moves first on $i\rhd G$, then Right will move last on $i$, which means that Left will move first on $G$, and since the final score of $i$ is always $0$, then $(i\rhd G)^{SL}_F=G^{SL}_F$.  Similarly for the game $G\rhd i$, the players will simply play through $G$, and regardless of what happens the game $i$ cannot change the score of $G$, and therefore $(G\rhd i)^{SL}_F=G^{SL}_F$.

To show that the set is a monoid and not a group we need to demonstrate that not all games are invertible, so consider the game $Y=\{\{c|b|.\}|a|.\}\}$, and the game $G=\{e|d|f\}$.  If $Y$ is invertible this means that there exists a game $Y^{-1}$ such that $Y\rhd Y^{-1}\rhd G\approx G$ for all games $G$.  $G^{SR}_F=f$, however $(Y\rhd Y^{-1}\rhd G)^{SR}_F=a+a'+d\neq f$ and so the theorem is proven.
\end{proof}

\begin{theorem}
For any outcome classes $\mathcal{X}$, $\mathcal{Y}$ and $\mathcal{Z}$, there is a game $G\in \mathcal{X}$ and $H\in \mathcal{Y}$ such that $G\rhd H\in \mathcal{Z}$.
\end{theorem}

\begin{proof}  To prove this let $G=\{\{c|b|.\}|a|.\}$ and $H=\{H^L|d|H^R\}$, where $H^L$ and $H^R\neq \emptyset$, then $G^{SL}_F=a$, $G^{SR}_F=b$, $(G\rhd H)^{SL}_F=a+d$ and $(G\rhd H)^{SR}_F=b+d$.  Since $d$ is not dependent on $H^{SL}_F$ and $H^{SR}_F$, and can be any real number, then we can pick $a$, $b$ and $d$, so that $G$ and $H$ are in any outcome class and $G\rhd H$ is any outcome class.  Therefore the theorem is proven.  
\end{proof}

\subsection{Sequential Join Examples}

Again to compare the difference between the different operators I will use the same games from Sections~\ref{conjex} and \ref{sse}, and consider them played under the sequential join.

\begin{figure}[htb]
\begin{center}
\begin{pspicture}(9,5)
\put(0,2){$BTF\triangledown TFB=$}\put(4,4){\cgtree{
  {BTF(0)}(|{FTB(-1)}({FBT(-1)}|))
  }}
   \put(9,4){\cgtree{
  {TFB(0)}({BTF(1)}(|{FBT(1)})|)}}
\put(6.5,2){$\rhd$}
\end{pspicture}
\end{center}
\caption{$BTF\rhd TFB=\{.|0|\{-1|-1|.\}\}\rhd\{\{.|1|1\}|0|.\}$}
\end{figure}

The game $BTF\rhd TFB=\{.|0|\{-1|-1|.\}\}\rhd\{\{.|1|1\}|0|.\}=\{.|0|\{\{.|0|0\}|-1.\}|-1|.\}$.  Which again gives a different result for the same played under the previous two operators.  It is also worth noting that the sequential join is not commutative and therefore if we played the same games in a different order we would get a different answer.

The game in Figure~\ref{h.ex2} has value $\{\{.|1|\{.|0|-1\}\}|0|\{\{.|0|1\}|-1|-3\}\}$ under the sequential join.  Again the difference is to due to the two red edges being played sequentially after the removal of the blue edge.

\begin{figure}[htb]
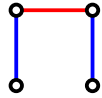

\begin{center}
\begin{graph}(1,1)

\graphnodecolour{1}\graphlinewidth{0.05}

\roundnode{4}(0,0)
\roundnode{5}(0,1)\roundnode{6}(1,0)\roundnode{7}(1,1)

\edge{4}{5}[\graphlinecolour(0,0,1)]
\edge{5}{7}[\graphlinecolour(1,0,0)]\edge{6}{7}[\graphlinecolour(0,0,1)]

\end{graph}
\end{center}
\caption{The game $\{\{.|1|-1,\{.|0|-1\}|0|\{\{.|0|1\}|-1|-3\}\}$}\label{h.ex2}
\end{figure}

\subsection{Impartial Games}

\begin{theorem}
Impartial games for a non-trivial monoid under the sequential join.
\end{theorem}

\begin{proof}  From the proof of Theorem~\ref{seqid} we know that there is a non-trivial identity set, so to prove this we simply need to show that there is a game $G$ that is not invertible.  So consider the game $G=\{1,\{0|0|0\}|0|\{0|0|0\},-1\}$.  Let $Y$ be the inverse of $G$, then this implies that $G\rhd Y\rhd P\approx P$ for all impartial games $P$.

So let $P=\{.|0|.\}$, and consider the game $G\rhd Y\rhd P$.  If Left moves first and moves to the game $1\rhd Y\rhd P$, then his implies that $-1$ is one of the Right options of $G$, since if Right moves to $-1$ on $Y$ then Left will move first on $P$ and $G\rhd Y$ will not change the final score of $P$.  But $Y$ is impartial, so this implies that $1$ is a Left option of $Y$.  So therefore if Left moves to the game $\{0|0|0\}\rhd Y\rhd P$, then this means that Right must move to the game $0\rhd Y\rhd P$, and Left will move first on $Y$, and Left can choose the option $1$ and hence win $G\rhd Y\rhd P$, i.e. $G\rhd Y\rhd P\not \approx P$ which is a contradiction.

So this means that $G$ is not invertible, and therefore the set of impartial games form a non-trivial monoid under the sequential join and the theorem is proven.
\end{proof}

\subsection{Sprague-Grundy Theory}

When consider the sequential join it does not really make sense to look at taking and breaking games, because once you break the heap into two or more smaller heaps we have to define the order that we play the two new heaps in.  Therefore we concentrate below on taking-no-breaking games.

\begin{definition}
Let $n\in O=(t_1t_2\dots t_f, p_1,\dots p_f)$ and $m\in P=(s_1s_2\dots s_e, q_1\dots q_e)$, be two taking-no-breaking games;

\begin{itemize}
\item{$\G_s(0)=0$.}
\item{$\G_s(n\rhd m)=\begin{cases}
&\max\{p_k-\G_s(n-k\rhd m)\} \text{ if $n\neq 0$}\\
&\max\{q_l-\G_s(n\rhd m-l)\} \text{ otherwise}
\end{cases}$}
\end{itemize}
\end{definition}

There is not really a lot to say about this operator, other than to make the following conjecture;

\begin{conjecture}
Let $O=(n_1n_2\dots n_t, p_1p_2\dots p_t)$ and $P=(m_1m_2\dots m_l, q_1q_2\dots q_l)$ be two finite octal games, then if $\G_s(n\plus m)$ eventually has period $p$, $\G_s(n\rhd m)$ also eventually has period $p$.
\end{conjecture}

This conjecture seems quite a reasonable one due to the nature of the operator.  By playing the heaps in order, it means that $m$ cannot change the period of $n$.  However since it is very hard to even prove that $\G_s(n+p)=\G_s(n)$ for all $n$ large enough, a proof of this conjecture will also be very difficult.

An example of octal games played under the sequential join is given in Table~\ref{sj1};

\begin{table}[htb]
\begin{center}
\begin{tabular}{c|rrrrrrrrrrrrr}
$n\rhd m$&0&1&2&3&4&5&6&7&8&9&10&11&12\\\hline
0&0&1&0&3&2&3&0&1&2&3&2&1&0\\
1&1&0&1&2&3&2&1&0&1&2&3&2&1\\
2&0&1&0&3&2&3&0&1&2&3&2&1&0\\
3&1&0&1&2&3&2&1&0&1&2&3&2&1\\
4&4&$-3$&4&5&6&$-1$&$-2$&3&4&5&0&$-1$&2\\
5&5&$-4$&5&6&7&$-2$&$-3$&4&5&6&$-1$&$-2$&3\\
6&4&$-3$&4&5&6&$-1$&$-2$&3&4&5&0&$-1$&2\\
7&5&$-4$&5&6&7&$-2$&$-3$&4&5&6&$-1$&$-2$&3\\
8&4&$-3$&4&5&6&$-1$&$-2$&3&4&5&0&$-1$&2\\
9&1&0&1&2&3&2&1&0&1&2&3&2&1\\
10&0&1&0&3&2&3&0&1&2&3&2&1&0\\
11&1&0&1&2&3&2&1&0&1&2&3&2&1\\
12&0&1&0&3&2&3&0&1&2&3&2&1&0\\
\end{tabular}
\end{center}
\caption{$n\in (123, 123)\rhd m\in (30033, 10045)$}\label{sj1}
\end{table}

\section{Summary}

In this chapter I looked at the general structure of scoring play games under three different operators, namely the conjunctive sum, selective sum and sequential join.  

I showed that under these operators scoring play games are still quite ``wild'', namely they is very little comparibility between games and they do not form a group.

I also looked at impartial games under each of these different operators and in particular the game of Nim and the function $\G_s(n)$.  What is interesting is that it appears that the function still has the same periodicity regardless of which operator is used.  I think a proof of that would be very interesting, and quite a remarkable and surprising result.
\chapter{Games and Complexity}

\section{Complexity Theory}

Complexity theory is primarily focused on decision problems, that is problems where for any given instance the outcome is either positive or negative.  Clearly this is directly related to combinatorial game theory, since for any given game position we want to know if a given player has a winning strategy or not.  In this case we want to know if there is a way to find out if a given player has a winning strategy quickly or not.  In other words how is the size of the input related to the number of steps that a computer or a person requires to determine the answer.

A Turing machine was first conceived by Alan Turing in 1936 \cite{AT}, it is not designed to be practical computing technology but a thought experiment representing a computer machine.  A Turing machine is the backbone for determining how hard, or complex a given decision problem is to solve, and it is defined as follows.

\begin{definition}\cite{AT}
A one tape Turing machine is a 7-tuple $M=<\mathcal{Q},\Gamma,b,\Sigma,\delta,q_0,F>$ where;

\begin{itemize}
\item{$\mathcal{Q}$ is a finite non-empty set of states.}
\item{$\Gamma$ is a finite, non-empty set of the tape alphabet or symbols.}
\item{$b\in \Gamma$ is the blank symbol.}
\item{$\Sigma\subseteq \Gamma\backslash\{b\}$ is the set of input symbols.}
\item{$q_0\in \mathcal{Q}$ is the initial state.}
\item{$F\subseteq \mathcal{Q}$ is the set of final or accepting states.}
\item{$\delta:\mathcal{Q}\backslash F\times\Gamma\rightarrow\mathcal{Q}\times\Gamma\times\{L,R\}$ is a transition function, where $L$ is a left shift and $R$ is a right shift.}
\end{itemize}
\end{definition}

\begin{definition}\cite{GJ}
A deterministic polynomial-time Turing machine is a deterministic Turing machine $M$ that satisfies the following conditions:

\begin{enumerate}
\item{$M$ halts on all input $w$.}
\item{There exists $k\in \mathbb{N}$ such that $T_M(n)\in O(n^k)$, where $T_m(n)=\hbox{max}\{t_M(w)|w\in \Sigma^*,|w|=n\}$, and $t_m(w)=\hbox{ the number of steps } M\hbox{ takes to halt on input }w$.}
\end{enumerate}
\end{definition}

Before we can define the set $NP$, we first need the set $P$, that is the set of all decision problems that can be solved in polynomial time or faster.  Formally it is defined as follows.

\begin{definition}\cite{GJ}$P=\{L|L=L(M)\hbox{ for some deterministic Turing machine }M\}$,\\ where $L(M)=\{w\in\Sigma^*|M\hbox{ accepts }w\}$.
\end{definition}

\begin{definition}\cite{GJ} $L\in NP$ if $\exists$ a binary relation $R\subset\Sigma^*\times\Sigma^*$ and a $k>0$ such that:

\begin{enumerate}
\item{If $x\in\Sigma^*$ then $x\in L$ if and only if $\exists y\in \Sigma^*$ such that $(x,y)\in R$ and $|y|\in O(|x|^k)$. }
\item{The language $L_k=\{x b y|(x,y)\in R\}$ over $\Sigma \cup \{b\}$ is decidable by a Turing machine in polynomial time.}
\end{enumerate}
\end{definition}

In other words, a decision problem $L\in NP$ if the solution is verifiable in polynomial time, which means that if I were to simply guess a solution to the decision problem $L$, then we can define a polynomial time Turing machine that will verify if my guess is correct or not.  Which brings us to the definition of $NP$-complete.

\begin{definition}\cite{GJ}
A decision problem $L'$ is polynomial time reducible to a decision problem $L$ (written as $L'\leq_p L$), if the following hold:
\begin{enumerate}
\item{$\exists$ $f:\Sigma^*\rightarrow \Sigma^*$ such that $\forall w\in\Sigma^*$, $w\in L'$ if and only if $f(w)\in L$.}
\item{$\exists$ a polynomial time Turing machine that halts with $f(w)$ steps on its tape on any input $w$.}
\end{enumerate}
\end{definition}

\begin{definition}\cite{GJ} $L$ is $NP$-complete if the following hold:
\begin{enumerate}
\item{$L\in NP$.}
\item{For all $L'\in NP$, $L'\leq_p L$.}
\end{enumerate}
\end{definition}

\begin{definition}\cite{GJ}
$L$ is $NP$-hard is $L'$ if $NP$-complete and $L'\leq_p L$.
\end{definition}

Note that for a problem $L$ to be $NP$-hard it does not have to be a decision problem.

\section{Normal Play Red-Blue Hackenbush}

In Winning Ways \cite{WW}, the authors give a proof that in normal play the game Red-Blue Hackenbush is $NP$-hard, I will give an outline of the proof here.

\begin{definition}
Redwood furniture, is a Hackenbush position where all blue edges are grounded, there are no grounded red edges, and every red edge is connected to exactly one blue edge.
\end{definition}

\noindent PROBLEM: \textbf{REDWOOD FURNITURE HACKENBUSH}

\noindent INSTANCE: A position of redwood furniture $G$, integer $m$.

\noindent QUESTION: Is $G\leq \frac{1}{2^m}$?

\begin{theorem}\cite{WW}
Redwood furniture is $NP$-hard.
\end{theorem}

\begin{proof}  The proof of this relies on the result that a ``redwood tree'' $T$ has value $\frac{1}{2}$.  A redwood tree is a piece of redwood furniture where there are no cycles in the red edges.  We will not be proving that these games have value $\frac{1}{2}$, but instead refer the reader to \cite{WW} pages 211-217.

The authors then show that any piece of redwood furniture $F$ has value $\frac{1}{2^{m+1}}$, and do a reduction from a minimal spanning tree problem in graph theory to obtain the value of $m$.  Thus demonstrating that the problem is $NP$-hard.
\end{proof}

Since this problem is a subset of the problem of determining the outcome class of a general position of Red-Blue Hackenbush, we can conclude that Red-Blue Hackenbush is also $NP$-hard.

\section{Mis\`ere Hackenbush}

\subsection{Red-Blue Mis\`ere Hackenbush}

One might expect when we look at Red-Blue Hackenbush under mis\`ere rules that the game is also $NP$-hard, however that is not the case, as was shown by Stewart in \cite{FS}, below.

\begin{theorem}\label{hackenbush}\cite{FS}  Let $G$ be a game of Red-Blue
mis\`ere Hackenbush, and let $B$ and $R$ be the number of grounded
blue and red edges respectively, then the outcome of $G$ can be
determined by the following formula:

$$G\in\begin{cases}
&\lef\text{, if $B>R$}\\
&\ri\text{, if $R>B$}\\
&\n\text{, if $B=R$}
\end{cases}$$

\end{theorem}

The proof can be found in \cite{FS}.

\subsection{Red-Blue-Green Mis\`ere Hackenbush}

We now come to the only new result in this chapter, showing that Red-Blue-Green Mis\`ere Hackenbush is NP-hard.  Red-Blue-Green Hackenbush is played identically to Red-Blue Hackenbush, apart from the addition of green edges that may be removed by either player.

\noindent PROBLEM: \textbf{RED-BLUE-GREEN MIS\`ERE HACKENBUSH}

\noindent INSTANCE: A position of Red-Blue-Green Mis\`ere Hackenbush
$G$.

\noindent QUESTION: What is the outcome of $G$?

\begin{theorem}
Red-Blue-Green Mis\`ere Hackenbush is NP-hard.
\end{theorem}

\begin{proof}
To prove this we will a do a transformation from Red-Blue Hackenbush under normal play rules.
First we note two things, as previously stated, it is known that determining the outcome of a general
position of normal play Red-Blue Hackenbush is NP-hard.  It is also known that we can think of the
ground in Hackenbush as being a single vertex, which is drawn as a ground with separate vertices
for clarity in diagrams.  With this in mind we will make our transformation.

The transformation will be as follows, start with a general Red-Blue Hackenbush position $G$.  Next
take the same position and replace the ground, and all the vertices that are on the ground with a single vertex
and call this game $G'$.  Lastly attach $G'$ to a single grounded green edge, and call this game $G_m$.  This process
is illustrated in Figure~\ref{transform}.

\begin{figure}[htb]
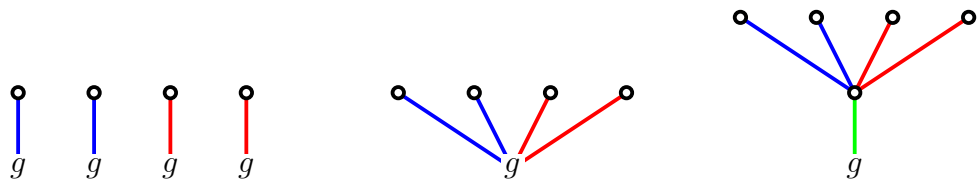

\begin{center}
\begin{graph}(12.5,2)

\graphnodecolour{1}\graphlinewidth{0.05}
\roundnode{1}(0,0)\roundnode{2}(0,1)\roundnode{3}(1,0)\roundnode{4}(1,1)
\roundnode{5}(2,0)\roundnode{6}(2,1)\roundnode{7}(3,0)\roundnode{8}(3,1)
\roundnode{9}(6.5,0)\roundnode{10}(11,0)\roundnode{11}(11,1)\roundnode{12}(7,1)
\roundnode{13}(6,1)\roundnode{14}(5,1)\roundnode{15}(8,1)\roundnode{16}(9.5,2)\roundnode{17}(10.5,2)\roundnode{18}(11.5,2)\roundnode{19}(12.5,2)

\edge{1}{2}[\graphlinecolour(0,0,1)]
\edge{3}{4}[\graphlinecolour(0,0,1)]
\edge{9}{13}[\graphlinecolour(0,0,1)]
\edge{9}{14}[\graphlinecolour(0,0,1)]\edge{11}{16}[\graphlinecolour(0,0,1)]\edge{11}{17}[\graphlinecolour(0,0,1)]

\edge{5}{6}[\graphlinecolour(1,0,0)]
\edge{7}{8}[\graphlinecolour(1,0,0)]
\edge{9}{12}[\graphlinecolour(1,0,0)]
\edge{9}{15}[\graphlinecolour(1,0,0)]
\edge{11}{18}[\graphlinecolour(1,0,0)]\edge{11}{19}[\graphlinecolour(1,0,0)]

\edge{10}{11}[\graphlinecolour(0,1,0)]

\nodetext{1}{$g$}\nodetext{3}{$g$}\nodetext{5}{$g$}\nodetext{7}{$g$}\nodetext{9}{$g$}\nodetext{10}{$g$}

\end{graph}
\end{center}
\caption{Transformation of $G$ to $G_m$.}\label{transform}
\end{figure}

If we are playing $G_m$ under mis\`ere rules, then neither player will want to cut the single green edge, since doing
so will remove every edge in the game, and thus the next player will be unable to move and therefore win under mis\`ere rules.
So both players will want to move last on the graph $G'$ that is attached to the single green edge, thus forcing your opponent
to remove the green edge, which will result in you winning the game.  In other words, whoever wins $G'$ under normal play rules, will
also win $G_m$ under mis\`ere play rules, and since determining the outcome of $G'$ is NP-hard, determing the outcome of $G_m$ is also NP-hard.
So the theorem is proven.
\end{proof}

\section{Summary}

In this chapter I simply proved that Red-Blue-Green Mis\`ere Hackenbush is NP-hard.  This is a new and interesting result, as it is the first complexity theory result for a mis\`ere game.
\chapter{Applications}

\section{Possible Applications}

This section is dedicated to explaining how scoring play theory can be applied to analyse ``real world'' games.  That is games that are played regularly by people, as opposed to games such as Hackenbush.  The two in particular which I am going to look at are Go and Sowing games.

\subsection{Mathematical Go Comparison}\label{Go}

Go is one of the oldest games in existence that is still played today in its original form.  Nobody knows precisely how old it is, or indeed who its creator was.  What is known is that it was invented in China at least 3000 years ago.  Its name in Chinese is Weiqi which literally translates as wrapping or surrounding chess. 

The game travelled to Japan around 1000 years ago where it became wildly popular and was played by Samurai warriors as a way of sharpening their minds to make them better fighters.  People in the western world often mistakingly cite Go as an ``Asian'' or ``Japanese'' game, because the name Go comes from the Japanese pronunciation of the game, which is Igo.  However one thing to note is that the characters used for the game in China, Japan and Korea are identical they just have different pronunciations, the Japanese use the Chinese name for the game, i.e. surrounding or wrapping chess.

Today this game is massively popular and enjoyed by people all over the world, particularly in Asia, where professionals can earn huge amounts of prize money in the numerous competitions that are held.  As a result the game has attracted a massive amount of mathematical interest, and in 1994 Berlekamp and Wolfe wrote a book called Mathematical Go \cite{MG}.

The idea behind this book was to take the theory of combinatorial games that was developed by Berlekamp, Conway and Guy \cite{WW, MG} and use it to analyse very specific types of Go endgame scenarios.  The idea is to try and determine who wins the last point, and therefore wins the game.

This idea was rather successful and has been used in actual game play, as well as further research to analyse more complicated Go positions.  However our understanding of the game of Go still remains highly limited, with very little real progress being made.

Scoring play theory however can offer us some hope.  Elwyn Berlekamp and David Wolfe compared Go to \emph{normal play} games, and used \emph{normal play} theory in their analysis.  By using scoring play theory and putting Go into the realm of scoring play games, we get different results.

The basic rules of the game are as follows;

\begin{enumerate}
\item{There are two players, Left plays with Black stones, and Right plays with White stones.}
\item{The game is played on a $n\times m$ grid, where $n$ and $m$ are the number of horizontal and vertical lines respectively.}
\item{Players take turns to place stones on the intersections of the lines.}
\item{Players can capture their opponent's stones or intersections of the board by placing stones on all adjacent intersections.}
\item{Under Chinese rules players get one point for every stone that they have on the board and every vacant point between those stones. Vacant points between both sides' stones are shared equally.}
\item{Under Japanese rules players get one point for every stone and intersection that is captured.}
\item{The game ends either when all areas of the board have been captured, a player concedes or a player makes two consecutive passes.}
\item{The player who gets the most points wins.}
\end{enumerate}

While I have given the rules for Japanese and Chinese scoring, and these rules do make the scoring quite different, in this section I will be concentrating purely on the Japanese scoring system, since this is the one that is used in the book Mathematical Go \cite{MG}.

The game has other rules and variations depending upon if you are using Japanese rules, Chinese rules, American rules, Ancient rules and so on.  For example under Japanese rules you may not place a stone in an area where it would be automatically captured, but in Chinese rules these moves are permitted.

There is also a special rule called a ``ko''.  All variations of the game use the rule that you \emph{may} place a stone in an area where it would normally be automatically captured if you are capturing one of the surrounding stones.  This rule can lead to situations involving loops or kos as they are called in Japanese.  I discuss this particular rule further in Section~\ref{app}.

\subsubsection{Zero}\label{zero}

The first thing I need to talk about is the idea of zero.  In normal play combinatorial game theory $0=\{.|.\}$ and is equivalent to all $\pre$ positions.  By using this approach on Go it is possible to dismiss a lot of positions because they are equivalent to zero and therefore make no difference to the overall game.

\begin{figure}[htb]
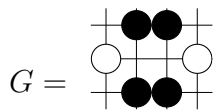

\begin{center}
\begin{tabular}{r@{\hspace{0.2in}}r@{}}
$G=$&\begin{psgopartialboard*}{(7,7)(10,9)}
\stone{white}{g}{8}
\stone{black}{h}{7}
\stone{black}{h}{9}
\stone{black}{j}{7}
\stone{black}{j}{9}
\stone{white}{k}{8}
\end{psgopartialboard*}\\
\end{tabular}
\end{center}
\caption{A game with normal play value $0$}\label{z}
\end{figure}

A good example of this is the game $G$ in Figure~\ref{z}.  Using the theory from Mathematical Go this position can be dismissed because it makes no difference to the overall game since it does not affect who moves last.

As an alternative using scoring play theory \emph{if} this game really makes no difference then we can understand \emph{why} it makes no difference.  As we know from Theorem~\ref{identity} the above game is not identical to zero and therefore not equal to zero.  This means that there may be Go positions where this game \emph{can} affect the outcome of an overall game.

By using scoring play theory one could perhaps show that no such Go position can exist, and by doing so we can really understand how positions where neither player can gain any points affect the game, if they do at all.  This can lead to a much greater understanding of the game.

\subsubsection{Dominated Options}

Another difference between scoring play theory and Mathematical Go theory are dominated options.  Using normal play theory the authors were able to dismiss many positions because under normal play theory they are dominated options.  However as an alternative using scoring play theory the same positions are not dominated, and as before if indeed these are bad moves when you play Go, scoring play theory can help us understand \emph{why}.

\begin{figure}[htb]
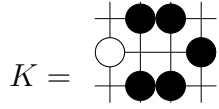

\begin{center}
\begin{tabular}{r@{\hspace{0.2in}}r@{}}
$K=$&\begin{psgopartialboard*}{(7,7)(10,9)}
\stone{white}{g}{8}
\stone{black}{h}{7}
\stone{black}{h}{9}
\stone{black}{j}{7}
\stone{black}{j}{9}
\stone{black}{k}{8}
\end{psgopartialboard*}\\
\end{tabular}
\end{center}
\caption{The Go Position $K$}\label{corridor}
\end{figure}

To demonstrate this consider the game $K$ shown in Figure~\ref{corridor}.  The Left options of $K$ are $G$ and $H$ and are shown in Figures~\ref{g} and \ref{h}.  According to the book Mathematical Go $G=*$ and $H=1$, therefore $H\geq G$.  This means that in the book Mathematical Go they do not consider the option $G$ from the game shown in Figure~\ref{corridor}.

Using scoring play theory however we can see that $H\not\geq G$.  The game $G=\{0|0|0\}$ and $H=\{0|1|.\}$.  So consider the game $X=\{1|-2|.\}$, Left wins $G\plus X$ moving second, but loses moving second on $H\plus X$, which means that $H\not\geq G$.  In other words we have to consider the option $G$ from $K$.

\begin{figure}[htb]
\begin{pspicture}(5,5)
\put(0,2){$G=$}\put(2,1.5){\begin{psgopartialboard*}{(7,7)(10,9)}
\stone{white}{g}{8}
\stone{black}{h}{7}
\stone{black}{h}{9}
\stone{black}{j}{7}
\stone{black}{j}{9}
\stone{black}{j}{8}
\stone{black}{k}{8}
\end{psgopartialboard*}}\put(4.5,2){$=$}\put(10,4){\cgtree{{\begin{psgopartialboard*}{(7,7)(10,9)}
\stone{white}{g}{8}
\stone{black}{h}{7}
\stone{black}{h}{9}
\stone{black}{j}{7}
\stone{black}{j}{9}
\stone{black}{j}{8}
\stone{black}{k}{8}
\end{psgopartialboard*}}(++{\begin{psgopartialboard*}{(7,7)(10,9)}
\stone{white}{g}{8}
\stone{black}{h}{7}
\stone{black}{h}{9}
\stone{black}{j}{7}
\stone{black}{j}{9}
\stone{black}{j}{8}
\stone{black}{h}{8}
\stone{black}{k}{8}
\end{psgopartialboard*}}|++{\begin{psgopartialboard*}{(7,7)(10,9)}
\stone{white}{g}{8}
\stone{black}{h}{7}
\stone{black}{h}{9}
\stone{black}{j}{7}
\stone{black}{j}{9}
\stone{black}{j}{8}
\stone{white}{h}{8}
\stone{black}{k}{8}
\end{psgopartialboard*}})
}}
\end{pspicture}
\caption{$G=\{0|0|0\}$}\label{g}
\end{figure}

\begin{figure}[htb]
\begin{pspicture}(5,5)
\put(0,2){$H=$}\put(2,1.5){\begin{psgopartialboard*}{(7,7)(10,9)}
\stone{white}{g}{8}
\stone{black}{h}{7}
\stone{black}{h}{9}
\stone{black}{j}{7}
\stone{black}{j}{9}
\stone{black}{h}{8}
\stone{black}{k}{8}
\end{psgopartialboard*}}\put(4.5,2){$=$}\put(10,4){\cgtree{{\begin{psgopartialboard*}{(7,7)(10,9)}
\stone{white}{g}{8}
\stone{black}{h}{7}
\stone{black}{h}{9}
\stone{black}{j}{7}
\stone{black}{j}{9}
\stone{black}{h}{8}
\stone{black}{k}{8}
\end{psgopartialboard*}}(++{\begin{psgopartialboard*}{(7,7)(10,9)}
\stone{white}{g}{8}
\stone{black}{h}{7}
\stone{black}{h}{9}
\stone{black}{j}{7}
\stone{black}{j}{9}
\stone{black}{h}{8}
\stone{black}{j}{8}
\stone{black}{k}{8}
\end{psgopartialboard*}}|)
}}
\end{pspicture}
\caption{$H=\{0|1|.\}$}\label{h}
\end{figure}

As with Section~\ref{zero}, this difference is important because while it is most likely the case that the game $H$ is always better than $G$ \emph{when you play Go}, scoring play theory can help us understand \emph{why}.

Again the key is understand the true nature and structure of the game of Go, and scoring play theory provides us with the necessary tools to gain that understanding.

\subsection{Application to Go in General}\label{app}

As I've shown scoring play theory differs from the theory presented in Mathematical Go in some key places, however the theory can be used to study Go much more generally.

First of all recall that $G\geq H$ if $H\plus X\in O$ implies that $G\plus X\in O$, where $O=L_\geq, R_\geq, L_>$ or $R_>$ for all games $X$.  The problem with this definition is that not every scoring play game can be represented by a Go position.  So a lot of the time we might have a situation where $G$ is a better option than $H$ when playing the game of Go, but $G\not\geq H$ in general.

The way to overcome a problem like this is to restrict our set to only Go positions.  To do that we have to know precisely what properties a Go position has, and define them clearly and mathematically.  There may of course be big differences between different rulesets e.g. Chinese or Japanese rules, however in this way one can prove conclusively that an option $G$ is indeed always better than an option $H$. 

\subsubsection{Ko}\label{ko}

Kos are an area that is very complicated, but extremely important when studying Go because they are a major part of the game. I'd like to demonstrate how scoring play theory can be used to study these positions as well.

A Ko is a position which repeats itself after two consecutive moves by the same player, i.e. Left plays, then Right plays, then Left plays again and we have returned to the original position.  An example is the game shown in Figure~\ref{ko}

\begin{figure}[htb]
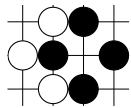

\begin{center}
\begin{psgopartialboard*}{(7,7)(10,9)}
\stone{white}{g}{8}
\stone{white}{h}{7}
\stone{black}{h}{8}
\stone{white}{h}{9}
\stone{black}{j}{7}
\stone{black}{j}{9}
\stone{black}{k}{8}
\end{psgopartialboard*}
\end{center}
\caption{An example of a ko}\label{ko}
\end{figure}

Interestingly in Chinese the word for ko is either qie or hukou, the first meaning a disaster, the second meaning a tiger's mouth.  The reason why a position of this nature is so undesirable is because if white takes the black stone, black can immediately retake the white stone, and this could repeat, potentially, forever.

However under Go rules if white captures the black stone, then black may not recapture the white stone until his next turn.  This rule is an attempt to prevent a potential cycle that may never end.

\begin{figure}[htb]
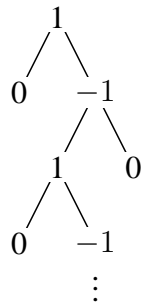

\begin{center}
\begin{graph}(1.5,3.5)

\roundnode{1}(0,0.5)\roundnode{2}(1,0.5)\roundnode{3}(0.5,1.5)\roundnode{4}(1.5,1.5)
\roundnode{5}(1,2.5)\roundnode{6}(0,2.5)\roundnode{7}(0.5,3.5)

\edge{1}{3}\edge{2}{3}\edge{3}{5}\edge{4}{5}\edge{5}{7}\edge{6}{7}

\nodetext{1}{0}\nodetext{2}{$-1$}\nodetext{3}{1}\nodetext{4}{0}\nodetext{5}{$-1$}\nodetext{6}{0}\nodetext{7}{1}

\freetext(1,0){$\vdots$}

\end{graph}
\end{center}
\caption{The game tree of a ko}\label{gtk}
\end{figure}

So how can scoring play theory be used to analyse such a position?  With scoring play theory the best thing to do is not to draw in ``loops'' on the game tree, but rather in the form shown in Figure~\ref{gtk}, which corresponds to the game given in Figure~\ref{ko}.  The reason for this is that looping can change the score through repeated cycles so it makes more sense to think of ``loopy'' games as a game tree of infinite depth.

With this approach using scoring play theory it is a matter of determining for any given position when the best time to move to 0 is, or as they say in Go, ``fill the ko''.  Using our braces and slashes notation this game is written as $\{0|1|\{0|1|\dots\}\}$.

So I have demonstrated how scoring play theory can be applied to Go, and of course there are many other possible ways one can use this theory to analyse the game of Go, and increase our understanding of the game and the subtle strategies involved in a rigorous mathematical way.

\subsection{Sowing Games}

Sowing games are among the oldest known combinatorial games in existence.  They originated in Africa over 7000 years ago.  There are hundreds of different variants of these games, but in this section I will only be describing one, called ``Mancala''.  The rules are given below;

\begin{enumerate}
\item{Each player has six cups, which contain beans, and one pot for keeping gathered beans, called a kala.}
\item{On a player's turn he removes all of the beans from one of his cups, and sows one into each of his cups, then his kala, then his opponent's cups.}
\item{A player may only sow beans into his own kala.}
\item{The game ends when one player is unable to move, i.e. all of his cups are empty, and the winner is the player who collected the most beans into his kala.}
\end{enumerate}

What is interesting is that despite the wide variety of games, there has been little interest in these games mathematically.  In fact I know of only one paper written about them and that is a paper by Jeff Erickson published in 1996 \cite{JE}.

Since these are clearly scoring games the theory can again be directly applied to analyse them.  What is noteworthy about these games is that unlike many combinatorial games they do not naturally split up into sums of separate components.  I feel that scoring play combinatorial game theory can help us to understand these games a lot more than using normal play theory.

The reason for this is similar to the reasons I gave for why I feel that scoring play theory works better than normal play theory when analysing Go.  We want to know why certain moves are better than others, and since sowing games are generally scoring play games, scoring play theory is the natural choice to help us gain that understanding.

\chapter{Conclusion}

\section{Areas for Further Research}

In this section I will be looking at what direction we can take the theory of scoring play games when researching it further.  This thesis has really laid the foundations, but the really exciting work is ahead.  Here are some of my ideas which I think would be good ways to take the theory of scoring play games forward.

\subsection{Scoring Play Hackenbush}

I discussed scoring play Hackenbush a little in chapter 2.  I used it to demonstrate different scoring games in each of the different outcome classes.  This game is important because in both normal and mis\`ere play it has been used to highlight and analyse many different aspects of the theory concerning both methods of play.  

Scoring play Hackenbush has many distinguishing features that make it more interesting to analyse than both normal and mis\`ere play Hackenbush.  As I demonstrated in chapter 2, the red-blue version of the game can belong to any outcome class, unlike normal and mis\`ere play Hackenbush where, if the game is played with red and blue edges only the game can only belong to $\lef$, $\ri$ or $\pre$ and $\lef$, $\ri$ or $\n$ under normal and mis\`ere play respectively.

The other interesting aspect of the game is that changing the way the game is scored changes the analysis of the game.  For example consider the following variations.

\begin{enumerate}
\item{A player gains points for all edges that they disconnect.}
\item{A player only gains points for disconnecting edges of his own colour.}
\item{A player gains points for disconnecting edges of his own colour and loses points for disconnecting edges of his opponent's colour.}
\end{enumerate}

These three variations change the situation dramatically, for example consider the following simple position.

\begin{figure}[htb]
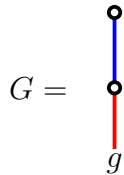

\begin{center}
\begin{graph}(1,2)

\graphnodecolour{1}\graphlinewidth{0.05}

\roundnode{1}(1,0)\roundnode{2}(1,1)\roundnode{3}(1,2)

\edge{1}{2}[\graphlinecolour(1,0,0)]\edge{2}{3}[\graphlinecolour(0,0,1)]

\freetext(0,1){$G=$}

\nodetext{1}{$g$}

\end{graph}
\end{center}
\caption{A simple Hackenbush position.}
\end{figure}

Under the first variation the game $G=\{\{.|1|0\}|0|-2\}$, under the second variation the game $G=\{\{.|1|0\}|0|-1\}$ and under the third variation the game $G=\{\{.|1|0\}|0|0\}$.  Under the first two variations $G\in \ri$, but under the third $G\in\ti$, which is a big difference in combinatorial game theory.

As I proved in chapter 5, mis\`ere play Hackenbush is $NP$-hard, so I have the following conjecture about scoring play Hackenbush.

\begin{question}
If a player gains points for all edges that they disconnect, then what is the complexity of scoring play Hackenbush?
\end{question}

\subsection{Zugzwang Positions}

A zugzwang position is a position in a combinatorial game where moving is worse than not moving.  These positions were first studied by Li and Yen in \cite{LY}, in relation to normal play combinatorial games.  For scoring play games a zugzwang position would be one where you start off with more points that your opponent, but if you move you end the game with fewer points than your opponent.

An example of a scoring play game with zugzwang positions would be the child's game of dots and boxes.  This is a game that is played on an $n\times m$ grid of dots, and each move you make you must draw a horizontal or vertical line between two of the dots.  You gain points by completing a box, if you complete a box you get an extra turn and at the end of the game the player who completes the most boxes wins.

This game has been well studied by combinatorial game theorists, and has a full solution since it is another example of a scoring game where the last player to move also happens to have the most points.  The reason this game has zugzwang positions is because if you complete a square you are forced to make another move, and this may mean that your opponent can now complete many squares.  In other words if you didn't have to make that move then your opponent would make it and you would win, but because you are forced to make the move you lose, which is zugzwang.

Of course there are many scoring play games with zugzwang positions, so this is an area that deserves further research.

\subsection{Passing}

Passing is not really something that is applicable to either normal play or mis\`ere play.  A little bit of work has been done on passing in relation to loopy games, but in general passing does not make sense in games where the winner is determined by who moves last.

However for scoring play games passing certainly makes sense since the winner is determined by the score.  Indeed in the game of Go passing is permitted.  So it would be nice to know how scoring play games behave if a player is allowed to pass if he only has bad moves available.

There are three basic ways  to pass;

\begin{enumerate}
\item{Optional passing.  Where the player can choose when he wishes to pass.}
\item{Limited passing.  Where players are given $n$ passes and may only use them when they do not have a move.}
\item{Unlimited passing.  Where players pass whenever they do not have a move, and the game ends when both players cannot move on all components.}
\end{enumerate}

I believe that further research into this aspect of scoring play games may yield some useful results.

\subsection{More Than Two Players}

Research into combinatorial game theory in the past has almost exclusively been concerned with two player games.  But what about perfect information games with $n$ players?  

Standard combinatorial game theory says that the winner is determined by who moved last, which makes the idea of an $n$ player game rather complicated, since determining the winner is not easy.  For example consider the following:

$G$ is a three player combinatorial game, where the last player to move is the winner.  The three players are called Left, Right and Middle.  Left and Right both have a single move, but Middle does not and it is Middle's turn to play.  Who wins this game?  If there were just two players, Left and Right, and Left has a move but Right does not, then we'd say Right loses playing first, or Left wins playing second.  But in the case above both Left and Right can move, so while Middle loses playing first, we must pick a winner from Left and Right.

This kind of confusion can readily be avoided with a scoring play game, since we no longer care who moves last, simply who has the highest number of points.

With two players we defined $G^S$ to be the difference between Left's points and Right's points.  This made the process of defining a scoring play game somewhat easier.  Of course with $n$ players we cannot do that, so instead we would define $G^S$ to be a co-ordinate in $n$ dimensional space.

In other words if $G$ is an $n$ player scoring play game, then $G^1,\dots,G^n$ are the set of moves that player $1,\dots,n$ can make and $G^S=(x_1,\dots,x_n)$, where $x_i$ is the number of points that player $i$ currently has.  Then clearly the winner is simply the player for which $x_i$ is the largest, and in the case where they are all equal, the game is a tie.

Similarly instead of a game tree, we could define an $n$-ary game tree, where each node is labelled by a $G^S$.  In this way we can extend combinatorial game theory to whole new areas of research, with potential applications in economics, and cross-overs with classical game theory.

\section{Conclusion}

In this thesis I have introduced a new theory for scoring play games.  I have looked at the basic structure of these games under the disjunctive sum, conjunctive sum, selective sum and sequential join.  I have also looked at impartial scoring play games and shown that there is an equivalent Sprague-Grundy theory that can be used to analyse many of these games.

I have also taken a look at some real life games such as Go, and sowing games and demonstrated how the theory presented in this thesis can be used to give a greater understanding of these games.

I feel that this work is substantial and very innovative, and I think that there are many new areas and ideas that can be explored further.  I think this will lead to some interesting new theories, with some potentially bold and far reaching applications.

\backmatter


\begin{thebibliography}{}

\end{thebibliography}


\begin{thebibliography}{99}

\bibitem{LIP}{M. Albert, R. Nowakowski, D. Wolfe, \emph{Lessons in Play},
A.K. Peters (2007).}

\bibitem{WW}{E. Berlekamp, J. Conway, R. Guy, \emph{Winning Ways for your
Mathematical Plays}, Volumes 1-4, A.K. Peters (2002).}

\bibitem{MG}{E. Berlekamp, D. Wolfe, \emph{Mathematical Go: Chilling Gets
the Last Point}, A.K. Peters (1994).}

\bibitem{ONAG}{J. Conway, \emph{On Numbers and Games}, A.K. Peters (2000).}

\bibitem{JE}{J. Erickson, \emph{Sowing Games}, Games of No Chance, MRSI Publications \textbf{29} (1996).}

\bibitem{FT}{A. S. Fraenkel, U. Tassa, \emph{Strategy for a class of games with dynamic ties},
Comput. Math. Appl. \textbf{1} (1975) 237--254.}

\bibitem{FT2}{A. S. Fraenkel, U. Tassa, \emph{Strategies for compounds of partizan games},
Math. Proc. Cambridge Philos. Soc. \textbf{92} (1982) 193--204.}

\bibitem{GJ}{M. R. Garey, D. S. Johnson, \emph{COMPUTERS AND INTRACTABILITY
A Guide to the Theory of NP-Completeness}, Freeman (1979).}

\bibitem{PG}{P. Grundy, \emph{Mathematics and Games}, Eureka \textbf{2} (1939) 6--8; reprinted ibid.\textbf{27} 1964 9--11.}

\bibitem{PS}{P. Grundy, C. Smith, \emph{Disjunctive games with the last player losing} Proc. Camb. Philos Soc. \textbf{52} (1956) 443--458.}

\bibitem{H}{P.G. Hinman, \emph{Finite termination games with tie}, Israel J. Math. \textbf{12} (1972) 17-22.}

\bibitem{LY}{S. Li,  R. Yen,
\emph{Sums of Zuchswang games.},
J. Combinatorial Theory Ser. \textbf{A 21} (1976) 52--67.}

\bibitem{LY2}{S. Li, R. Yen,
\emph{Generalized impartial games},
Internat. J. Game Theory \textbf{3} (1974) 169--184.}

\bibitem{GAMO}{G. A. Mesdal, P. Ottaway, \emph{Simplifcation of
Partizan Games in Mis\`ere Play}, Integers \textbf{7} (2007).}

\bibitem{rjn}{R. Nowakowski, Personal Communication.}

\bibitem{O}{P. Ottaway, Personal Communication.}

\bibitem{P}{T. Plambeck, \emph{Taming the Wild in Impartial Combinatorial
Games}, Integers \textbf{5} (2005).}

\bibitem{PS}{T. Plambeck, A. Siegel, \emph{Mis\`ere Quotients for Impartial
Games} \url{http://www.integraldomain.net/aaron/}.}

\bibitem{RS}{R. Sprague, \emph{Uber Mathematische Kampfspiele}, Tohoku Math. J. \textbf{41} (1935-1936) 438--444; Zbl. \textbf{13} 290.}

\bibitem{FS}{F. Stewart, \emph{Finding Structure in Mis\`ere Games and
Studying Variants of Mis\`ere Hackenbush}, MSc Thesis, Dalhousie
University, Department of Mathematics and Statistics (2006).}

\bibitem{FS2}{F. Stewart, \emph{The Sequential Join of Combinatorial Games}, Integers \textbf{7} (2007).}

\bibitem{SU}{W. Stromquist, D. Ullman, \emph{Sequential Join of a
Combinatorial Game}, Journal of Theoretical Computer Science
\textbf{119} (1993) 311--321.}

\bibitem{AT}{A. Turing, \emph{On Computable Numbers, with an Application to the Entscheidungsproblem}, Proc. London Math. Soc. \textbf{42} 230--265 and \textbf{43} 544--546.}

\bibitem{IV}{I.P. Varvak,
\emph{Games on the sum of graphs},
Cybernetics \textbf{4} (1968) 49--51 (trans. of  Kibernetika  \textbf{4} (1968) 63--66).}


\end{thebibliography}
\end{document}